\documentclass[12pt]{amsart}
\usepackage{amsmath,amssymb,amsfonts,amsthm}
\usepackage[all,cmtip]{xy}
\usepackage{fullpage}
\usepackage[OT2]{fontenc}
\usepackage[T5,T1]{fontenc}
\DeclareSymbolFont{cyrletters}{OT2}{wncyr}{m}{n}
\DeclareMathAlphabet{\mathpzc}{OT1}{pzc}{m}{it}
\usepackage{hyperref}
\usepackage[francais,english]{babel}
\begin{document}
\theoremstyle{plain}
\newtheorem{thm}{Theorem}[section]
\newtheorem{lem}[thm]{Lemma}
\newtheorem{cor}[thm]{Corollary}
\newtheorem{prop}[thm]{Proposition}
\newtheorem{rem}[thm]{Remark}
\newtheorem{defn}[thm]{Definition}
\newtheorem{ex}[thm]{Example}
\newtheorem{ques}[thm]{Question}
\newtheorem{fact}[thm]{Fact}
\newtheorem{conj}[thm]{Conjecture}
\newtheorem{nota}[thm]{Notation}
\numberwithin{equation}{section}
\def\theequation{\thesection.\arabic{equation}}
\newcommand{\mc}{\mathcal}
\newcommand{\mb}{\mathbb}
\newcommand{\surj}{\twoheadrightarrow}
\newcommand{\inj}{\hookrightarrow}
\newcommand{\red}{{\rm red}}
\newcommand{\codim}{{\rm codim}}
\newcommand{\rank}{{\rm rank}}
\newcommand{\Ker}{{\rm Ker \,}}
\newcommand{\Hom}{{\rm Hom}}
\newcommand{\im}{{\rm im \,}}
\newcommand{\Spec}{{\rm Spec \,}}
\newcommand{\Sing}{{\rm Sing}}
\newcommand{\Char}{{\rm char}}
\newcommand{\Tr}{{\rm Tr}}
\newcommand{\CH}{{\rm CH}}
\newcommand{\pr}{{\rm pr}}
\newcommand{\gr}{{\rm Gr }}
\newcommand{\Coker}{{\rm Coker \,}}
\newcommand{\id}{{\rm id}}
\newcommand{\sA}{{\mathcal A}}
\newcommand{\sB}{{\mathcal B}}
\newcommand{\sC}{{\mathcal C}}
\newcommand{\sD}{{\mathcal D}}
\newcommand{\sE}{{\mathcal E}}
\newcommand{\sF}{{\mathcal F}}
\newcommand{\sG}{{\mathcal G}}
\newcommand{\sH}{{\mathcal H}}
\newcommand{\sI}{{\mathcal I}}
\newcommand{\sJ}{{\mathcal J}}
\newcommand{\sK}{{\mathcal K}}
\newcommand{\sL}{{\mathcal L}}
\newcommand{\sM}{{\mathcal M}}
\newcommand{\sN}{{\mathcal N}}
\newcommand{\sO}{{\mathcal O}}
\newcommand{\sP}{{\mathcal P}}
\newcommand{\sQ}{{\mathcal Q}}
\newcommand{\sR}{{\mathcal R}}
\newcommand{\sS}{{\mathcal S}}
\newcommand{\sT}{{\mathcal T}}
\newcommand{\sU}{{\mathcal U}}
\newcommand{\sV}{{\mathcal V}}
\newcommand{\sW}{{\mathcal W}}
\newcommand{\sX}{{\mathcal X}}
\newcommand{\sY}{{\mathcal Y}}
\newcommand{\sZ}{{\mathcal Z}}
\newcommand{\A}{{\mathbb A}}
\newcommand{\B}{{\mathbb B}}
\newcommand{\C}{{\mathbb C}}
\newcommand{\D}{{\mathbb D}}
\newcommand{\E}{{\mathbb E}}
\newcommand{\F}{{\mathbb F}}
\newcommand{\G}{{\mathbb G}}
\renewcommand{\H}{{\mathbb H}}
\newcommand{\I}{{\mathbb I}}
\newcommand{\J}{{\mathbb J}}
\newcommand{\M}{{\mathbb M}}
\newcommand{\N}{{\mathbb N}}
\renewcommand{\P}{{\mathbb P}}
\newcommand{\Q}{{\mathbb Q}}
\newcommand{\R}{{\mathbb R}}
\newcommand{\T}{{\mathbb T}}
\newcommand{\V}{{\mathbb V}}
\newcommand{\W}{{\mathbb W}}
\newcommand{\X}{{\mathbb X}}
\newcommand{\Y}{{\mathbb Y}}
\newcommand{\Z}{{\mathbb Z}}
\newcommand{\Image}{\rm Im \,}
\title[Motivic Stable Homotopy]{E-Motives and motivic stable homotopy}
\author{{Nguyen Le Dang Thi}}
\email{nguyen.le.dang.thi@gmail.com}
\date{19. 05. 2015}          
\begin{abstract}
We introduce in this work the notion of the category of pure $\mathbf{E}$-Motives, where $\mathbf{E}$ is a motivic strict ring spectrum and construct twisted $\mathbf{E}$-cohomology by using six functors formalism of J. Ayoub. In particular, we construct the category of pure Chow-Witt motives $CHW(k)_{\Q}$ over a field $k$ and show that this category admits a fully faithful embedding into the geometric stable $\A^1$-derived category $D_{\A^1,gm}(k)_{\Q}$.  
\end{abstract}
\maketitle
\tableofcontents
\section{Introduction}
One of the main motivations for this work is the embedding theorem of Voevodsky \cite{Voe00}, which asserts that there is a fully faithful embedding of the category of Grothendieck-Chow pure motives $\underline{Chow}(k)$ into the category of geometric motives $\bold{DM}_{gm}(k)$, hence also into the category of motives $\mathbf{DM}^-_{Nis}(k)$
\[\underline{Chow}(k)^{op} \rightarrow \mathbf{DM}_{gm}(k),\]
if $k$ is a perfect field, which admits resolution of singularities (see e.g. \cite[Prop. 20.1 and Rem. 20.2]{MVW06}, the assumption on resolution of singularities can be removed by using Poincar\'e duality). 
In this note, we construct a category $CHW(k)_{\Q}$, which we call the category of pure Chow-Witt motives over a field $k$ and show that $CHW(k)_{\Q}$ admits a fully faithful embedding into the geometric $\P^1$-stable $\A^1$-derived category $D_{\A^1,gm}(k)_{\Q}$ rationally. Our work can be viewed as an $\A^1$-version for Voevodsky's embedding theorem. The advantage here is that by using duality formalism for $\P^1$-stable $\A^1$-derived category $D_{\A^1}(k)$ (see \cite[App. A]{Hu05} for stable $\A^1$-homotopy categories) and the six operations formalism of J. Ayoub \cite{Ay08}, we do not have to assume the resolution of singularities. However, unlike in motivic setting, one of the main problems here is that we don't have cancellation theorem for the effective $\A^1$-derived category in general, see \cite[Rem. 3.2.4]{AH11}, that is the reason why we can prove the embedding result only for $\Q$-coefficient. F. Morel conjectured in general that (see \cite{Mor04}): 
\begin{conj}\label{conjst}\cite{Mor04} 
Let $S$ be a regular Noetherian scheme of finite Krull dimension. One has a direct decomposition in the rationally motivic stable homotopy category $\mathbf{StHo}_{\A^1,\P^1}(S)$:  
\[[S^i,\G_m^{\wedge j}]_{\P^1}\otimes \Q = H^{j-i}_B(S,\Q(j)) \oplus H^{-i}_{Nis}(S,\underline{\bold{W}}\otimes \Q),\]
where $H^*_B(-,\Q(*))$ denotes the Beilinson motivic cohomology, $\underline{\bold{W}}$ is the unramified Witt sheaf and $[-,-]_{\P^1} \otimes \Q$ denotes $\Hom_{\mathbf{StHo}_{\A^1,\P^1}(S)}(-,-)_{\Q}$. 
\end{conj}
Over a general base scheme $S$ one can split the rational motivic sphere spectrum $\mathbf{1}_{\Q} = \mathbf{1}_{\Q+} \vee \mathbf{1}_{\Q-}$. The identifcation of the plus part $\mathbf{1}_{\Q+} = H_B$ has been done in \cite[Thm 16.2.13]{CD10} over any Noetherian scheme of finite Krull dimension $S$. The minus part $\mathbf{1}_{\Q-} = H\underline{\mathbf{W}}_{*\Q}$, where $H\underline{\mathbf{W}}_{*\Q}$ denotes the Eilenberg-Maclane spectrum associated to the rational Witt homotopy module $\underline{\mathbf{W}}_{*\Q}$, is given in the work of A. Ananyeskiy, M. Levine and I. Panin (\cite[\S3, Thm. 5]{ALP15}) over fields $S = \Spec k$. In general, the conjecture \ref{conjst} over a regular Noetherian scheme $S$ of finite Krull dimension is still widely open, as far as I know. On the other hand, our interest started originally from the study of the existence of $0$-cycles of degree one on algebraic varieties. More precisely, H\'el\`ene Esnault asked (cf. \cite{Lev10}): Given a smooth projective variety $X$ over a field $k$, such that $X$ has a zero cycle of degree one. Are there "motivic" explanations which give the (non)-existence of a $k$-rational point? In \cite{AH11}, A. Asok and C. Haesemeyer show that the existence of zero cycles of degree one over an infinite perfect field of $char(k) \neq 2$ is equivalent to the assertion that the structure map $\mathbf{H}^{st\A^1}_0(X) \rightarrow \mathbf{H}^{st \A^1}_0(\Spec k)$ is a split epimorphism, where $\mathbf{H}^{st\A^1}_i(X)$ denotes the $\P^1$-stable $\A^1$-homology sheaves, while in an earlier work \cite{AH11a} they also showed that the existence of a $k$-rational point over an arbitrary field $k$ is equivalent to the condition that the structure map $\bold{H}_0^{\A^1}(X) \rightarrow \bold{H}_0^{\A^1}(\Spec k)$ is split surjective. So roughly speaking, the obstruction to the lifting of a zero cycle of degree one to a rational point arises by passing from $S^1$-spectra to $\P^1$-spectra. As remarked by M. Levine, it is not to expect that the category of Chow-Witt motives $CHW(k)$ contains any information about the existence of rational points. Now we state our main theorem in this work:               
\begin{thm}\label{mainthm}
Let $k$ be a field. There exists a category of pure Chow-Witt motives $CHW(k)_{\Q}$, which admits a fully faithful embedding 
\[CHW(k)_{\Q} \rightarrow D_{\A^1,gm}(k)_{\Q}.\]
\end{thm} 
In fact, one of the main steps in the work of \cite{AH11} is to exhibit a natural isomorphism $H_0^{st\A^1}(X)(L) \rightarrow \widetilde{\CH}_0(X_L)$ for any separable, finitely generated field extension $L/k$. So one may relate this step to our work as evaluating at a generic point, but much weaker than expected, since we can only prove the result for $\Q$-coefficient. Now our paper is organized as follows: we will review shortly $\A^1$-homotopy theory in section \S 2. Section \S 3 is devoted for $\A^1$-derived categories, in fact we will define the geometric $\P^1$-stable $\A^1$-derived category $D_{\A^1,gm}(k)$ over a field $k$ in \ref{defgmA1} at the end of \S 3. In fact, this is the subcategory of compact objects $D_{\A^1,c}(k)$ of $D_{\A^1}(k)$ (see \cite[Ex. 5.3.43]{CD10}). In these \S 2 and \S 3 we simply steal everything which is needed from the presentation of \cite{AH11}. For a complete treatment we strongly recommend the reader to \cite{Ay08}, \cite{CD10} and \cite{Mor12}. In section \S 4 we introduce the notion of pure $\mathbf{E}$-motives, where $\mathbf{E}$ is a motivic strict ring spectrum and relate several categories of $\mathbf{E}$-correspondences with each other via the twisted $\mathbf{E}$-cohomology. The twisted $\mathbf{E}$-cohomology appears since we will not assume the motivic ring spectrum $\mathbf{E}$ to be orientable. In topology, if $\mathbf{E}$ is a multiplicative cohomology theory and $V$ is an $\mathbf{E}$-orientable vector bundle of rank $r$, then one has a Thom-Dold isomorphism
\[\mathbf{E}^*(X) \stackrel{\cong}{\longrightarrow} \widetilde{\mathbf{E}}^{*+r}(Th(V)),\]
where $Th(V)$ is the Thom space of $V$ and the right hand side is the reduced cohomology. If $\mathbf{E}$ is a ring spectrum, then one can intepret this isomorphism as following: Via the Thom diagonal 
\[Th(V) \rightarrow Th(V) \wedge X_+,\]
which is induced by the diagonal 
\[X_+ \rightarrow X_+ \wedge X_+\]
one can express $Th(V)$ as a comodule over $X_+$ and the comodule map is the natural map 
\[X_+ \rightarrow Th(V).\] 
The geometric Thom isomorphism is the homotopy equivalence 
\[\mathbf{E} \wedge Th(V) \rightarrow \mathbf{E} \wedge Th(V) \wedge X_+ \rightarrow \mathbf{E} \wedge \Sigma^n\mathbf{E} \wedge X_+ \stackrel{\mu_{\mathbf{E}} \wedge \id}{\longrightarrow} \mathbf{E} \wedge \Sigma^n X_+.\]
The composition is an $\mathbf{E}$-module map, hence one may take function spectrum 
\[F_{\mathbf{E}}(\mathbf{E} \wedge \Sigma^nX_+ ,\mathbf{E}) \simeq F(\Sigma^nX_+,\mathbf{E}) \stackrel{\simeq}{\longrightarrow} F_{\mathbf{E}}(\mathbf{E} \wedge Th(V),\mathbf{E}) \simeq F(Th(V),\mathbf{E}),\]
which induces the Thom isomorphism on $\mathbf{E}$-cohomology. In algebraic geometry one has a similar result. For an oriented motivic ring spectrum $\mathbf{E} \in SH(S)$, where $S$ is a regular base, one has (\cite[Thm. 2.12]{NSO09})
\[\mathbf{E}^{*,*}(X) \stackrel{\cong}{\longrightarrow} \mathbf{E}^{*+2r,*+r}(Th(V)),\]
where $V$ is a vector bundle of rank $r$ on a smooth $S$-scheme $X$. The key point is that since $\mathbf{E}$ is oriented one can define the first Chern class and then prove the projective bundle theorem \cite[Thm. 2.11]{NSO09}. The situation becomes much more difficult, even in topology, if $\mathbf{E}$ is not necessary oriented. One has to introduce twisted cohomology. Again in topology, by Atiyah duality one has a commutative diagram in the $(\infty,1)$-category $\mathbb{S}Mod$: 
\[\xymatrix{ & Th(-T_X) \ar[d]^{\simeq} \\ \mathbb{S} \ar[ru]^{PT} \ar[r] & X^{\vee} }\]   
where $PT : \mathbb{S} \rightarrow Th(-T_X)$ is the Pontryagin-Thom collapse map. Let $\mathbf{E}$ be an $E_{\infty}$-ring spectrum. By taking $-\wedge_{\mathbb{S}} \mathbf{E}$ one obtains a map in the $(\infty,1)$-category $\mathbf{E}Mod$
\[\mathbf{E} \rightarrow X^{\vee} \wedge_{\mathbb{S}} \mathbf{E}.\] 
Taking function spectrum we have the (twisted) Umkehr map 
\[F_{\mathbf{E}}(Th(-T_X)\wedge_{\mathbb{S}} \mathbf{E},\mathbf{E}) \simeq F(Th(-T_X),\mathbf{E}) \rightarrow \mathbf{E}.\]
If $\mathbf{E}$ is non-oriented, there is no geometric Thom isomorphism. However, $F_{\mathbf{E}}(Th(-T_X) \wedge_{\mathbb{S}}\mathbf{E},\mathbf{E})$ will give the twisted cohomology. This is the motivation from topology for us, since in algebraic geometry we also have the Atiyah-Spanier-Whitehead duality, but I do not know any $\infty$-categorical approach to twisted cohomology like the one in topology \cite{ABGHR14}. So I introduce in section \S 4 the twisted $\mathbf{E}$-cohomology rather through the guide of the six functors formalism of J. Ayoub. The reader may recognize that the notion of $\mathbf{E}$-correspondences is similar to the construction of Jack Morava in topology. While it is very simple to define the category of $\mathbf{E}$-correspondences $Corr_{\mathbf{E}}(k)$, it is quite difficult to construct the cateogry $\widetilde{Corr}_{\mathbf{E}}(k)$ via twisted $\mathbf{E}$-cohomology. This category exists only up to a number of natural $2$-isomorphisms. This phenomenon reflexes the fact that we rely on six functors formalism, where Thom transformations are only $2$-isomorphic to each other. The composition in $\widetilde{Corr}_{\mathbf{E}}(k)$ is associative only up to a natural isomorphism induced by a natural $2$-isomorphism. In \S 5 we give the proof of the main theorem. In the appendix we give a minimal list of well-known facts and definitions of model categories. We fix now some notations throughout this work. For a pair of adjoint functors $F: \sA \rightarrow \sB$ and $G: \sB \rightarrow \sA$, we will adopt the notation in \cite{CD10} 
\[F: \sA \rightleftarrows \sB : G,\]
where $F$ is left adjoint to $G$ and $G$ is right adjoint to $F$. Sometime we will write
\[\varepsilon_{(F,G)} : FG \rightarrow \id, \quad \eta_{(F,G)} : \id \rightarrow GF\]
for the counit and unit of the adjunction repsectively. For every morphism $f: FY \rightarrow X$ in $Mor(\sB)$, there is a unique morphism $g: Y  \rightarrow GX$ in $Mor(\sA)$ such that the following diagram commutes:
\[\xymatrix{FY \ar[d]_{F(g)} \ar[rrd]^f \\ FG(X) \ar[rr]_{\varepsilon_{(F,G)}(X)} && X }\]
For every morphism $g : Y \rightarrow GX$ in $Mor(\sA)$, there is a unique morphism $f: FY \rightarrow X$ in $Mor(\sB)$, such that the following diagram commutes:
\[\xymatrix{Y \ar[rr]^{\eta_{(F,G)}(Y)} \ar[rrd]_g && GF(Y) \ar[d]^{G(f)} \\ && GX }\]
In a symmetric monoidal category $(\sC,\wedge,\mathbf{1})$, an object $A$ is called strongly dualizable if there exists an object $A^{\vee}$ and morphisms
\[coev_A: \mathbf{1} \rightarrow A \wedge A^{\vee}, \quad ev_A:  A^{\vee} \wedge A \rightarrow \mathbf{1},\]
such that the following compositions
\[A \cong \mathbf{1} \wedge A \stackrel{coev_A \wedge \id}{\longrightarrow} A \wedge A^{\vee} \wedge A \stackrel{\id \wedge ev_A}{\longrightarrow} A \wedge \mathbf{1} \cong A \]
and
\[A^{\vee} \cong A^{\vee} \wedge \mathbf{1} \stackrel{\id \wedge coev_A}{\longrightarrow} A^{\vee} \wedge A \wedge A^{\vee} \stackrel{ev_A \wedge \id}{\longrightarrow} \mathbf{1} \wedge A^{\vee} \cong A^{\vee} \]
are the identities $\id_A$ and $\id_{A^{\vee}}$. The natural isomorphism
\[\alpha: \Hom_{\sC}(-,A) \stackrel{\cong}{\longrightarrow} \Hom_{\sC}(A^{\vee} \wedge - , \mathbf{1})\]
is given by 
\[\alpha(\phi) = ev_{A^{\vee}} \circ (\id_{A^{\vee}} \wedge \phi),\]
and its inverse $\alpha^{-1}$ is given by 
\[\alpha^{-1}(\varphi) = (\id_A \wedge \varphi) \circ (coev_{A^{\vee}} \wedge \id_{-}).\]
Given two smooth $k$-schemes $X,Y \in Sm/k$ and two vector bundles $\sE, \sE'$ over $X$ resp. $Y$, we write $\sE \times \sE'/ X\times Y$ for the external sum over $X \times_k Y$. The $\P^1$- stable homotopy category over a base scheme $S$ will be denoted by $\mathbf{StHo}_{\A^1,\P^1}(S)$ and we write $\mathbf{StHo}_{\A^1,S^1}(S)$ for the $S^1$-stable homotopy category. Sometime when it is clear which category we are talking about, we just abbreviate our $\P^1$-stable homotopy category by $SH(S)$.    
\section{$\A^1$-homotopy category}
\subsection{Unstable $\A^1$-homotopy category}
Let $Sm/k$ denote the category of separated smooth schemes of finite type over a field $k$. We write $Spc/k$ for the category $\Delta^{op}Sh_{Nis}(Sm/k)$ consisting of simplicial Nisnevich sheaves of sets on $Sm/k$. An object in $Spc/k$ is simply called a $k$-space, which is usually denoted by calligraphic letter $\sX$. The Yoneda embedding $Sm/k \rightarrow Spc/k$
is given by sending a smooth scheme $X \in Sm/k$ to the corresponding representable sheaf $\Hom_{Sm/k}(-,X)$ then by taking the associated constant simplicial object, where all face and degeneracy maps are the identity. We will identify $Sm/k$ with its essential image in $Spc/k$. Denote by $Spc_+/k$ the category of pointed $k$-space, whose objects are $(\sX,x)$, where $\sX$ is a $k$-space and $x: \Spec k \rightarrow \sX$ is a distinguished point. One has an adjoint pair 
\[ Spc/k \rightleftarrows Spc_+/k, \]
which means that the functor $Spc/k \rightarrow Spc_+/k$ sending $\sX \rightarrow \sX_+ = \sX \coprod \Spec k$ is left-adjoint to the forgetful functor $Spc_+/k \rightarrow Spc/k$. The category $Spc/k$ can be equipped with the injective local model structure $(C_s,W_s,F_s)$, where cofibrations are monomorphisms, weak equivalences are stalkwise weak equivalences of simplicial sets and fibrations are morphisms with right lifting property wrt. morphisms in $C_s \cap W_s$. Denote by $\bold{Ho}_s^{Nis}(k)$ the resulting unpointed homotopy category as constructed by Joyal-Jardine (cf. \cite[\S 2 Thm.  1.4]{MV01}). We will write $\bold{Ho}_{s,+}^{Nis}(k)$ for the pointed homotopy category. 
\begin{defn}\cite{MV01}\label{def2.1}
\begin{enumerate}
\item A $k$-space $\sZ \in Spc/k$ is called $\A^1$-local if and only for any object $\sX \in Spc/k$, the projection $\sX \times \A^1 \rightarrow \sX$ induces a bijection 
\[\Hom_{\mathbf{Ho}_s^{Nis}(k)}(\sX,\sZ) \stackrel{\simeq}{\rightarrow} \Hom_{\mathbf{Ho}_s^{Nis}(k)}(\sX \times \A^1,\sZ). \]
\item Let $\sX \rightarrow \sY \in Mor(Spc/k)$ be a morphism of $k$-spaces. It is an $\A^1$-weak equivalence if and only for any $\A^1$-local object $\sZ$, the induced map 
\[\Hom_{\mathbf{Ho}_s^{Nis}(k)}(\sY,\sZ) \rightarrow \Hom_{\mathbf{Ho}_s^{Nis}(k)}(\sX,\sZ) \]
is bijective.  
\end{enumerate}
\end{defn}       
In \cite[\S 2 Thm. 3.2]{MV01}, F. Morel and V. Voevodsky proved that $Spc/k$ can be endowed with the $\A^1$-local injective model structure $(C,W_{\A^1},F_{\A^1})$, where cofibrations are monomorphisms, weak equivalences are $\A^1$-weak equivalences. The associated homotopy category obtained from $Spc/k$ by inverting $\A^1$-weak equivalences is denoted by $\mathbf{Ho}_{\A^1}(k) \stackrel{def}{=} Spc/k[W_{\A^1}^{-1}] $. This category is called the unstable $\A^1$-homotopy category of smooth $k$-schemes. Let $\mathbf{Ho}_{s,\A^1-loc}^{Nis}(k) \subset \mathbf{Ho}_s^{Nis}(k)$ be the full subcategory consisting of $\A^1$-local objects. In fact, one has an adjoint pair (cf. \cite{MV01})
\[L_{\A^1}: \mathbf{Ho}_s^{Nis}(k) \rightleftarrows \mathbf{Ho}_{s,\A^1-loc}^{Nis}(k) : i,\]
where $L_{\A^1}$ is the $\A^1$-localization functor sending $\A^1$-weak equivalences to isomorphisms. $L_{\A^1}$ induces thus an equivalence of categories $\mathbf{Ho}_{\A^1}(k) \rightarrow \mathbf{Ho}_{s,\A^1-loc}^{Nis}(k)$. This will imply that if $\sX \in Spc/k$ is any object and $\sY$ is an $\A^1$-local object, then one has a canonical bijection 
\[\Hom_{\mathbf{Ho}_s^{Nis}(k)}(\sX,\sY) \stackrel{\simeq}{\rightarrow} \Hom_{\mathbf{Ho}_{\A^1}(k)}(\sX,\sY).\]
We will write $\mathbf{Ho}_{\A^1,+}(k)$ for the unstable pointed $\A^1$-homotopy category of smooth $k$-schemes. Recall 
\begin{defn}\label{defThomsp} 
Let $X \in Sm/k$ and $E$ be a vector bundle over $X$. The Thom space of $E$ is the pointed sheaf 
\[Th(E/X) = E/E-s_0(X),\]
where $s_0: X \rightarrow E$ is the zero section of $E$.
\end{defn}
Let $T \in Spc_+/k$ be the quotient sheaf $\A^1/(\A^1- \{0\})$ pointed by the image of $\A^1-\{0\}$. Then $T \cong S^1_t \wedge S^1_s$ in $\bold{Ho}_{\A^1,+}(k)$ (\cite[Lem. 2. 15]{MV01}). For a pointed space $\sX \in Spc_+/k$, we denote by $\Sigma_T(\sX,x) = T \wedge (\sX,x)$. Remark that $\P^n/\P^{n-1} \cong T^n \stackrel{def}{=} T^{\wedge n}$ is an $\A^1$-equivalence. In particular, we have $(\P^1,*) \cong T$ (\cite[Cor. 2.18]{MV01}). Recall 
\begin{prop}\cite[\S 3 Prop. 2. 17]{MV01}\label{propthom}
Let $X, Y \in Sm/k$ and $E, E'$ be vector bundles on $X$ and $Y$ respectively. One has 
\begin{enumerate} 
\item There is a canonical isomorphism of pointed sheaves 
\[Th(E \times E' / X \times Y) = Th(E/X) \wedge Th(E'/Y).\]
\item There is a canonical isomorphism of pointed sheaves 
\[Th(\sO_X^n) = \Sigma^n_T X_+\]
\item The canonical morphism of pointed sheaves 
\[\P(E \oplus  \sO_X)/\P(E) \rightarrow Th(E)\]
is an $\A^1$-weak equivalence. 
\end{enumerate}
\end{prop}
The following theorem due to Voevodsky will play an essential role for our purpose. However, as pointed out by M. Levine, the identities in $K_0(-)$ are not enough for us to construct maps between twisted $\mathbf{E}$-cohomology. Following a suggestion by M. Levine, we will refine this result of Voevodsky later (see \ref{refinement1}).  
\begin{thm}\cite[Thm. 2.11]{Voe03}\label{ThmVoeThom}
Let $X \in SmProj/k$ a smooth projective variety of pure dimension $d_X$ over a field $k$. There exists an integer $n_X$ and a vector bundle $V_X$ over $X$ of rank $n_X$, such that  
\[ V_X \oplus T_X = \sO_X^{n_X+d_X} \in K_0(X),\]
where $T_X$ denotes the tangent bundle of $X$. Moreover, there exists a morphism $T^{\wedge n_X + d_X} \rightarrow Th(V_X)$ in $\bold{Ho}_{\A^1,+}(k)$, such that the induced map $H^{2d_X}_{\sM}(X,\Z(d_X)) \rightarrow \Z$ coincides with the degree map $\deg: \CH_0(X) \rightarrow \Z$, where $T = S^1_s \wedge \G_m$.
\end{thm}
\begin{rem}\rm{ One can always add a trivial bundle to $V_X$ in Voevodsky's theorem \ref{ThmVoeThom} to increase $n_X$ appropriately. 
}
\end{rem}
\subsection{Stable $\A^1$-homotopy category} 
Let $\mathbf{Spect}^{\Sigma}(k)$ be the category of symmetric spectra in $k$-spaces, which can be viewed as category of Nisnevich sheaves of symmetric spectra. By applying the construction in \cite[Def. 4.4.40, Cor. 4.4.42, Prop. 4.4.62]{Ay08}, $\mathbf{Spect}^{\Sigma}(k)$ has the structure of a monoidal model category. Let $\mathbf{StHo}_{S^1}(k)$ be the resulting homotopy category. The stable $\A^1$-homotopy category of $S^1$-spectra $\mathbf{StHo}_{\A^1,S^1}(k)$ is obtained from $\mathbf{StHo}_{S^1}(k)$ by Bousfield localization. Equivalently, the category $\mathbf{Spect}^{\Sigma}(k)$ can be equipped with an $\A^1$-local model structure (cf. \cite[Def. 4.5.12]{Ay08}). The homotopy category of this $\A^1$-local model structure is $\mathbf{StHo}_{\A^1,S^1}(k)$, which is also known to be equivalent to the category $\mathbf{StHo}^{S^1}_{\A^1-loc}(Sm/k)$ constructed by F. Morel in \cite[Def. 4.1.1]{Mor05}. The $\A^1$-local symmetric sphere spectrum is defined by taking the functor 
\[\underline{n} \mapsto L_{\A^1}(S_s^{1 \wedge n})\]
with an action of symmetric groups, where $L_{\A^1}$ denotes the $\A^1$-localization functor. For a pointed space $(\sX,x)$, its $\A^1$-local symmetric suspension spectrum is defined as the symmetric sequence 
\[\underline{n} \mapsto L_{\A^1}(S_s^{1 \wedge n} \wedge \sX)\]
together with symmetric groups actions. Let $\sE$ be an $\A^1$-local symmetric
spectrum in $Spc/k$. One defines (\cite[Def. 2.1.11]{AH11}) the $i$-th
$S^1$-stable $\A^1$-homotopy sheaf $\pi_i^{st\A^1,S^1}(\sE)$ of $\sE$ as the
Nisnevich sheaf on $Sm/k$ associated to the presheaf 
\[U \mapsto \Hom_{\mathbf{StHo}_{\A^1,S^1}(k)}(S^{1 \wedge i}_s \wedge \Sigma^{\infty}_s U_+,\sE).\]
Now we consider the symmetric $T$-spectra or $\P^1$-spectra (\cite{Jar00}). $\P^1$ is pointed
with $\infty$ and $\P^{1 \wedge n}$ has a natural action of $\Sigma_n$ by
permutation of the factors, so the association $\underline{n} \mapsto \P^{1
\wedge n}$ is a symmetric sequence. A symmetric $\P^1$-spectrum is a symmetric
sequence with a module structure over the sphere spectrum $\mathbb{S}^0$. Denote
by $\mathbf{Spect}_{\P^1}^{\Sigma}(k)$ the full subcategory of the category of
symmetric sequence in $k$-spaces $\mathbf{Fun}(\mathpzc{Sym},Spc_+/k)$
consisting of
symmetric $\P^1$-spectra, which also has a model structure \cite[Def.
4.5.21]{Ay08}. Here we denote by $\mathpzc{Sym}$ the groupoid, whose objects
are $\underline{n}$ and morphisms are given by bijections. Let
$\mathbf{StHo}_{\A^1,\P^1}(k)$ be the resulting homotopy
category, which is called $\P^1$-stable $\A^1$-homotopy category. For a pointed
space $(\sX,x)$, we will write $\Sigma^{\infty}_{\P^1}(\sX,x)$ for the
suspension symmetric $\P^1$-spectrum, i.e., it is given by the functor
$\underline{n} \mapsto \P^{1 \wedge n} \wedge \sX$ equipped with an action of
symmetric group by permuting the first $n$-factors. Let $\mathbf{S}^i$ be the
suspension symmetric $\P^1$-spectrum of $S^i_s$. If $\sE$ is a symmetric
$\P^1$-spectrum, then the $i$-th $\P^1$-stable $\A^1$-homotopy sheaf
$\pi_i^{st\A^1,\P^1}(\sE)$ is defined as the Nisnevich sheaf on $Sm/k$
associated to the presheaf (cf. \cite[Def. 2.1.14]{AH11}) 
\[U \mapsto \Hom_{\mathbf{StHo}_{\A^1,\P^1}(k)}(\bold{S}^i \wedge
\Sigma^{\infty}_{\P^1}U_+,\sE).\] 
\begin{thm}\label{thmMorS1}\cite[Thm. 6.1.8 and Cor. 6.2.9]{Mor05}
 Let $\sE$ be an $\A^1$-local symmetric $S^1$-spectrum. The homotopy
sheaves $\pi_i^{st\A^1,S^1}(\sE)$ are strictly $\A^1$-invariant.    
\end{thm}
One has a canonical isomorphism \cite[Prop. 2.1.16]{AH11}
\begin{multline*}
colim_n\Hom_{\mathbf{StHo}_{\A^1,S^1}(k)}(\Sigma^{\infty}_s \G_m^{\wedge n}
\wedge \Sigma^{\infty}_s(U_+), \Sigma^{\infty}_s \G_m^{\wedge n} \wedge
\Sigma^{\infty}_s(\sX,x)) \stackrel{\cong}{\rightarrow} \\
\Hom_{\mathbf{StHo}_{\A^1,\P^1}(k)}(\Sigma^{\infty}_{\P^1}(U_+),\Sigma^{\infty}
_{\P^1}(\sX,x)).
\end{multline*}
So one may view that $\bold{StHo}_{\A^1,\P^1}(k)$ is obtained from
$\mathbf{StHo}_{\A^1,S^1}(k)$ by formally inverting the $\A^1$-localized
suspension spectrum of $\G_m$. So from \ref{thmMorS1}, we see that for a
pointed $k$-space $(\sX,x)$, the homotopy sheaves $\pi_i^{st\A^1,\P^1}(\sX)$
are also strictly $\A^1$-invariant. By the computation of F. Morel (\cite{Mor04}, \cite{Mor12}), one can
identify the Milnor-Witt $K$-theory sheaves with stable homotopy sheaves of spheres 
\[\mathbf{K}^{MW}_n \cong
\pi_0^{st\A^1,\P^1}(\Sigma^{\infty}_{\P^1}(\G_m^{\wedge n})).\]
This identification allows us to conclude that $\bold{K}^{MW}_n$ are strictly
$\A^1$-invariant sheaves.  
\section{$\A^1$-homological algebra}
\subsection{Effective $\A^1$-derived category}
Let $Ch_{-}(\sA b_k)$ be the category of chain complexes over the category $\sA
b_k$ of abelian Nisnevich sheaves. Denote by $Ch_{\geq 0}(\sA b_k)$ the
category of chain complexes of abelian Nisnevich sheaves on $Sm/k$, whose homoglocial
degree $\geq 0$. The sheaf-theoretical Dold-Kan correspondence 
\[N: \Delta^{op}\sA b_k \rightleftarrows Ch_{\geq 0}(\sA b_k) : K,\]
where $\Delta^{op} \sA b_k$ is the cateogry of simplicial abelian Nisnevich
sheaves, gives us via the inclusion functor $Ch_{\geq 0}(\sA b_k) \inj
Ch_{-}(\sA b_k)$, a functor 
\[\Delta^{op}(\sA b_k) \rightarrow Ch_{-}(\sA b_k).\]
By applying this functor on the Eilenberg-Maclane spectrum $H\Z$, we obtain a
ring spectrum $\widetilde{H\Z}$ in $\mathbf{Fun}(\mathpzc{Sym},Ch_{-}(\sA b_k))$.
Let $\mathbf{Spect}^{\Sigma}(Ch_{-}(\sA b_k))$ be the full subcategory of the
category $\mathbf{Fun}(\mathpzc{Sym},Ch_{-}(\sA b_k))$ consisting of modules over
$\widetilde{H\Z}$. On the other hand, by composing with the free abelian group
functor 
\[\Z (-): Spc/k \rightarrow \Delta^{op}(\sA b_k),\]
one obtains a functor 
\[\mathbf{Fun}(\mathpzc{Sym},Spc_+/k) \rightarrow
\mathbf{Fun}(\mathpzc{Sym},Ch_{-}(\sA b_k)),\]
which sends the sphere symmetric sequence to $\widetilde{H\Z}$. This induces
then a functor between categories of symmetric spectra 
\[\mathbf{Spect}^{\Sigma}(Spc/k) \rightarrow \mathbf{Spect}^{\Sigma}(Ch_{-}(\sA
b_k)).\] 
In fact, by \cite[Thm. 9.3]{Hov01}, this induces a Quillen functor, which one refers as Hurewicz functor 
\[\mathfrak{H}^{ab}: \mathbf{StHo}_{S^1}(k) \rightarrow D_{-}(\sA b_k).\]
Now the effective $\A^1$-derived category $D^{eff}_{\A^1}(k)$ is constructed by applying $\A^1$-localization on the category $\mathbf{Spect}^{\Sigma}(Ch_{-}(\sA b_k))$. By the work of Cisinski and D\'eglise (cf. \cite[\S 5]{CD10}), this category is equivalent to the $\A^1$-derived category constructed by F. Morel in \cite{Mor12}. Let $(\sX,x) \in Spc_+/k$ be a pointed space, and $\Sigma^{\infty}_s(\sX,x)$ its suspension symmetric spectrum. We apply the Hurewicz functor on $\Sigma^{\infty}_s(\sX,x)$ and then $L_{\A^1}^{ab}(-)$, so we may define a functor 
\[\widetilde{C}_*^{\A^1}: \mathbf{StHo}_{S^1}(k) \rightarrow D^{eff}_{\A^1}(k), \quad \Sigma^{\infty}_s(\sX,x) \mapsto L_{\A^1}^{ab}(\mathfrak{H}^{ab}(\Sigma^{\infty}_s(\sX,x))).\]
Here we write $L_{\A^1}^{ab}$ for the $\A^1$-localization functor on chain complexes to distinguish from the $\A^1$-localization $L_{\A^1}$ on spaces.  If $\sX \in Spc/k$ is not pointed, then we write $C^{\A^1}_*(\sX) \stackrel{def}{=} \widetilde{C}^{\A^1}_*(\sX_+)$. Define $\Z[n] = \mathfrak{H}^{ab}(\Sigma^{\infty}_s S^{n}_s)$. 
\begin{defn}\label{defA1homology}
Let $\sX \in Spc/k$ be a $k$-space. Its $i$-th $\A^1$-homology sheaf is the Nisnevich sheaf $\mathbf{H}^{\A^1}_i(\sX)$ associated to the presheaf 
\[U \mapsto \Hom_{D^{eff}_{\A^1}(k)}(C_*^{\A^1}(U)[i],C_*^{\A^1}(\sX)) \stackrel{def}{=} \Hom_{D^{eff}_{\A^1}(k)}(C_*^{\A^1}(U) \otimes \Z[i],C_*^{\A^1}(\sX)).\]
\end{defn} 
Consider $(\P^1,\infty)$ pointed by $\infty$. According to \cite[Cor. 2.18]{MV01}, we have $\P^1 = S^1_s \wedge \G_m$, so we have an identification $\widetilde{C}_*^{\A^1}(\P^1) = \widetilde{C}^{\A^1}_*(S^1_s \wedge \G_m)$. We define the $\A^1$-Tate complex (called enhanced Tate (motivic) complex by A. Asok and C. Haesemeyer \cite[Def. 2.1.25 and Def. 3.2.1 and Lem. 3.2.2]{AH11}) as 
\[\Z_{\A^1}(n) \stackrel{def}{=} \widetilde{C}_*^{\A^1}(\P^{1 \wedge n})[-2n] = \Z_{\A^1}(1)^{\otimes n}.\]
\begin{defn}\label{defA1unstcoh}
Let $\sX \in Spc/k$ be a $k$-space. The bigraded unstable $\A^1$-cohomology group $H^{p,q}_{\A^1}(\sX,\Z)$ is defined as 
\[H^{p,q}_{\A^1}(\sX,\Z) = \Hom_{D^{eff}_{\A^1}(k)}(C_*^{\A^1}(\sX),\Z_{\A^1}(q)[p]).\]
\end{defn} 
The relationship between unstable $\A^1$-cohomology and Nisnevich hypercohomology with coefficient $\Z_{\A^1}(n)$ is given by the following 
\begin{prop}\label{propunstcoh}\cite[Prop. 3.2.5]{AH11} 
Let $k$ be a field and $\sX \in Spc/k$ be a $k$-space. One has 
\begin{enumerate}
\item For any $p,q$, there is a canonical isomorphism 
\[\H^p_{Nis}(\sX,\Z_{\A^1}(q)) \stackrel{\simeq}{\rightarrow} H^{p,q}_{\A^1}(\sX,\Z).\]
\item The cohomology sheaves $\underline{H}^p(\Z_{\A^1}(q)) = 0$, if $p > q$. 
\item There is a canonical isomorphism $\underline{H}^p(\Z_{\A^1}(p)) \cong \bold{K}^{MW}_p$, for all $p >0$. 
\end{enumerate}
\end{prop} 
\begin{rem}\label{remHI}{\rm
By construction the complex $\Z_{\A^1}(n)$ is $\A^1$-local, hence by definition (cf. \cite[Def. 5.17]{Mor12}) one has immediately that the sheaves $\underline{H}^p(\Z_{\A^1}(q))$ are strictly $\A^1$-invariant. 
}
\end{rem}
\subsection{$\P^1$-stable $\A^1$-derived category}
Having defined an $\A^1$-Tate complex, the way that we stabilize the category $D^{eff}_{\A^1}(k)$ is to invert formally the $\A^1$-Tate complex to obtain the $\P^1$-stable $\A^1$-derived category $D_{\A^1}(k)$. This can be done by following the construction detailed in \cite[\S 5]{CD10}. As before, we take $D_{\A^1}(k)$ as the resulting homotopy category of the model category $\bold{Spect}^{\Sigma}_{\P^1}(Ch_{-}(\sA b_k))$ consisting of modules over the $\A^1$-localization of the normalized chain complex of the free abelian group on the sphere symmetric $\P^1$-spectrum. For a pointed space $(\sX,x) \in Spc_+/k$, the stable $\A^1$-complex $\widetilde{C}_*^{st\A^1}(\sX)$ of $(\sX,x)$ is defined as $L_{\A^1}^{ab}(N\Z(\Sigma^{\infty}_{\P^1}(\sX,x)))$ and if $\sX \in Spc/k$ is an unpointed $k$-space, then we write $C_*^{st\A^1}(\sX) $ for $\widetilde{C}_*^{st\A^1}(\sX_+)$. The category $D_{\A^1}(k)$ has an unit object, denoted by $\bold{1}_k$, which is the complex $\widetilde{C}^{st\A^1}_*(\bold{S}^0)$. Define $\bold{1}_k[n] = \bold{1}_k \otimes \widetilde{C}_*^{st\A^1}(S^n_s) $ and $\widetilde{C}_*^{st\A^1}(\sX)[n] = \widetilde{C}_*^{st\A^1}(\sX) \otimes \bold{1}_k[n]$ for a $k$-space $(\sX,x) \in Spc_+/k$. 
\begin{defn}\label{defA1sthomology}
Let $\sX \in Spc/k$ be a $k$-space. The $i$-th $\P^1$-stable $\A^1$-homology sheaf $\bold{H}^{st\A^1}_i(\sX)$ is the Nisnevich sheaf associated to the presheaf 
\[U \mapsto \Hom_{D_{\A^1}(k)}(C_*^{st\A^1}(U)[i],C_*^{st\A^1}(\sX)). \]
\end{defn}     
Just like in case of stable $\A^1$-homotopy categories, one has the following result 
\begin{prop}\label{propstabilize}\cite[Prop. 2.1.29]{AH11}
Let $U \in Sm/k$ and $(\sX,x) \in Spc_+/k$. One has a canonical isomorphism 
\begin{multline}\label{eqstabilize} 
colim_n \Hom_{D^{eff}_{\A^1}(k)}(C_*^{\A^1}(U) \otimes \Z_{\A^1}(n) [i], \widetilde{C}_*^{\A^1}(\sX) \otimes \Z_{\A^1}(n)[i]) \stackrel{\cong}{\rightarrow} \\  \Hom_{D_{\A^1}(k)}(C_*^{st\A^1}(U),\widetilde{C}_*^{st\A^1}(\sX)).
\end{multline}
\end{prop}
The Hurewicz formalism induces the following functors, which one still calls Hurewicz functors (or abelianization functors) 
\begin{align*}
\mathbf{StHo}_{\A^1,S^1}(k) \rightarrow D^{eff}_{\A^1}(k), \\ \mathbf{StHo}_{\A^1,\P^1}(k) \rightarrow D_{\A^1}(k),
\end{align*}
which give rise to morphisms of sheaves 
\begin{align*}
\pi_i^{st\A^1,S^1}(\Sigma^{\infty}_s(\sX_+)) \rightarrow \bold{H}_i^{\A^1}(\sX), \\ 
\pi_i^{st\A^1,\P^1}(\Sigma^{\infty}_{\P^1}(\sX_+)) \rightarrow \bold{H}_i^{st\A^1}(\sX). 
\end{align*}
\begin{defn}\label{defnstA1coh}
Let $\sX \in Spc/k$ be a $k$-space. The bigraded $\P^1$-stable $\A^1$-cohomology group $H^{p,q}_{st\A^1}(\sX,\Z)$ is defined as 
\[H^{p,q}_{st\A^1}(\sX,\Z) = \Hom_{D_{\A^1}(k)}(C_*^{st\A^1}(\sX),\Z_{\A^1}(q)[p]).\] 
\end{defn}
The advantage of $\P^1$-stable $\A^1$-derived category $D_{\A^1}(k)$ is that one has duality formalism. In the context of stable $\A^1$-homotopy theory, it was done in \cite[App. A]{Hu05} and also \cite{Rio05}. Firstly we recall that Deligne introduced in \cite[\S 4]{Del87} virtual categories. If $f: X \rightarrow \Spec k$ is a smooth $k$-scheme, the category $\mathcal{V}(X)$ of virtual bundles on $X$ is identified to the fundamental groupoid of $\mathcal{K}(X)$ where $\mathcal{K}$ is some $\A^1$-fibrant genuine model of algebraic $K$-theory. An actual vector bundle $\xi$ defines an object $\xi$ in $\mathcal{V}(X)$ whose isomorphism class corresponds to $[\xi] \in K_0(X)$. A short exact sequence of vector bundles
\[0 \rightarrow \xi' \rightarrow \xi \rightarrow \xi'' \rightarrow 0\]
gives not just an equality $[\xi] = [\xi'] + [\xi'']$ in $K_0(X)$ but also a specific isomorphism $\xi \cong \xi' \oplus \xi''$ of objects in $\mathcal{V}(X)$. By using universal property of $\mathcal{V}(X)$ as a Picard category, one can define an isomorphism (see \cite[\S 4]{Rio10})
\[Th(\xi/X) \cong Th(\xi'/X) \wedge Th(\xi''/X).\]
We haven't yet introduced in this section the $6$ operations formalism of J. Ayoub, however we should mention that the construction of Thom spectrum extends to a functor (cf. \cite[Prop. 4.1.1, Def. 4.2.1]{Rio10} and \cite[Thm. 1.5.18]{Ay08})
\[Th_X: \mathcal{V}(X) \rightarrow \mathbf{StHo}_{\A^1,\P^1}(X) \stackrel{f_\#}{\rightarrow} \mathbf{StHo}_{\A^1, \P^1}(k).\]
We discuss a little bit more about the Thom spectrum of virtual bundles. If $\xi$ is a virtual vector bundle on an affine variety $U$, then there exist an actual vector bundle $\xi'$ on $U$ and an integer $n \geq 0$, such that 
\[\xi \oplus \sO_X^n = \xi'.\]
So one may define 
\[\Sigma^{\infty}_{\P^1}Th(\xi/U) = \Sigma^{\infty}_{\P^1}Th(\xi'/U) \wedge S^{-2n,-n}.\]
If $X$ is a projective variety, one can define an affine variety, which is $\A^1$-weak equivalent to $X$ (see \cite[p. 10]{Hu05}): Consider first of all the projective space $\P^N$. One defines
\[U = \P^N \times \P^N \setminus Proj \, k[x_0,\cdots, x_N,y_0,\cdots,y_N]/(\sum_{i=0}^N x_iy_i = 0),\]
which is an $\A^N$-bundle $pr_1: U \rightarrow \P^N$. If $X$ is a projective variety, one has $i : X \inj \P^N$ and the affine variety $\pi: i^*U \rightarrow X$ is an $\A^1$-weak equivalence, where $\pi$ is the pullback of $pr_1$ along the closed immersion $i$: 
\[\xymatrix{i^*U \ar[d]_{\pi} \ar[r] & U \ar[d]^{pr_1} \\ X \ar@{(->}[r]_i & \P^N}\]  
If $-T_X$ is the virtual normal bundle on $X$ of the diagonal embedding $\Delta_X : X \inj X \times_k X$, which is the virtual tangent bundle, then its Thom spectrum is defined to be the Thom spectrum $Th(\mu/i^*U)$, where $\mu$ is the complement of the pullback of the tangent bundle of $X$ along $\pi$. We state the following result in $D_{\A^1}(k)$, although the proof in case of $\mathbf{StHo}_{\A^1,\P^1}(k)$ is given in \cite[Thm. A1]{Hu05} or see \cite[Thm. 2.2]{Rio05} 
\begin{prop}\label{propdual}\cite[Prop. 3.5.2 and Lem. 3.5.3]{AH11}
Let $X \in SmProj/k$, then $C_*^{st\A^1}(X)$ is a strong dualizable object in $D_{\A^1}(k)$ and its dual is $C_*^{st\A^1}(X)^{\vee} = \widetilde{C}_*^{st\A^1}(Th(-T_X))$. Consequently, one has a canonical isomorphism 
\begin{equation}\label{eqdual}
\Hom_{D_{\A^1}(k)}(\mathbf{1}_k,C_*^{st\A^1}(X)) \stackrel{\cong}{\rightarrow} \Hom_{D_{\A^1}(k)}(C_*^{st\A^1}(X)^{\vee},\mathbf{1}_k). 
\end{equation}
\end{prop}
Fortunately, we will use later duality via $6$ operations formalism of J. Ayoub, which is good enough for our main purpose. We end up this section by a definition: 
\begin{defn}\label{defgmA1}
Let $k$ be a field. One defines the geometric stable $\A^1$-derived category $D_{\A^1,gm}(k)$ over $k$ as the thick subcategory of $D_{\A^1}(k)$ generated by $C_*^{st\A^1}(X)$, where $X \in Sm/k$. 
\end{defn}
\section{$\mathbf{E}$-Motives}
\subsection{$\mathbf{E}$-Correspondences}
Let $k$ be a field and we denote by $SH(k)$ the motivic stable homotopy category. Throughout this section we fix a motivic spectrum $\mathbf{E} \in SH(k)$ together with a multiplication map 
\[\mu_{\mathbf{E}} : \mathbf{E} \wedge^{\mathbb{L}}_S \mathbf{E} \rightarrow \mathbf{E}\]
and a unit map 
\[\varphi_{\mathbf{E}}: S \rightarrow \mathbf{E},\]
such that the fowlling diagrams commute 
\[\xymatrix{\mathbf{E} \ar[rr]^{\id \wedge \varphi_{\mathbf{E}}} \ar@{=}[drr] && \mathbf{E} \wedge^{\mathbb{L}}_S \mathbf{E} \ar[d]^{\mu_{\mathbf{E}}} && \ar[ll]_{\varphi_{\mathbf{E}} \wedge \id} \mathbf{E} \ar@{=}[lld] \\ && \mathbf{E}  }\]
\[\xymatrix{ \mathbf{E} \wedge^{\mathbb{L}}_S \mathbf{E} \wedge^{\mathbb{L}}_S \mathbf{E} \ar[d]_{\mu_{\mathbf{E}} \wedge \id} \ar[rr]^{\id \wedge \mu_{\mathbf{E}}} && \mathbf{E} \wedge^{\mathbb{L}}_S \mathbf{E} \ar[d]^{\mu_{\mathbf{E}}} \\ \mathbf{E} \wedge^{\mathbb{L}}_S \mathbf{E} \ar[rr]^{\mu_{\mathbf{E}}} && \mathbf{E}  }\]
Such a triple $(\mathbf{E},\mu_{\mathbf{E}},\varphi_{\mathbf{E}})$ is called a motivic ring spectrum. 
\begin{prop}\label{propEcor}
Let $X,Y,Z,W \in SmProj(k)$. Let $\alpha \in SH(k)[\Sigma^{\infty}_{T,+}X,\Sigma^{\infty}_{T,+}Y \wedge^{\mathbb{L}}\mathbf{E}]$, $\beta \in SH(k)[\Sigma^{\infty}_{T,+}Y,\Sigma^{\infty}_{T,+}Z \wedge^{\mathbb{L}}_S \mathbf{E}]$ and $\gamma \in SH(k)[\Sigma^{\infty}_{T,+}Z,\Sigma^{\infty}_{T,+}W \wedge^{\mathbb{L}}_S \mathbf{E}]$. Let's denote 
\[\beta \circ_M \alpha: \Sigma^{\infty}_{T,+}X \stackrel{\alpha}{\rightarrow} \Sigma^{\infty}_{T,+}Y \wedge^{\mathbb{L}}_S \mathbf{E} \stackrel{\beta \wedge \id_{\mathbf{E}}}{\longrightarrow} \Sigma^{\infty}_{T,+}Z \wedge^{\mathbb{L}}_S \mathbf{E} \wedge^{\mathbb{L}}_S \mathbf{E} \stackrel{\id_Z \wedge \mu_{\mathbf{E}}}{\longrightarrow} \Sigma^{\infty}_{T,+}Z \wedge^{\mathbb{L}}_S \mathbf{E}\]
and similarly for $\gamma \circ_M \beta$. Then $\circ_M$  is associative and unital.
\end{prop}
\begin{proof}
Both $\gamma \circ_M (\beta \circ_M \alpha)$ and $(\gamma \circ_M \beta) \circ_M \alpha$ are equal to the following composition
\begin{multline*}\Sigma^{\infty}_{T,+} X \stackrel{\alpha}{\rightarrow} \Sigma^{\infty}_{T,+} Y \wedge^{\mathbb{L}}_S \mathbf{E} \stackrel{\beta \wedge \id_{\mathbf{E}}}{\longrightarrow} \Sigma^{\infty}_{T,+}Z \wedge^{\mathbb{L}}_S \mathbf{E} \wedge^{\mathbb{L}}_S \mathbf{E}  \stackrel{\id_Z \wedge \mu_{\mathbf{E}}}{\longrightarrow} \Sigma^{\infty}_{T,+} Z \wedge^{\mathbb{L}}_S \mathbf{E} \stackrel{\gamma \wedge \id_{\mathbf{E}}}{\longrightarrow} \Sigma^{\infty}_{T,+}W \wedge^{\mathbb{L}}_S \mathbf{E} \wedge^{\mathbb{L}}_S \mathbf{E} \\ \stackrel{\id_W \wedge \mu_{\mathbf{E}}}{\longrightarrow} \Sigma^{\infty}_{T,+}W \wedge^{\mathbb{L}}_S \mathbf{E}.
\end{multline*} 
\end{proof}
\begin{defn}{\rm
The category of $\mathbf{E}$-correspondences $Corr_{\mathbf{E}}(k))$ is defined as: 
\[Obj(Corr_{\mathbf{E}}(k)) = Obj(SmProj(k))\] 
and 
\[Corr_{\mathbf{E}}(k)(X,Y) = SH(k)[\Sigma^{\infty}_{T,+}X,\Sigma^{\infty}_{T,+}Y \wedge^{\mathbb{L}}_S \mathbf{E}],\] 
where the composition 
\[\circ_M: Corr_{\mathbf{E}}(k)(X,Y) \otimes Corr_{\mathbf{E}}(k)(Y,Z) \rightarrow Corr_{\mathbf{E}}(k)(X,Z), \quad (\alpha,\beta) \mapsto \beta \circ_M \alpha\] 
is defined as 
\[\beta \circ_M \alpha: \Sigma^{\infty}_{T,+}X \stackrel{\alpha}{\rightarrow} \Sigma^{\infty}_{T,+}Y \wedge^{\mathbb{L}}_S \mathbf{E} \stackrel{\beta \wedge \id_{\mathbf{E}}}{\longrightarrow} \Sigma^{\infty}_{T,+}Z \wedge^{\mathbb{L}}_S \mathbf{E} \wedge^{\mathbb{L}}_S \mathbf{E} \stackrel{\id_Z \wedge \mu_{\mathbf{E}}}{\longrightarrow} \Sigma^{\infty}_{T,+}Z \wedge^{\mathbb{L}}_S \mathbf{E}. \]
}
\end{defn}
\begin{prop}
There is a functor 
\[\mathfrak{h}: SmProj(k) \rightarrow Corr_{\mathbf{E}}(k), \quad X \mapsto X,\]
which sends a morphism $f: X \rightarrow Y$ of $k$-schemes to 
\[\Sigma^{\infty}_{T,+}(f) \wedge \varphi_{\mathbf{E}}: \Sigma^{\infty}_{T,+}X = \Sigma^{\infty}_{T,+}X \wedge^{\mathbb{L}}_S S \rightarrow \Sigma^{\infty}_{T,+}Y \wedge^{\mathbb{L}}_S \mathbf{E}.\]
\end{prop}
\begin{proof}
The identity morphism in $Corr_{\mathbf{E}}(k)(X,X)$ is given by
\[\id_X \wedge \varphi_{\mathbf{E}}: \Sigma^{\infty}_{T,+}X = \Sigma^{\infty}_{T,+}X \wedge^{\mathbb{L}}_S S \rightarrow \Sigma^{\infty}_{T,+}X \wedge^{\mathbb{L}}_S \mathbf{E}.\]
Let $\alpha \in Corr_{\mathbf{E}}(k)(X,X)$ be an arbitrary $\mathbf{E}$-correspondence. By definition we have 
\[\alpha \circ_M (\id_X \wedge \varphi_{\mathbf{E}}) = (\id_X \wedge \mu_E) \circ (\alpha \wedge \id_E) \circ (\id_X \wedge \varphi_{\mathbf{E}}).\] 
Since $\mathbf{E}$ is a ring spectrum, we must have $\alpha \circ_M (\id_X \wedge \varphi_{\mathbf{E}}) = \alpha$. Similarly, $(\id_X \wedge \varphi_{\mathbf{E}}) \circ_M \alpha = \alpha$.
We check the compatibility of the composition laws. Let $X \stackrel{f}{\rightarrow} Y \stackrel{g}{\rightarrow} Z$ be morphisms of $k$-schemes. By definition we have 
\[\mathfrak{h}(g \circ f) = \Sigma^{\infty}_{T,+}(g \circ f) \wedge \varphi_{\mathbf{E}}\]
and 
\[\mathfrak{h}(g) \circ_M \mathfrak{h}(f) = (\id_Z \wedge \mu_{\mathbf{E}}) \circ (\Sigma^{\infty}_{T,+}(g) \wedge \varphi_{\mathbf{E}} \wedge \id_{\mathbf{E}}) \circ (\Sigma^{\infty}_{T,+}(f) \wedge \varphi_{\mathbf{E}}).\]  
The equality $\mathfrak{h}(g \circ f) = \mathfrak{h}(g) \circ_M \mathfrak{h}(f)$ follows from the fact that $\mathbf{E}$ is a ring spectrum. 
\end{proof}
Let $\mathbf{Spect}_T^{\Sigma}(k)$ be the model category of symmetric motivic $T$-spectra (\cite{Jar00}). Following \cite{CD10}, \cite[\S 2.2]{Deg13} we call  $\mathbf{E} \in \mathbf{Spect}_T^{\Sigma}(k)$ a strict motivic ring spectrum, if $\mathbf{E}$ is a commutative monoid object in $\mathbf{Spect}_T^{\Sigma}(k)$. An $\mathbf{E}$-module spectrum is a pair $(M,\gamma_M)$, where $M \in \mathbf{Spect}_{T}^{\Sigma}(k)$ and $\gamma_M: M \wedge \mathbf{E} \rightarrow M$, such that the following diagrams commute: 
\[\xymatrix{S \wedge M \ar[rr]^{\varphi_{\mathbf{E}}\wedge \id_M} \ar@{=}[drr] && \mathbf{E} \wedge M \ar[d]^{\gamma_M} \\ && M }\]
\[\xymatrix{\mathbf{E} \wedge \mathbf{E} \wedge M \ar[rr]^{\mu_{\mathbf{E}} \wedge \id_M} \ar[d]_{\id_{\mathbf{E}} \wedge \gamma_M} && \mathbf{E} \wedge M \ar[d]^{\gamma_M} \\ \mathbf{E} \wedge M \ar[rr]_{\gamma_M} && M }\]
Given two $\mathbf{E}$-modules $(M,\gamma_)$ and $(N,\gamma_N)$, an $\mathbf{E}$-module map is a map $f: M \rightarrow N$, such that the following diagram commutes:
\[\xymatrix{\mathbf{E} \wedge M \ar[rr]^{\id_{\mathbf{E}}\wedge f} \ar[d]_{\gamma_M} && \mathbf{E} \wedge N \ar[d]^{\gamma_N} \\ M \ar[rr]_f && N}\]
Given a strict motivic ring spectrum $\mathbf{E}$ one can form the model category $\mathbf{E}-Mod^{\Sigma}$ of $\mathbf{E}$-modules with respect to the symmetric monoidal model category $\mathbf{Spect}_{T}^{\Sigma}(k)$ (see e.g \cite{SS00})and there is a Quillen adjunction of model categories (we will return to this point in the last discussion in the Appendix):
\[-\wedge \mathbf{E} : \mathbf{Spect}_{T}^{\Sigma}(k) \leftrightarrows \mathbf{E}-Mod^{\Sigma} : U, \]
where $U$ denotes the forgetful functor. This Quillen adjunction induces an adjunction between homotopy categories: 
\begin{equation}\label{eqCorrE} -\wedge^{\mathbb{L}}_S \mathbf{E}: SH(k) \rightleftarrows Ho_k(\mathbf{E}-Mod) : RU,
\end{equation}
where we denote by $Ho_k(\mathbf{E}-Mod)$ the homotopy category associated to the category of strict $\mathbf{E}$-modules.  
\begin{thm}
Let $k$ be a field and $\mathbf{E} \in \mathbf{Spect}_T^{\Sigma}(k)$ be a strict motivic ring spectrum. There is a functor 
\[Corr_{\mathbf{E}}(k) \rightarrow Ho(\mathbf{E}-Mod), \quad X \mapsto \Sigma^{\infty}_{T,+}X \wedge^{\mathbb{L}}_S \mathbf{E}.\]
\end{thm} 
\begin{proof}
 Recall that we may regard $\Sigma^{\infty}_{T,+} X \wedge^{\mathbb{L}}_S \mathbf{E}$ as an $\mathbf{E}$-modules via the map 
\[\Sigma^{\infty}_{T,+}X \wedge^{\mathbb{L}}_S \mathbf{E} \wedge^{\mathbb{L}}_S \mathbf{E} \stackrel{\id_X \wedge \mu_{\mathbf{E}}}{\longrightarrow} \Sigma^{\infty}_{T,+} X \wedge^{\mathbb{L}}_S \mathbf{E}.\]
Let us denote the assocation above by 
\[F: Corr_{\mathbf{E}}(k) \rightarrow Ho(\mathbf{E}-Mod), \quad X \mapsto \Sigma^{\infty}_{T,+}X \wedge^{\mathbb{L}}_S \mathbf{E}.\]
$F$ maps on morphisms as following: Given $\alpha: \Sigma^{\infty}_{T,+}X \rightarrow \Sigma^{\infty}_{T,+}Y \wedge^{\mathbb{L}}_S \mathbf{E}$, we associate 
\[\Sigma^{\infty}_{T,+} X \wedge^{\mathbb{L}}_S \mathbf{E} \stackrel{\alpha \wedge \id_{\mathbf{E}}}{\longrightarrow} \Sigma^{\infty}_{T,+}Y \wedge^{\mathbb{L}}_S \mathbf{E} \wedge^{\mathbb{L}}_S \mathbf{E} \stackrel{\id_Y \wedge \mu_{\mathbf{E}}}{\longrightarrow} \Sigma^{\infty}_{T,+} Y \wedge^{\mathbb{L}}_S \mathbf{E}. \]
We have to check firstly, that $(\id_Y \wedge \mu_{\mathbf{E}}) \circ (\alpha \wedge \id_{\mathbf{E}})$ is a morphism of $\mathbf{E}$-modules. Since $\mathbf{E}$ is a ring spectrum, there is a commutative diagram 
\[\xymatrix{\Sigma^{\infty}_{T,+}X \wedge^{\mathbb{L}}_S \mathbf{E} \wedge^{\mathbb{L}}_S \mathbf{E} \ar[d]_{\id_X \wedge \mu_{\mathbf{E}}} \ar[rr]^{\alpha \wedge \id_{\mathbf{E}} \wedge \id_{\mathbf{E}} \quad \quad} && \Sigma^{\infty}_{T,+}Y \wedge^{ \mathbb{L}}_S \mathbf{E} \wedge^{\mathbb{L}}_S \mathbf{E} \wedge^{\mathbb{L}}_S \mathbf{E} \ar[rr]^{\quad \quad \id_Y \wedge \mu_{\mathbf{E}} \wedge \id_{\mathbf{E}}} && \Sigma^{\infty}_{T,+}Y \wedge^{\mathbb{L}}_S \mathbf{E} \wedge^{\mathbb{L}}_S \mathbf{E} \ar[d]^{\id_Y \wedge \mu_{\mathbf{E}}} \\ \Sigma^{\infty}_{T,+}X \wedge^{\mathbb{L}}_S \mathbf{E} \ar[rrrr]_{(\id_Y \wedge \mu_{\mathbf{E}}) \circ (\alpha \wedge \id_{\mathbf{E}})} &&&&  \Sigma^{\infty}_{T,+}Y \wedge^{\mathbb{L}}_S \mathbf{E} }\]
Now we have to check the compatibility of the composition laws. Given $\alpha: \Sigma^{\infty}_{T,+}X \rightarrow \Sigma^{\infty}_{T,+}Y \wedge^{\mathbb{L}}_S \mathbf{E}$ and $\beta: \Sigma^{\infty}_{T,+}Y \rightarrow \Sigma^{\infty}_{T,+}Z \wedge^{\mathbb{L}}_S \mathbf{E}$. Then $F(\beta)\circ F(\alpha)$ is the following composition 
\begin{multline*} \Sigma^{\infty}_{T,+}X \wedge^{\mathbb{L}}_S \mathbf{E} \stackrel{\alpha \wedge \id_{\mathbf{E}}}{\longrightarrow} \Sigma^{\infty}_{T,+}Y \wedge^{\mathbb{L}}_S \mathbf{E} \wedge^{\mathbb{L}}_S \mathbf{E} \stackrel{\id_Y \wedge \mu_{\mathbf{E}}}{\longrightarrow} \Sigma^{\infty}_{T,+}Y \wedge^{\mathbb{L}}_S \mathbf{E} \stackrel{\beta \wedge \id_{\mathbf{E}}}{\longrightarrow} \Sigma^{\infty}_{T,+}Z \wedge^{\mathbb{L}}_S \mathbf{E} \wedge^{\mathbb{L}}_S \mathbf{E} \\ \stackrel{\id_Z \wedge \mu_{\mathbf{E}}}{\longrightarrow} \Sigma^{\infty}_{T,+}Z \wedge^{\mathbb{L}}_S \mathbf{E}.
\end{multline*}
The composition $F(\beta \circ_M \alpha)$ is 
\[\Sigma^{\infty}_{T,+} X \wedge^{\mathbb{L}}_S \mathbf{E} \stackrel{(\beta \circ_M \alpha) \wedge \id_{\mathbf{E}}}{\longrightarrow} \Sigma^{\infty}_{T,+}Z \wedge^{\mathbb{L}}_S \mathbf{E} \wedge^{\mathbb{L}}_S \mathbf{E} \stackrel{\id_Z \wedge \mu_{\mathbf{E}}}{\longrightarrow} \Sigma^{\infty}_{T,+}Z \wedge^{\mathbb{L}}_S \mathbf{E},\]
where 
\[ \beta \circ_M \alpha: \Sigma^{\infty}_{T,+}X \stackrel{\alpha}{\rightarrow} \Sigma^{\infty}_{T,+}Y \wedge^{\mathbb{L}}_S \mathbf{E} \stackrel{\beta \wedge \id_{\mathbf{E}}}{\longrightarrow} \Sigma^{\infty}_{T,+}Z \wedge^{\mathbb{L}}_S \mathbf{E} \wedge^{\mathbb{L}}_S \mathbf{E} \stackrel{\id_Z \wedge \mu_{\mathbf{E}}}{\longrightarrow} \Sigma^{\infty}_{T,+}Z \wedge^{\mathbb{L}}_S \mathbf{E}.\]
Hence, $F(\beta \circ_M \alpha)$ is the following composition 
\begin{multline*}
\Sigma^{\infty}_{T,+} X \wedge^{\mathbb{L}}_S \mathbf{E} \stackrel{\alpha \wedge \id_{\mathbf{E}}}{\longrightarrow} \Sigma^{\infty}_{T,+}Y \wedge^{\mathbb{L}}_S \mathbf{E} \wedge^{\mathbb{L}}_S \mathbf{E} \stackrel{\beta \wedge \id_{\mathbf{E}} \wedge \id_{\mathbf{E}}}{\longrightarrow}  \Sigma^{\infty}_{T,+}Z \wedge^{\mathbb{L}}_S \mathbf{E} \wedge^{\mathbb{L}}_S \mathbf{E} \wedge^{\mathbb{L}}_S \mathbf{E} \\ \stackrel{\id_Z \wedge \mu_{\mathbf{E}} \wedge \id_{\mathbf{E}}}{\longrightarrow} \Sigma^{\infty}_{T,+}Z \wedge^{\mathbb{L}}_S \mathbf{E} \wedge^{\mathbb{L}}_S \mathbf{E} \stackrel{\id_Z \wedge \mu_{\mathbf{E}}}{\longrightarrow} \Sigma^{\infty}_{T,+}Z \wedge^{\mathbb{L}}_S \mathbf{E}.
\end{multline*}
Since $\mathbf{E}$ is a ring spectrum, we have a commutative diagram 
\[\xymatrix{\Sigma^{\infty}_{T,+}Y \wedge^{\mathbb{L}}_S \mathbf{E} \wedge^{\mathbb{L}}_S \mathbf{E} \ar[d]_{\id_Y \wedge \mu_{\mathbf{E}}} \ar[rr]^{\beta \wedge \id_{\mathbf{E}} \wedge \id_{\mathbf{E}} \quad \quad} && \Sigma^{\infty}_{T,+}Z \wedge^{ \mathbb{L}}_S \mathbf{E} \wedge^{\mathbb{L}}_S \mathbf{E} \wedge^{\mathbb{L}}_S \mathbf{E} \ar[rr]^{\quad \quad \id_Z \wedge \mu_{\mathbf{E}} \wedge \id_{\mathbf{E}}} && \Sigma^{\infty}_{T,+}Z \wedge^{\mathbb{L}}_S \mathbf{E} \wedge^{\mathbb{L}}_S \mathbf{E} \ar@{=}[d] \\ \Sigma^{\infty}_{T,+}Y \wedge^{\mathbb{L}}_S \mathbf{E} \ar[rrrr]_{\beta \wedge \id_{\mathbf{E}}} &&&&  \Sigma^{\infty}_{T,+}Y \wedge^{\mathbb{L}}_S \mathbf{E} } \]
This implies that $F(\beta \circ_M \alpha) = F(\beta) \circ F(\alpha)$.
\end{proof}
\begin{defn}{\rm
Let $k$ be a field and $\mathbf{E} \in \mathbf{Spect}_T^{\Sigma}(k)$ be a strict motivic ring spectrum. We define the category $\mathbf{Mot}_{\mathbf{E}}(k)$  of pure $\mathbf{E}$-motives over $k$ to be the smallest pseudo-abelian subcategory of $Ho_k(\mathbf{E}-Mod)$ generated as an additive category by $\{\Sigma^{\infty}_{T,+} X \wedge^{\mathbb{L}}_S \mathbf{E}| X \in SmProj(k) \}$. 
}
\end{defn}
\begin{rem}{\rm
We know that if $char(k) = 0$ then there is an equivalence of categories 
\[Ho_k(H\Z-Mod) \cong \mathbf{DM}(k),\] 
where $\mathbf{DM}(k)$ denotes the category of big Voevodsky's motives (cf. \cite{RO08}). As the category of pure Grothendieck-Chow motives $\underline{Chow}(k) \inj \mathbf{DM}(k)$ is embedded fully faithful into $DM(k)$, we raise a question: is $\mathbf{Mot}_{H\Z}(k)$ equivalent to $\underline{Chow}(k)$ via the equivalence above? We only know that
\begin{multline*}
\mathbf{Mot}_{H\Z}(k)(X,Y) \cong SH(k)[\Sigma^{\infty}_{T,+}X,\Sigma^{\infty}_{T,+}Y \wedge^{\mathbb{L}}_S H\Z] \cong H^{2(n_Y+d_Y),(n_Y+d_Y)}_M(X_+ \wedge Th(V_Y),\Z) \\ \cong H_M^{2d_Y,d_Y}(X\times Y,\Z) \cong \CH^{d_Y}(X\times Y),
\end{multline*}
where the first isomorphism comes from the adjunction 
\[ -\wedge^{\mathbb{L}}_S H\Z : SH(k) \leftrightarrows Ho_k(H\Z-Mod).\]
The second isomorphism comes from duality, the third isomorphism is the Thom isomorphism for motivic cohomology and the last isomorphism is the comparison isomorphism of Voevodsky (\cite[Cor. 19.2]{MVW06}). The question is, if these isomorphisms are compatible with the equivalence $Ho_k(H\Z-Mod) \cong \mathbf{DM}(k)$? It seems the problem with the first three isomorphisms is not difficult, however it seems that the problem with the last isomorphism is hard.  
}
\end{rem}
\begin{cor}\label{corMotE}
Let $k$ be a field and $\mathbf{E}$ be a strict motivic ring spectrum. There is a functor
\[\mathbf{Mot}_{\mathbf{E}}(k) \rightarrow SH(k).\]
\end{cor}
\begin{proof}
This follows from the adjunction \ref{eqCorrE}. 
\end{proof}
\subsection{Functoriality in motivic stable homotopy}
Following \cite{Ay08}, we recall that the stable homotopy category of schemes defines a $2$-functor from category of quasi-projective smooth schemes over a field $QSProjSm/k$ to the category of symmetric monoidal closed triangulated categories. Remark that the six operations formalism works much more general than what we here require. However we restrict ourselves only to $QSProjSm/k$, since it is already enough for our aim. We will list now a minimal list of properties of the six operations formalism: for any morphism of schemes $f: T \rightarrow S$, there is a pullback functor 
\[f^*: SH(S) \rightarrow SH(T),\]
such that $(f \circ g)^* = g^* \circ f^*$. Moreover, 
\begin{enumerate}
\item One has an adjunction for any morphism of schemes $f: T \rightarrow S$ 
\[f^*: SH(S) \leftrightarrows SH(T) : f_*.\]
If $f$ is smooth, then one has an adjunction 
\[f_{\#}: SH(T) \leftrightarrows SH(S) : f^* \]
\item Given a cartesian square 
\begin{equation*} \xymatrix{ Y \ar[r]^q \ar[d]_g & X \ar[d]^f \\ T \ar[r]_p & S}
\end{equation*}
and assume $f$ is smooth, then 
\[f_{\#}p^* \stackrel{\cong}{\rightarrow} g_{\#}q^*\]
\item Let $f: Y \rightarrow X$ be a smooth morphism, $\sE \in SH(Y)$ and $\sF \in SH(X)$, the natural transformation 
\[f_{\#}(\sE \wedge f^* \sF) \stackrel{\cong}{\rightarrow} f_{\#}\sE \wedge \sF\]
is an isomorphism.
\item Let $i: Z \inj X$ be a closed immersion with complement $j: U \inj X$, then there is a distinguished triangle 
\[j_{\#}j^* \rightarrow Id \rightarrow i_*i^* \stackrel{+1}{\rightarrow}\]
\item For any closed immersion $i: Z \inj X$, one has an adjunction 
\[i_*: SH(Z) \leftrightarrows SH(X) : i^! \]
\item Given a cartesian square 
\[\xymatrix{ T \ar[r]^k \ar[d]_g & Y \ar[d]^f \\ Z \ar[r]_i & X}\]
where $i: Z \inj X$ is a closed immersion, then one has an isomorphism 
\[f^*i_* \stackrel{\cong}{\rightarrow} k_*g^*\]
\item Let $i: Z \inj X$ be a closed immersion, $\sE \in SH(Z)$ and $\sF \in SH(X)$, the natural transformation 
\[i_*(\sE \wedge i^*\sF) \stackrel{\cong}{\rightarrow} i_* \sE \wedge \sF\]
is an isomorphism.
\item For any separated morphism of finite type $f: Y \rightarrow X$, there is an adjunction 
\[f_! : SH(Y) \leftrightarrows SH(X) : f^!.\]
\item For a smooth separated morphism of finite type $f: Y \rightarrow X$ with the relative tangent bundle $T_f$ there are canonical natural isomorphisms, which are dual to each other
\[f_{\#} \stackrel{\cong}{\longrightarrow} f_!(Th_Y(T_f) \wedge_Y -), \quad f^* \stackrel{\cong}{\longrightarrow} Th_Y(-T_f) \wedge_Y f^!.\]
Moreover, for any separated morphism of finite type $f : Y \rightarrow X$, there exist natural isomorphisms
\[Ex(f^*_!,\wedge) : (f_!K) \wedge_X L \stackrel{\cong}{\longrightarrow} f_!(K \wedge_Y f^*L),\] 
\[\underline{\Hom}_X(f_!L,K) \stackrel{\cong}{\longrightarrow} f_*\underline{\Hom}_Y(L,f^!K), \]
\[f^!\underline{\Hom}_X(L,M) \stackrel{\cong}{\longrightarrow} \underline{\Hom}_Y(f^*L,f^!M).\]
\item If $f : Y \rightarrow X$ is a smooth projective morphism then $f_{\#}(\mathbf{1}_Y)$ is strongly dualizable in $SH(X)$ with the dual 
\[D_X(f_{\#}(\mathbf{1}_Y)) = f_{\#}Th_Y(-T_f).\]
Furthermore, one has $D_X(f_*K) \cong f_*D_Y(K \wedge_Y Th_Y(T_f)), \forall K\in SH(Y)$. 
\end{enumerate}
We will need some facts about cohomology with supports in the next subsection. 
\subsection{Cohomology with supports}
Let $S = \Spec k$. We consider the category $SH(k)$. For a ring spectrum $\mathbf{E} \in SH(k)$ and a closed pair $(X,Z)$, where $\pi_X: X \rightarrow S$ is a smooth quasi-projective $k$-scheme and $i: Z \inj X$ a smooth closed subscheme, one defines the cohomology with support as 
\[\mathbf{E}^{p,q}_{Z}(X) = SH(S)[X/X-Z,\mathbf{E}\wedge S^{p,q}] \cong SH(X)[i_*(\mathbf{1}_Z),\mathbf{E}_X \wedge S^{p,q}] \cong SH(Z)[\mathbf{1}_Z,i^! \mathbf{E}_X \wedge S^{p,q}] ,\]
where we write $\mathbf{E}_X = \pi_X^*\mathbf{E}$. As $\Sigma^{\infty}_{T,+}X/X-Z := \pi_{X\#}i_*(\mathbf{1}_Z)$ in $SH(k)$, so the first isomorphism follows from the adjunction
\[\pi_{X \#} : SH(X) \leftrightarrows SH(S) : \pi_X^* \]
and the last isomorphism comes from the adjunction
\[i_* : SH(Z) \leftrightarrows SH(X) : i^!.\]
If $f: Y \rightarrow X$ is a smooth morphism of smooth quasi-projective $S$-schemes we have a canonical homomorphism 
\[f^*:  \mathbf{E}^{p,q}_Z(X) \rightarrow \mathbf{E}^{p,q}_T(Y),\]
where $T = Y\times_X Z$ defined as following: Consider the commutative diagram
\[\xymatrix{  T \ar[r]^j \ar[d]_g & Y \ar[d]^f \\ Z \ar[r]_i & X}\]
For a morphism 
\[\alpha: i_*\mathbf{1}_Z \rightarrow \mathbf{E}_X \wedge S^{p,q}\] 
we can associate to a morphism 
\[f^* \alpha: j_*\mathbf{1}_T \cong j_*g^*\mathbf{1}_Z \stackrel{(Ex^*_*)^{-1}}{\cong} f^*i_*\mathbf{1}_Z  \rightarrow f^*\mathbf{E}_X \wedge S^{p,q} \stackrel{\cong}{\rightarrow} \mathbf{E}_Y \wedge S^{p,q}.\] 
If $T \stackrel{j}{\inj} Z \stackrel{i}{\inj} X$ are closed immersions, we can define a pushforward on cohomology with supports 
\[j_! : \mathbf{E}^{p,q}_T(X) \rightarrow \mathbf{E}^{p,q}_Z(X) \]
as following: Given a morphism $\alpha: X/X-T \rightarrow \mathbf{E}_X \wedge S^{p,q}$ we associate $j_!(\alpha) = \alpha \circ \bar{j}$, where $\bar{j}: X/X-Z \rightarrow X/X-T$ is the canonical morphism in $SH(X)$ induced by the immersion $X-Z \inj X-T$.  If $\alpha \in \mathbf{E}^{p,q}_Z(X)$ and $\beta \in \mathbf{E}^{m,n}_Z(X)$ we define their product in $\mathbf{E}^{p+m,q+n}_Z(X)$ as a morphism 
\[\alpha \cup \beta: X/X-Z \stackrel{\Delta}{\rightarrow}  X/(X-Z) \wedge X/(X-Z) \stackrel{\alpha \wedge \beta}{\longrightarrow} \mathbf{E} \wedge^{\mathbb{L}} \mathbf{E} \wedge S^{p+m,q+n} \stackrel{\mu_{\mathbf{E}}}{\rightarrow} \mathbf{E} \wedge S^{p+m,q+n}.\]
If $\xi/X$ is a vector bundle over a smooth $k$-scheme $X$ with the zero section $s_0: X \rightarrow \xi$, then the $\mathbf{E}$-cohomology of the Thom spectrum $Th(\xi)$ is 
\[\mathbf{E}^{p,q}(Th(\xi)) = \mathbf{E}^{p,q}_{X}(\xi).\]
The pushforward defined as above works only for closed immersions. We will define later pushforward on $\mathbf{E}$-cohomology of Thom spectrum for projective smooth morphism using duality. 
\subsection{Relation to the category of twisted $E$-correspondences}
\begin{nota}{\rm
For a quasi-projective smooth $k$-scheme $\pi_X: X \rightarrow \Spec k$ and a vector bundle $p_{\xi}: \xi \rightarrow X$ with $0$-section $s_X: X \rightarrow \xi$ we will write $Th_X(\xi) = p_{\xi\#}s_{X!}(\mathbf{1}_X)$ for the Thom transformation $Th(s_X,p_{\xi}) = p_{\xi\#}s_{X!}$ applying on $\mathbf{1}_X$. $Th_X(\xi)$ is an object in $SH(X)$ and $Th_X(-\xi) = s_X^!p_{\xi}^*(\mathbf{1}_X)$ for its inverse as the inverse Thom transformation $Th^{-1}(s_X,p_{\xi})$ applying on $\mathbf{1}_X$. The Thom spectrum will be denoted by $Th(\xi/X)$, which means
\[Th(\xi/X) = \pi_{X\#}Th_X(\xi) = \pi_{X\#} p_{\xi\#} s_{X!}(\mathbf{1}_X) \cong \pi_{X\#}p_{\xi\#}s_{X!}\pi_X^*\mathbf{1}_k.\]   
Sometime we only write $Th(\xi)$ for the Thom spectrum, if it is clear which scheme $X$ we talk about. One can see easily that this definition coincides with the traditional definition of Thom spectrum as follow: Let $j: \xi - s_X(X) \inj \xi$ be the open immersion with the complement $s_X : X \rightarrow \xi$. One has a localization sequence 
\[j_{\#}j^*(\mathbf{1}_{\xi}) \rightarrow \mathbf{1}_{\xi} \rightarrow s_{X*}s_X^*(\mathbf{1}_{\xi}).\] 
Applying $\pi_{\#}p_{\xi \#}$ and as $s_X{*} \cong s_{X!}$ is a natural $2$-isomorphism one has a natural isomorphism in $SH(k)$: 
\[Th_X(\xi) \cong \Sigma_{T,+}^{\infty}\xi/\xi-s_X(X).\]
}
\end{nota}
Let $E \in SH(k)$ be a ring spectrum and $X/k$ a quasi-projective smooth $k$-scheme. Let $p_{\xi}: \xi \rightarrow X$ be a vector bundle of rank $r$ with the zero section $s: X \rightarrow \xi$. We define $E$-cohomology of $X$ twisted by a vector bundle as 
\[E^{p,q}(X, \xi) = SH(X)[\mathbf{1}_X,s^!p_{\xi}^*E_X^{2r,r} \wedge S^{p,q}] = SH(X)[\mathbf{1}_X,Th_X(-\xi)\wedge_X E^{2r,r}_X \wedge S^{p,q}].\]
where we write $E_X^{2r,r} = E_X \wedge S^{2r,r}$. We denote by $E^{*,*}(X,\xi)$ the bigraded ring
\[E^{*,*}(X,\xi) = \oplus_{p,q} E^{p,q}(X,\xi).\]
Remark that $E^{*,*}(X,\xi)$ is bigraded ring. Even if $E$ is a commutative ring spectrum, $E^{*,*}(X,\xi)$ is never bigraded commutative. If $\xi \in \sV(X)$ is a virtual vector bundle of rank $r < 0$ then $p_{-\xi}: -\xi \rightarrow X$ is an actual vector bundle, so we define 
\[E^{p,q}(X,\xi) = SH(X)[\mathbf{1}_X,p_{-\xi}s_!E_X^{2r,r}\wedge S^{p,q}].\]
This group has the following interpretation by Jouanolou trick: As $X$ is quasi-projective, so we have an immersion $i: X \inj \P^N$. Via the Segre embedding $\P^N \times \P^N \inj \P^{N^2+2N}$, $U$ is an affine variety. Let 
\[U = \P^N\times \P^N - Proj \, k[x_0,\cdots,x_N,y_0,\cdots,y_N]/\sum_{i=0}^N x_iy_i.\]
$pr_1: U \rightarrow \P^N$ is an $\A^N$-bundle. Consider the pullback diagram
\[\xymatrix{i^*U \ar[d]_{\pi} \ar[r] & U \ar[d]^{pr_1} \\ X \ar[r]_i & \P^N} \]
Then 
\[E^{p,q}(X,\xi) \cong SH(k)[Th(\zeta/U),E^{2(r+n),(r+n)} \wedge S^{p,q}],\]
where $\zeta$ is an actual vector bundle on $U$, such that $\pi^*\xi \oplus \sO^n = \zeta$. 
\begin{prop}
Let $f: \xi \stackrel{\cong}{\longrightarrow} \xi'$ be an isomorphism of vector bundles on $X$ 
\[\xymatrix{ \xi \ar[rr]^{f}_{\cong} \ar[dr]_{p_{\xi}} && \xi' \ar[dl]^{p_{\xi'}} \\ & X  }\]
There is a natural isomorphism 
\[E^{p,q}(X, \xi) \cong E^{p,q}(X, \xi').\]
\end{prop}
\begin{proof}
Consider the Cartesian squares
\[\xymatrix{X \ar@{=}[d] \ar[r]^{s} & \xi \ar[d]^{\cong}_f \ar[r]^{p_{\xi}} & X \ar@{=}[d] \\ X \ar[r]_{s'} & \xi' \ar[r]_{p_{\xi'}} & X }\]
One has two $2$-isomorphisms (\cite[\S 1.5.5]{Ay08})
\[Th_X(s,p_{\xi}) \stackrel{\cong}{\longrightarrow} Th_X(s',p'_{\xi}),\]
and
\[Th_X^{-1}(s',p'_{\xi}) \stackrel{\cong}{\longrightarrow} Th_X^{-1}(s,p_{\xi}),\]
which prove the Proposition. 
\end{proof}
\begin{prop}
If $E$ is orientable in sense of \cite[Def. 12.2.2]{CD10}, then there is a natural isomorphism 
\[E^{p,q}(X, \xi) \stackrel{\cong}{\longrightarrow} E^{p,q}(X).\]
\end{prop}
\begin{proof}
Since $E$ is orientable, one has by \cite[Thm. 2.4.50 (3)]{CD10} a canonical natural isomorphism 
\[p_{\xi}^*E_X \stackrel{\cong}{\rightarrow} p_{\xi}^!E_X \wedge S^{-2r,-r}.\]
This induces a natural isomorphism 
\[E^{p,q}(X,\xi) = SH(X)[\mathbf{1}_X,s^!p_{\xi}^*E_X^{2r,r} \wedge S^{p,q}] \stackrel{\cong}{\longrightarrow} SH(X)[\mathbf{1}_X,s^!p_{\xi}^!E_X \wedge S^{p,q}] = E^{p,q}(X). \]
\end{proof}
\begin{prop}(twisted Thom isomorphism) \label{isothomE} Let $X/k$ be a quasi-projective smooth $k$-scheme and $p_{\xi}: \xi \rightarrow X$ be a vector bundle of rank $r$ with the zero section $s : X \rightarrow \xi$. One has a natural isomorphism 
\[th_E^X(\xi): E^{p,q}(X,\xi) \cong E^{p+2r,q+r}(Th(\xi)),\]
which we call the twisted Thom isomorphism. 
\end{prop}
\begin{proof}
We have two adjunctions 
\[s_!: SH(X) \rightleftarrows SH(\xi): s^!, \quad p_{\xi \#} : SH(\xi) \rightleftarrows SH(X): p_{\xi}^* .\]
Hence, we have 
\begin{multline*}
E^{p,q}(X,\xi) = SH(X)[\mathbf{1}_X,s^!p_{\xi}^*E_X^{2r,r} \wedge S^{p,q}] \cong \\ \cong SH(X)[p_{\xi \#}s_!(\mathbf{1}_X),E_X^{2r,r} \wedge S^{p,q} ] \cong E^{p+2r,q+r}(Th(\xi)), 
\end{multline*}
where the last natural isomorphism is induced by the adjunction $(\pi_{X\#},\pi_X^*)$, where $\pi_X: X \rightarrow \Spec k$ is the structure morphism.  
So for a morphism 
\[\alpha : \mathbf{1}_X \rightarrow s^!p_{\xi}^*E_X^{2r,r} \wedge S^{p,q}\]
the twisted Thom isomorphism is explicitly given by 
\[th_E^X(\xi)(\alpha) = \varepsilon_{(\pi_{X\#},\pi_X^*)}\circ \pi_{X\#} \circ \varepsilon_{(p_{\xi \#},p_{\xi}^*)} \circ p_{\xi \#} \circ \varepsilon_{(s_!,s^!)} \circ s_!(\alpha). \]
\end{proof}
\begin{ex}{\rm
The twisted Chow-Witt group $\widetilde{\CH}^p(X,\det \xi)$ defined by J. Fasel (cf. \cite{Fas07} and \cite{Fas08}) and also by F. Morel (\cite{Mor12}) is an example of twisted cohomology. One has a natural isomorphism 
\[\widetilde{\CH}^p(X,\det \xi) \stackrel{defn}{=} H^p_{Nis}(X,\underline{\mathbf{K}}^{MW}_p(\det \xi)) \cong H(\underline{\mathbf{K}}^{MW}_*)^{2p,p}(X,\xi),\]
where $H(\underline{\mathbf{K}}^{MW}_*)$ denotes the Eilenberg-Maclane spectrum associated to the homotopy module $\underline{\mathbf{K}}^{MW}_*$. We will discuss later about $H\underline{\mathbf{K}}^{MW}_*$ after introducing the homotopy t-structure. 
}
\end{ex}
Before going further, we want to give a list of properties of Thom transformations that we will need for our constructions. 
\begin{prop}\cite[Prop. 2.3.19]{Ay08}\label{Ayoub1}
Let $X$ be a quasi-projective $k$-scheme and $\xi/X$ be a vector bundle. Let $f: Y \rightarrow X$ be a morphism. Then one has two natural $2$-isomorphisms
\[f^*Th_X(\xi) \stackrel{\cong}{\longrightarrow} Th_Y(f^*\xi)f^*, \quad f^*Th_X(-\xi) \stackrel{\cong}{\longrightarrow} Th_Y(-f^*\xi) f^* ,\] 
which satisfy: For all $(K,L) \in Obj(SH(X)^2)$, there are two commutative diagrams
\[\xymatrix{f^*K\wedge_Y(f^*Th_X(\xi)L) \ar[d]_{\cong} \ar[r]^{\cong} & f^*(K\wedge_XTh_X(\xi) L) \ar[r]^{\cong} & f^*Th_X(\xi)(K\wedge_{X}L) \ar[d]^{\cong} \\ f^*K \wedge_Y Th_Y(f^*\xi)f^*L \ar[r]_{\cong} & Th_Y(f^*\xi)(f^*K \wedge_Y f^*L) \ar[r]_{\cong} & Th_Y(f^*\xi)f^*(K\wedge_X L) }\] 
and 
\[\xymatrix{f^*K\wedge_Y(f^*Th_X(-\xi)L) \ar[d]_{\cong} \ar[r]^{\cong} & f^*(K\wedge_XTh_X(-\xi) L) \ar[r]^{\cong} & f^*Th_X(-\xi)(K\wedge_{X}L) \ar[d]^{\cong} \\ f^*K \wedge_Y Th_Y(-f^*\xi)f^*L \ar[r]_{\cong} & Th_Y(-f^*\xi)(f^*K \wedge_Y f^*L) \ar[r]_{\cong} & Th_Y(-f^*\xi)f^*(K\wedge_X L) }\]
\end{prop}
\begin{prop}\cite[Prop. 2.3.20]{Ay08}\label{Ayoub2}
Let $f: Y \rightarrow X$ be a $k$-morphism of quasi-projective schemes and $\xi/X$ be a vector bundle. There are two natural $2$-isomorphisms 
\[Th_X(\xi)f^* \stackrel{\cong}{\longrightarrow} f_*Th_X(f^*\xi), \quad Th_X(-\xi)f^* \stackrel{\cong}{\longrightarrow} f_*Th_Y(-f^*\xi),\] 
such that the following diagrams commute for all $(K,L) \in SH(X) \times SH(Y)$:
\[\xymatrix{K\wedge_X Th_X(\xi)f_*L \ar[d]_{\cong} \ar[r]^{\cong} & K\wedge_X f_*Th_Y(f^*\xi) L \ar[r]^{\cong} & f_*(f^*K \wedge_Y Th_Y(f^*\xi)L) \ar[d]^{\cong} \\ Th_X(\xi)(K \wedge_X f_*L) \ar[r]_{\cong} & Th_X(\xi)f_* (f^*K \wedge_Y L) \ar[r]_{\cong} & f_*Th_Y(f^*\xi)(f^*K\wedge_X L) }\] 
\[\xymatrix{K\wedge_X Th_X(-\xi)f_*L \ar[d]_{\cong} \ar[r]^{\cong} & K\wedge_X f_*Th_Y(-f^*\xi) L \ar[r]^{\cong} & f_*(f^*K \wedge_Y Th_Y(-f^*\xi)L) \ar[d]^{\cong} \\ Th_X(-\xi)(K \wedge_X f_*L) \ar[r]_{\cong} & Th_X(-\xi)f_* (f^*K \wedge_Y L) \ar[r]_{\cong} & f_*Th_Y(-f^*\xi)(f^*K\wedge_X L) }\] 
\end{prop} 
Let $f: Y \rightarrow X$ be any morphism of finite type and separated of quasi-projective smooth $k$-schemes. In the following we define a pullback map on twisted $E$-cohomology 
\[E^{p,q}(X,\xi) \rightarrow E^{p,q}(Y,f^* \xi).\] 
Consider the functor $f^*: SH(X) \rightarrow SH(Y)$. $f^*$ induces a map 
\begin{multline*}
E^{p,q}(X,\xi) = SH(X)[\mathbf{1}_X,s^!p_{\xi}^*E_X^{2r,r} \wedge S^{p,q}] \rightarrow SH(Y)[f^*\mathbf{1}_X,f^*s^!p_{\xi}^*E_X^{2r,r} \wedge S^{p,q}] = \\ = SH(Y)[\mathbf{1}_Y,f^*s^!p_{\xi}^*E_X^{2r,r} \wedge S^{p,q}].
\end{multline*}
Let $s_Y$ be the $0$-section of the vector bundle $p_{f^*\xi}: f^*\xi \rightarrow Y$ and we write $f_{\xi}: f^*\xi \rightarrow \xi$. One has an exchange transformation (see \cite[Prop. 1.4.15]{Ay08}) 
\[Ex^{*!}:f^*s^! \rightarrow s_Y^!f_{\xi}^*,\]
which is the following composition ($s_* \cong s_!, s_{Y!} \cong s_{Y*}$ since $s$ and $s_Y$ are closed immersion): 
\[f^*s^! \stackrel{\eta_{(s_{Y*},s_Y^!)}}{\longrightarrow} s_Y^!s_{Y*}f^*s^! \stackrel{Ex^*_*(-)^{-1}}{\longrightarrow} s_Y^!f_{\xi}^*s_*s^! \stackrel{\varepsilon_{(s_*,s^!)}}{\longrightarrow} s_Y^!f_{\xi}^*,\]
where $f_{\xi}: f^*\xi \rightarrow \xi$ is the induced map on vector bundles. Note that the exchange transformation $Ex^{*!}$ is an isomorphism, when $f$ is smooth (\cite[Cor. 1.4.17]{Ay08}). At this point we also notice that for an actual bundle $\xi$, the Thom transformation $Th_X(\xi)$ behaves well under pullback of a general morphism, since
\[f^*Th_X(\xi) = f^*p_{\xi \#}s_{X!} \cong p_{f^*\xi\#}f^*_{\xi}s_{X*} \cong p_{f^*\xi\#}s_{Y*}f^* \cong p_{f^*\xi \#}s_{Y!}f^* = Th_Y(f^*\xi),\]
and we have a natural transformation 
\[Ex^{*!} : f^*Th_X(-\xi) = f^*s_X^!p_{\xi}^* \rightarrow s_Y^!f^*_{\xi}p_{\xi}^* \cong s_Y^!p_{f^*\xi}^*f^* = Th_Y(-f^*\xi),\]
which is an isomorphism, if $f$ is smooth (see \cite[Lem. 1.5.4]{Ay08}). However, he showed that $Th_X(\xi)$ and $Th_X(-\xi)$ are inverse to each other \cite[Thm. 1.5.7]{Ay08}, hence $Th_Y(-f^*\xi) \stackrel{\cong}{\longrightarrow} f^*Th_X(-\xi)$ (cf. \cite[Rem. 1.5.10]{Ay08}) for all morphism not necessary smooth $f$. That is a very crucial point. Now consider the pullback diagram 
\[\xymatrix{ f^*\xi \ar[d]_{p_{f^* \xi}} \ar[r]^{f_{\xi}} & \xi \ar[d]^{p_{\xi}} \\ Y \ar[r]_f & X }\]
We have a natural isomorphism $f^*_{\xi}p^*_{\xi} \cong p^*_{f^*\xi}f^*$. Hence, we obtain a map 
\[E^{p,q}(X,\xi) \rightarrow E^{p,q}(Y,f^* \xi),\]
which we define as pullback of twisted $E$-cohomology. 
\begin{rem}{\rm
The composition of pullback on twisted $E$-cohomology $g^* \circ f^*$ is only defined up to the natural isomorphism $(f\circ g)^* \stackrel{\cong}{\longrightarrow} g^* f^*$.
}
\end{rem}
\begin{rem}{\rm
Let $a : \xi \stackrel{\cong}{\longrightarrow} \xi$ be an automorphism of a vector bundle $\xi$ of rank $r$ on $X$. Then one has the cartesian squares
\[\xymatrix{X \ar[r]^{s_X'} \ar@{=}[d]_{\id} & \xi \ar[d]_{a}^{\cong} \ar[r]^{p'_{\xi}} & X \ar@{=}[d]^{\id} \\  X \ar[r]_{s_X} & \xi \ar[r]_{p_{\xi}} & X }\]  
As in \cite[\S 1.5.5 p. 84]{Ay08} $a$ induces two $2$-isomorphisms between the Thom transformations
\[\omega(a) : Th(s'_X,p'_{\xi}) = p'_{\xi\#}s'_{X!} \stackrel{\cong}{\longrightarrow} Th(s_X,p_{\xi}) = p_{\xi\#}s_{X!} \]
and
\[\omega_{-1}(a): Th^{-1}(s_X,p_{\xi}) = s^!_X \pi^*_{\xi} \stackrel{\cong}{\longrightarrow} Th^{-1}(s'_{X},p'_{\xi}) = {s'}^!_{X} {p'}^*_{\xi}.\]
This induces an isomorphism, which is not necessary identity
\[\bar{\omega}(a):  SH(X)[\mathbf{1}_X,s_X^!p_{\xi}^*E_X^{2r,r} \wedge S^{p,q}] \stackrel{\cong}{\longrightarrow} SH(X)[\mathbf{1}_X,{s'}_X^!{p'}_{\xi}^*E_X^{2r,r} \wedge S^{p,q}].\]
However, the two pullbacks induced on twisted $E$-cohomology along a morphism $f: Y \rightarrow X$ must not be on the same target $E^{p,q}(Y,f^*\xi)$, as there are two different pullback diagrams 
\[\xymatrix{f^*\xi \ar[d]_{p_{f^*\xi}} \ar[r] & \xi \ar[d]^{p_{\xi}} &&&&  f^*\xi \ar[r] \ar[d]_{p'_{f^*\xi}} & \xi \ar[d]^{p'_{\xi}} \\ Y \ar[r]_f & X &&&& Y \ar[r]_f & \xi }\]
Consequently, there is no problem with maps between $E$-cohomology created by automorphisms of $\xi$.   
}
\end{rem}
\begin{rem}{\rm
Thanks to the Proposition \ref{Ayoub1}. The pullback of cohomology of virtual vector bundles is defined in the same way.
}
\end{rem}
\begin{prop}\label{functorialityPullback}
Let $Z \stackrel{g}{\rightarrow} Y \stackrel{f}{\rightarrow} X$ be morphisms of quasi-projective smooth $k$-schemes. Let $\xi/X$ be a vector bundle. Then one has up to a natural isomorphism induced by a natural $2$-isomorphism 
\[(f \circ g)^* = g^* \circ f^*: E^{p,q}(X,\xi) \rightarrow E^{p,q}(Z,g^*f^* \xi).\]
\end{prop}
\begin{proof}
Consider the chain of pullback bundles 
\[\xymatrix{g^*f^*\xi \ar[d]_{p_{g^*f^*\xi}} \ar[r]^{g_{\xi}} & f^*\xi \ar[d]_{p_{f^*\xi}} \ar[r]^{f_{\xi}} & \xi \ar[d]^{p_{\xi}} \\ Z \ar[r]_g & Y \ar[r]_f & X }\] 
Let $s_X, s_Y$ and $s_Z$ be the $0$-sections of $\xi, f^*\xi$ and $g^*f^*\xi$ respectively. The functoriality up to a natural isomorphism follows easily from the natural $2$-isomorphism 
\[(f\circ g)^*s_X^! \cong g^*f^*s^!_X \stackrel{Ex^{*!}}{\longrightarrow} g^*s_Y^!f_{\xi}^* \stackrel{Ex^{*!}}{\longrightarrow} s_Z^!g_{\xi}^*f_{\xi}^* \cong s_Z^!(f_{\xi} \circ g_{\xi})^*.\]
\end{proof}
This motivates us to give the following definition: 
\begin{defn}{\rm
A twisted $E$-cohomology pre-theory is an association, which is contravariant in both variables:  
\[E^{*,*}(-,-): QSProjSm(k)\times \sV \supset \sA \rightarrow \mathbf{Ring}^*,\] 
where $\mathbf{Ring}^*$ denotes the category of bigraded rings and $\sV$ is the 2-category, where objects are categories of virtual vector bundles $\sV(X)$ for $X \in QSProjSm(k)$ and 
\[1-Mor_{\sV}(\sV(X),\sV(Y)) = Fun(\sV(X),\sV(Y))\]
\[2-Mor_{\sV}(F,G) = \mathbf{Nat}(F,G).\]
$\sA$ is the full subcategory of $QSProjSm(k) \times \sV$ consisting of those pairs $(X,\xi)$, where $X \in QSProjSm(k)$ and $\xi \in \sV(X)$. $Mor_{\sA}((X,\xi),Y(,\eta))$ consists of pair $(f,g)$, where $f: X \rightarrow Y$ is a morphism of quasi-projective smooth $k$-schemes and $g: \xi \rightarrow \eta$ is a bundle map
\[\xymatrix{\xi \ar[r]^g \ar[d] & \eta \ar[d] \\ X \ar[r]_f & Y}\]
such that $\xi \rightarrow f^*\eta$ is a monomorphism in $\sV(X)$. $E^{*,*}(-,-)$ sends such a pair $(X,\xi)$ to $E^{*,*}(X,\xi)$. Given an $\sA$-morphism $(f,g): (X,\xi) \rightarrow (Y,\eta)$, $E^{*,*}(-,-)$ sends $(f,g)$ to the following composition
\[E^{*,*}(Y,\eta) \stackrel{f^*}{\longrightarrow} E^{*,*}(X,f^*\eta) \longrightarrow E^{*,*}(X,\xi),\]
where $f^*$ is the pullback map on twisted $E$-cohomology constructed as above and the last map is induced by 
\[Th_X(\xi) \rightarrow Th_X(f^*\eta),\]
as $\xi \rightarrow f^*\eta$ is a monomorphism in $\sV(X)$.   
}
\end{defn}
\begin{prop}\label{propthEpullback}
Let $f : Y \rightarrow X$ be a $k$-morphism of quasi-projective smooth $k$-schemes and $p_{\xi}: \xi \rightarrow X$ be a vector bundle of rank $r$ on $X$. There is a commutative diagram up to a natural isomorphism induced by a natural $2$-isomorphism
\[\xymatrix{E^{p,q}(X,\xi) \ar[r]^{f^*} \ar[d]^{\cong}_{th^X_E(\xi)} & E^{p,q}(Y,f^*\xi) \ar[d]_{\cong}^{th^Y_E(f^*\xi)} \\ E^{p+2r,q+r}(Th(\xi)) \ar[r]_{f^*} & E^{p+2r,q+2r}(Th(f^* \xi)) }\]
where $f^*: E^{p+2r,q+r}(Th(\xi)) \rightarrow E^{p+2r,q+r}(Th(f^*\xi))$ is the pullback given by 
\begin{multline*}
SH(X)[Th_X(\xi),E_X^{2r,r} \wedge S^{p,q}] \stackrel{f^*}{\rightarrow} SH(Y)[f^*Th_X(\xi),f^*E_X^{2r,r} \wedge S^{p,q}] \stackrel{Ex^{!*}}{\longrightarrow} \\ SH(Y)[Th_Y(f^*\xi),E_Y^{2r,r} \wedge S^{p,q}].
\end{multline*}
\end{prop}
\begin{proof}
It is obvious by construction. Remark that for general morphism $f$ we always have the pullback 
\[SH(X)[Th_X(\xi),E_X^{2r,r} \wedge S^{p,q}] \rightarrow SH(Y)[Th_Y(f^*\xi),E_Y^{2r,r} \wedge S^{p,q}].\] 
Since $\pi_X : X \rightarrow \Spec k$ and $\pi_Y : Y \rightarrow \Spec k$ are smooth, then one has the natural isomorphisms via the adjunctions $(\pi_{X\#},\pi_X^*)$ and $(\pi_{Y\#},\pi_Y^*)$:
\[SH(X)[Th_X(\xi),E_X^{2r,r} \wedge S^{p,q}] \cong SH(k)[\pi_{X\#}Th_X(\xi),E^{2r,r} \wedge S^{p,q}] = E^{p,q}(Th(\xi)),\]
and
\[SH(Y)[Th_Y(f^*\xi),E_Y^{2r,r} \wedge S^{p,q}] \cong SH(k)[\pi_{Y\#}Th_Y(f^*\xi),E^{2r,r}\wedge S^{p,q}] = E^{p,q}(Th(f^*\xi)).\]
Explicitly, given a morphism 
\[\alpha: \mathbf{1}_X \rightarrow s_X^!p_{\xi}^*E_X^{2r,r} \wedge S^{p,q}\]
we have 
\[f^*th_E^X(\xi)(\alpha) = \varepsilon_{(s_{X*},s_X^!)} \circ Ex_*^*(-)^{-1} \circ \eta_{(s_{Y*},s_Y^!)} f^* \circ \varepsilon_{(p_{\xi \#},p_{\xi}^*)} \circ p_{\xi \#} \circ \varepsilon_{(s_{X!},s_X^!)} \circ s_{X!}(\alpha)\]
and 
\[th_E^Y(f^*\xi)(f^*\alpha) = \varepsilon_{(p_{f^*\xi \#},p_{f^*\xi}^*)}  \circ p_{f^*\xi \#} \varepsilon_{(s_{Y!},s_Y^!)} s_{Y!} \circ \varepsilon_{(s_{X*},s_X^!)} \circ Ex_*^*(-)^{-1} \circ \eta_{(s_{Y*},s_Y^!)} \circ f^*(\alpha).\]
The two composition are natural isomorphism to each other, as we have the natural $2$-isomorphisms: 
\[f^*s_{X!} \stackrel{\cong}{\longrightarrow} s_{Y!}f_{\xi}^*, \quad p_{f^*\xi \#}f_{\xi}^* \stackrel{\cong}{\longrightarrow} f^*p_{\xi \#}.\]
\end{proof}
Let $f: Y \rightarrow X$ be a smooth projective morphism of projective smooth $k$-schemes of relative dimension $d = \dim(Y) - \dim(X)$ and $p_{\xi}: \xi \rightarrow X$ be a vector bundle of rank $r$ with the zero section $s: X \rightarrow \xi$. We define in the following a pushforward on twisted $E$-cohomology: Consider 
\[SH(Y)[\mathbf{1}_Y,Th_Y(T_f)\wedge_Y s^!_{Y}p_{f^*\xi}^* E_Y^{2(r-d),(r-d)} \wedge S^{p,q}],\]
where $T_f$ is the normal bundle of the diagonal immersion $\delta: Y \rightarrow Y \times_X Y$.
\begin{multline*}
SH(Y)[\mathbf{1}_Y,Th_Y(T_f)\wedge_Y s^!_{Y}p^*_{f^*\xi}E_Y^{2(r-d),(r-d)} \wedge S^{p,q}] \cong \\
SH(Y)[Th_Y(f^*\xi) \wedge_Y Th_Y(-T_f),E_Y^{2(r-d),(r-d)} \wedge S^{p,q}],
\end{multline*}
where $Th_Y(-T_f) \in SH(Y)$ is the inverse of $Th_Y(T_f) \in SH(Y)$.  Since $E_Y = f^*E_X$, the adjunction 
\[f_{\#} : SH(Y) \leftrightarrows SH(X) : f^*\]
gives us a natural isomorphism 
\begin{multline*}
SH(Y)[Th_Y(f^*\xi) \wedge_Y Th_Y(-T_f),E_Y^{2(r-d),(r-d)} \wedge S^{p,q}] \cong \\ SH(X)[f_{\#}(Th_Y(f^*\xi) \wedge_Y Th_Y(-T_f)), E_X^{2(r-d),(r-d)}\wedge S^{p,q}]
\end{multline*}
By the projection formula $Pr^*_{\#}$ and since $Th_Y(f^*\xi) \cong f^*Th_X(\xi)$ as $\xi$ is an actual bundle, we have then a natural isomorphism 
\[f_{\#}(Th_Y(f^*\xi) \wedge_Y Th_Y(-T_f)) \cong Th_X(\xi) \wedge_X f_{\#}Th_Y(-T_f).\]
So we have then a natural isomorphism
\begin{multline*}
SH(X)[f_{\#}(Th_Y(f^*\xi) \wedge_Y Th_Y(-T_f)), E_X^{2(r-d),(r-d)}\wedge S^{p,q}] \cong \\ SH(X)[Th_X(\xi) \wedge_X f_{\#}Th_Y(-T_f),E_X^{2(r-d),(r-d)} \wedge S^{p,q}]. 
\end{multline*}
By \cite[Prop. 2.4.31]{CD10} we have 
\[f_{\#}Th_Y(-T_f) \cong D_X(f_{\#}\mathbf{1}_Y),\]
where $D_X(f_{\#}\mathbf{1}_Y)$ means the dual of $f_{\#}\mathbf{1}_Y$ in $SH(X)$. Hence there is a natural isomorphism
\begin{multline*}
SH(X)[Th_X(\xi) \wedge_X f_{\#}Th_Y(-T_f),E_X^{2(r-d),(r-d)} \wedge S^{p,q}] \cong \\ SH(X)[Th_X(\xi),f_{\#}\mathbf{1}_Y \wedge_X E_X^{2(r-d),(r-d)} \wedge S^{p,q}].
\end{multline*}
From the counit of the adjunction $(f_{\#},f^*)$ 
\[f_{\#}\mathbf{1}_Y \cong f_{\#}f^*\mathbf{1}_X \rightarrow \mathbf{1}_X\]
we have an induced map 
\[SH(X)[Th_X(\xi),f_{\#}\mathbf{1}_Y \wedge_X E_X^{2(r-d),(r-d)} \wedge S^{p,q}] \rightarrow SH(X)[Th_X(\xi),E_X^{2(r-d),(r-d)}\wedge S^{p,q}].\]
By the twisted Thom isomorphism, the later group is  
\[SH(X)[Th_X(\xi),E_X^{2(r-d),(r-d)}\wedge S^{p,q}] \cong E^{p-2d,q-d}(X,\xi).\] 
Now we define formally:
\begin{defn}{\rm
Let $f: Y\rightarrow X$ be a smooth projective morphism of projective smooth $k$-schemes of relative dimension $d = \dim(Y) - \dim(X)$ and $p_{\xi}: \xi \rightarrow X$ be a vector bundle. We define
\[E^{p,q}(Y,f^*\xi - T_f) \stackrel{defn}{=} SH(Y)[\mathbf{1}_Y,Th_Y(T_f) \wedge_Y s_Y^!p_{f^*_{\xi}}^*E_Y^{2(r-d),(r-d)}\wedge S^{p,q}], \]
where $\pi_{f^*\xi}: f^*\xi \rightarrow Y$ is the pullback bundle and $s_Y$ is its $0$-section. The pushforward map is the induced map constructed as above
\[E^{p,q}(Y,f^*\xi - T_f) \stackrel{f_*}{\rightarrow} E^{p-2d,q-d}(X,\xi).\]
}
\end{defn} 
Remark that our definition of projective pushforward $f_*$ doesn't require $E$ to be an oriented cohomology theory, however we need the assumption on smoothness of $f$. The reason that we choose the notation $E^{p,q}(Y,f^*\xi-T_f)$ is that this group should behave like the so-called cohomology twisted by formal difference of vector bundles. The shifting in the definition $(-2d,-d)$ reminds us that the inverse Thom transformation $Th_Y(-T_f)$ should behave like the Thom spectrum of the virtual bundle $-T_f$ after taking $\pi_{Y\#}$, where $\pi_Y: Y \rightarrow \Spec k$ is the structure morphism of $Y$, as the rank of the virtual bundle $-T_f$ is $-d$. As already mentioned in $\S 3$ we refer the reader to \cite[\S 4]{Rio10} and \cite{Del87} for the discussion on Picard category of virtual bundles. But we remind the reader again that we always work with an actual bundle $\xi$. 
\begin{rem}
Let $Z \stackrel{g}{\longrightarrow} Y \stackrel{f}{\longrightarrow} X$ be a sequence of composable morphisms. Then $f_* \circ g_*$ is not defined for a trivial reason: One has only a natural $2$-isomorphism 
\[e_{\sigma}^{\vee}: Th^{-1}_Z(s_Z,p_{T_{fg}}) \stackrel{\cong}{\longrightarrow} Th^{-1}_Z(s_Z,p_{T_g}) \circ g^*Th^{-1}_Y(s_Y,p_{T_g}),\]
which comes from the exact sequence 
\[0 \rightarrow g^*T_f \rightarrow T_{fg} \rightarrow T_g \rightarrow 0.\]
The $2$-isomorphism $e_{\sigma}^{\vee}$ is not an identity. Consequently $f_*\circ g_*$ is only defined up to this specific natural $2$-isomorphism. 
\end{rem}
\begin{rem} Another variant to construct pushforward can be obtained as follows: One has a natural isomorphism via Thom transformation adjunctions
Apply the functor $f_*: SH(Y) \rightarrow SH(X)$ we obtain a map 
\begin{multline*}
SH(Y)[\mathbf{1}_Y,Th_Y(T_f)\wedge_YTh(-f^*\xi)\wedge_YE_Y^{2(r-d),(r-d)} \wedge S^{p,q}] \stackrel{f_*(-)}{\longrightarrow} \\ SH(X)[f_*\mathbf{1}_Y,f_*(Th_Y(T_f)\wedge_Y Th_Y(-f^*\xi) \wedge E_Y^{2(r-d),(r-d)} \wedge S^{p,q})].
\end{multline*}
By projection formula $Pr_*^*(f)$ we have a canonical isomorphism 
\begin{multline*}
Ex^{*!}(f,s_X)\circ Ex_*^*(f,\wedge)^{-1}: f_*(Th_Y(T_f)) \wedge_Y Th_Y(-f^*\xi) \wedge_Y E_Y) \stackrel{\cong}{\longrightarrow} \\ f_*(Th_Y(T_f)) \wedge_X Th_X(-\xi) \wedge_X E_X,
\end{multline*}
which induces a canonical isomorphism 
\begin{multline*}
Ex^{*!}(f,s_X)\circ Ex^*_*(f,\wedge)^{-1} \circ: SH(X)[f_*\mathbf{1}_Y,f_*(Th_Y(T_f)\wedge_Y Th_Y(-f^*\xi) \wedge E_Y^{2(r-d),(r-d)} \wedge S^{p,q})] \\ \stackrel{\cong}{\longrightarrow} SH(X)[f_*\mathbf{1}_Y,f_*(Th_Y(T_f)) \wedge_X Th_X(-\xi) \wedge_X E_X^{2(r-d),(r-d)} \wedge S^{p,q}].  
\end{multline*}
From the unit $\eta_{(f^*,f_*)}(X): \mathbf{1}_X \rightarrow f_*f^*\mathbf{1}_X \cong f_*\mathbf{1}_Y$ of the adjunction $(f^*,f_*)$ one obtains a map 
\begin{multline*}
-\circ \eta_{(f^*,f_*)}(X):  SH(X)[f_*\mathbf{1}_Y,f_*(Th_Y(T_f)) \wedge_X Th_X(-\xi) \wedge_X E_X^{2(r-d),(r-d)} \wedge S^{p,q}] \longrightarrow \\ SH(X)[\mathbf{1}_X,f_*(Th_Y(T_f)) \wedge_X Th_X(-\xi) \wedge_X E_X^{2(r-d),(r-d)} \wedge S^{p,q}].
\end{multline*}
As $f$ is projective we have $f_* \cong f_!$ and since $f$ is smooth we have the canonical purity isomorphism 
\[\mathfrak{p}_f : f_\#(-) \stackrel{\cong}{\longrightarrow} f_*(Th_Y(T_f) \wedge_Y -),\]
which induces a canonical isomorphism
\begin{multline*}
-\circ \mathfrak{p}_f^{-1} : SH(X)[\mathbf{1}_X,f_*(Th_Y(T_f)) \wedge_X Th_X(-\xi) \wedge_X E_X^{2(r-d),(r-d)} \wedge S^{p,q}] \stackrel{\cong}{\longrightarrow} \\
SH(X)[\mathbf{1}_X,f_{\#}\mathbf{1}_Y \wedge_X Th_X(-\xi) \wedge_X E_X^{2(r-d),(r-d)} \wedge S^{p,q}].
\end{multline*}
The counit $\varepsilon_{(f_\#,f^*)}(X):  f_\#\mathbf{1}_Y \cong f_\#f^*\mathbf{1}_X \rightarrow \mathbf{1}_X$ of the adjunction $(f_\#,f^*)$ induces then a map 
\begin{multline*}
\varepsilon_{(f_\#,f^*)}(X) \circ : SH(X)[\mathbf{1}_X,f_{\#}\mathbf{1}_Y \wedge_X Th_X(-\xi) \wedge_X E_X^{2(r-d),(r-d)} \wedge S^{p,q}] \longrightarrow \\ SH(X)[\mathbf{1}_X,Th_X(\xi) \wedge_X E_X^{2(r-d),(r-d)} \wedge S^{p,q}] = E^{p-2d,q-d}(X,\xi). 
\end{multline*}
So we obtain a map defined as the composition of the maps above
\[SH(Y)[\mathbf{1}_Y,Th_Y(T_f) \wedge_Y s_Y^!p^*_{f^*\xi}E_Y^{2(r-d),(r-d)} \wedge S^{p,q}] \rightarrow E^{p-2d,q-d}(X,\xi).\]
One can check the two constructions are equivalent. And again the Proposition \ref{Ayoub2} tells us that the composition of pushforward maps $f_* \circ g_*$ is only defined up to a specific natural $2$-isomorphism. 
\end{rem}
\begin{prop}\label{functorialitypushforward}
Let $Z \stackrel{g}{\rightarrow} Y \stackrel{f}{\rightarrow} X$ be smooth projective morphisms of projective smooth $k$-schemes of relative dimension $e$ resp. $d$. Let $\xi/X$ be a vector bundle of rank $r$. Then one has up to a natural $2$-isomorphism 
\[(f\circ g)_* = f_* \circ g_*: E^{p,q}(Z, g^*f^* \xi -T_{fg}) \rightarrow E^{p-2(d+e),q-(d+e)}(X,\xi).\]
\end{prop}
\begin{proof}
We have an exact sequence of vector bundles on $Z$ \cite[17.2.3]{EGA4}
\[0 \rightarrow g^*T_f \rightarrow T_{fg} \rightarrow T_g \rightarrow 0.\]
So we have an isomorphism (cf. \cite[Rem. 2.4.52]{CD10})
\[e_{\sigma}: Th_Z(T_{fg}) \cong Th_Z(T_g) \wedge_Z Th_Z(g^*T_f) \cong Th_Z(T_g) \wedge_Z g^*Th_Y(T_f),\]
where $-\wedge_Z-$ means relative wedge product over $Z$. Since $g^*$ is strong monoidal and since all $f$ and $g$ are smooth, which means that the $\wedge_Z$-inverse object of $Th_Z(T_{fg})$ is $Th_Z(-T_{fg})$ and $Th_Z(T_g)^{-1} = Th_Z(-T_g)$ and $(g^*Th_Y(T_f))^{-1} = g^*(Th_Y(-T_f))$ (cf. \cite[Thm. 2.4.50 (3)]{CD10}). Hence we have
\[e_{\sigma}^{\vee}: Th_Z(-T_{fg}) \cong Th_Z(-T_g) \wedge_Z g^*Th_Y(-T_f).\]  
Functoriality of pushforward follows from this isomorphism as follow: We write $h = f \circ g$. Let $s_Z : Z \rightarrow g^*f^*\xi$ be the $0$-section of the vector bundle $p_{g^*f^*\xi}: g^*f^*\xi \rightarrow Z$. Let us recall the notation now: For an adjunction between categories
\[L : \mathcal{A} \leftrightarrows \mathcal{B} : R,\]
we denote
\[\varepsilon_{(L,R)}: LR \rightarrow \id, \quad \eta_{(L,R)}: \id \rightarrow RL\]
the counit and unit of the adjunction $(L,R)$ respectively. The composition $f_* \circ g_*$ is by construction the following composition:
\begin{multline*}
E^{p,q}(Z, g^*f^*\xi - T_{fg}) \stackrel{def}{=} SH(Z)[\mathbf{1}_Z,Th_Z(T_{fg}) \wedge_Z s_Z^!p^*_{g^*f^*\xi}E_Z^{2r-2(d+e),r-(d+e)} \wedge S^{p,q}] \stackrel{(1)}{\rightarrow} \\ SH(Z)[Th_Z(g^*f^* \xi)\wedge_Z Th_Z(-T_{fg}),E_Z^{2r-2(d+e),r-(d+e)} \wedge S^{p,q}] \stackrel{(2)}{\rightarrow} \\ SH(Y)[g_{\#}(Th_Z(g^* f^* \xi)\wedge_Z Th_Z(-T_{fg})),E_Y^{2r-2(d+e),r-(d+e)} \wedge S^{p,q}] \stackrel{(3)}{\rightarrow} \\ SH(Y)[Th_Y(f^*\xi)\wedge_Y g_{\#}Th_Z(-T_{fg}),E_Y^{2r-2(d+e),r-(d+e)} \wedge S^{p,q}] \stackrel{(4)}{\rightarrow} \\ SH(Y)[Th_Y(f^*\xi)\wedge_Yg_{\#}(Th_Z(-T_g) \wedge_Z g^*Th_Y(-T_f)),E_Y^{2r-2(d+e),r-(d+e)} \wedge S^{p,q}] \stackrel{(5)}{\rightarrow} \\ SH(Y)[Th_Y(f^*\xi)\wedge_Y Th_Y(-T_f) \wedge_Y g_{\#}Th_Z(-T_g),E_Y^{2r-2(d+e),r-(d+e)} \wedge S^{p,q}] \stackrel{(6)}{\rightarrow} \\ SH(Y)[Th_Y(f^*\xi) \wedge_Y Th_Y(-T_f)\wedge D_Y(g_{\#}\mathbf{1}_Z),E_Y^{2r-2(d+e),r-(d+e)} \wedge S^{p,q}] \stackrel{(7)}{\rightarrow} \\ SH(Y)[Th_Y(f^*\xi) \wedge_Y Th_Y(-T_f),g_{\#}(\mathbf{1}_Z)\wedge_Y E_Y^{2r-2(d+e),r-(d+e)} \wedge S^{p,q}] \stackrel{(8)}{\rightarrow} \\ SH(Y)[Th_Y(f^*\xi)\wedge_Y Th_Y(-T_f),E_Y^{2r-2(d+e),r-(d+e)} \wedge S^{p,q}] \stackrel{(9)}{\rightarrow} \\ SH(X)[f_{\#}(Th_Y(f^*\xi)\wedge_Y Th_Y(-T_f)),E_X^{2r-2(d+e),r-(d+e)} \wedge S^{p,q}] \stackrel{(10)}{\rightarrow} \\ SH(X)[Th_X(\xi) \wedge_X f_{\#}Th_Y(-T_f),E_X^{2r-2(d+e),r-(d+e)} \wedge S^{p,q}] \stackrel{(11)}{\rightarrow} \\ SH(X)[Th_X(\xi) \wedge_X D_X(f_{\#}\mathbf{1}_Y),E_X^{2r-2(d+e),r-(d+e)} \wedge S^{p,q}] \stackrel{(12)}{\rightarrow}\\ SH(X)[Th_X(\xi),f_{\#}(\mathbf{1}_Y) \wedge_X E_X^{2r-2(d+e),r-(d+e)} \wedge S^{p,q}] \stackrel{(13)}{\rightarrow}\\ SH(X)[Th_X(\xi),E_X^{2r-2(d+e),r-(d+e)} \wedge S^{p,q}] \cong E^{p-2(d+e),q-(d+e)}(X,\xi),
\end{multline*}
where $(1)$ is the natural isomorphism given by the adjunction of the Thom transformations $Th(s_Z,p_{g^*f^*\xi})$ and $Th_Z(T_{fg})$, $(2)$ is the natural isomorphism given by the adjunction $(g_{\#},g^*)$
\[(2)(-) = \varepsilon_{(g_\#,g_*)}(-)\circ g_{\#}(-),\]
$(3)$ is the natural isomorphism given by the projection formula $Pr^*_{\#}(g)$
\[(3)(-) = \varepsilon_{(g_{\#},g^*)}(-)\circ g_{\#}(\eta_{(g_{\#},g^*}(-) \wedge_Z \id)(-),\]
$(4)$ is the natural isomorphism given by 
\[e_{\sigma}^{\vee}: Th_Z(-T_{fg}) \cong Th_Z(-T_g)\wedge_Z g^*Th_Y(-T_f),\]
$(5)$ is the natural isomorphism given by the projection formula $Pr_{\#}^*(g)$, $(6)$ is the natural isomorphism given by duality in $SH(Y)$:
\[g_{\#}Th_Z(-T_g) \cong D_Y(g_{\#}\mathbf{1}_Z),\]
$(7)$ is the natural isomorphism given by adjunction of duality in $SH(Y)$
\[(7)(-) = (\id_{g_{\#}\mathbf{1}_Z} \wedge - ) \circ (coev_{D_Y(g_{\#}\mathbf{1}_Z)} \wedge \id_{-}),\]
$(8)$ is the pushforward induced by the counit $g_{\#}\mathbf{1}_Z \cong g_{\#}g^*\mathbf{1}_Y \rightarrow \mathbf{1}_Y$ 
\[(8) = \varepsilon_{(g_{\#},g^*)}(-) \wedge_Y - ,\]
$(9)$ is the natural isomorphism given by the adjunction $(f_{\#},f^*)$ 
\[(9) = \varepsilon_{(f_{\#},f^*)}(-) \circ f_{\#}(-),\]
$(10)$ is the natural isomorphism given by the projection formula $Pr_{\#}^*(f)$ 
\[(10)(-) = \varepsilon_{(f_{\#},f^*)}(-) \circ f_{\#}(\eta_{(f_{\#},f^*)} \wedge_Y \id)(-),\]
$(11)$ is the natural isomorphism given by duality in $SH(X)$:
\[f_{\#}Th_Y(-T_f) \cong D_X(f_{\#}\mathbf{1}_Y),\]
$(12)$ is the natural isomorphism given by the adjunction of duality in $SH(X)$
\[(12)(-) = (\id_{f_{\#}\mathbf{1}_Y} \wedge - ) \circ (coev_{D_X(f_{\#}\mathbf{1}_Y)} \wedge \id_{-}),\]
and finally $(13)$ is the pushforward induced by the counit $f_{\#}\mathbf{1}_Y \cong f_{\#}f^*\mathbf{1}_X \rightarrow \mathbf{1}_X$:
\[(13)(-) = \varepsilon_{(f_{\#},f^*)}(-) \wedge_X - . \]  
$h_* = (f \circ g)_*$ is the following composition: 
\begin{multline*}
E^{p,q}(Z,g^*f^*\xi) \stackrel{defn}{=} SH(Z)[\mathbf{1}_Z,Th_Z(T_{fg})\wedge_Z s_Z^!p^*_{g^*f^*\xi}E_Z^{2r-2(d+e),r-(d+e)} \wedge S^{p,q}] \stackrel{(1')}{\rightarrow} \\ SH(Z)[Th_Z(g^*f^*\xi) \wedge_Z Th_Z(-T_{fg}),E_Z^{2r-2(d+e),r-(d+e)} \wedge S^{p,q}] \stackrel{(2')}{\rightarrow} \\ SH(X)[h_{\#}(Th_Z(g^*f^*\xi)\wedge_Z Th_Z(-T_{fg})),E_X^{2r-2(d+e),r-(d+e)} \wedge S^{p,q}] \stackrel{(3')}{\rightarrow} \\ SH(X)[Th_X(\xi) \wedge_X h_{\#}Th_Z(-T_{fg}),E_X^{2r-2(d+e),r-(d+e)} \wedge S^{p,q}]\stackrel{(4')}{\rightarrow} \\ SH(X)[Th_X(\xi) \wedge_X D_X(h_{\#}\mathbf{1}_Z),E_X^{2r-2(d+e),r-(d+e)} \wedge S^{p,q}] \stackrel{(5')}{\rightarrow} \\ SH(X)[Th_X(\xi),h_{\#}(\mathbf{1}_Z) \wedge_X E_X^{2r-2(d+e),r-(d+e)} \wedge S^{p,q}] \stackrel{(6')}{\rightarrow} \\ SH(X)[Th_X(\xi),E_X^{2r-2(d+e),r-(d+e)} \wedge S^{p,q}] \cong E^{p-2(d+e),q-(d+e)}(X,\xi),
\end{multline*}
where $(1')$ is the natural isomorphism given by the adjunction of the Thom transformations $Th(s_Z,p_{g^*f^*\xi})$ and $Th_Z(T_{fg})$, $(2')$ is the natural isomorphism given by the adjunction $(h_{\#},h^*)$ 
\[(2')(-) = \varepsilon_{(h_{\#},h^*)}(-) \circ h_{\#}(-),\]
$(3')$ is the natural isomorphism given by the projection formula $Pr_{\#}^*(h)$ 
\[(3')(-) = \varepsilon_{(h_{\#},h^*)}(-) \circ h_{\#}(\eta_{(h_{\#},h^*)}(-) \wedge_Z \id )(-),\]
$(4')$ is the natural isomorphism by duality in $SH(X)$:
\[h_{\#}Th_Z(-T_{fg}) \cong D_X(h_{\#}\mathbf{1}_Z),\]
$(5')$ is the natural isomorphism given by the adjunction of duality in $SH(X)$:
\[(5')(-) = (\id_{h_{\#}\mathbf{1}_Z} \wedge_X - ) \circ (coev_{D_X(h_{\#}\mathbf{1}_Z)} \wedge \id_{-} )  \]
and finally $(6')$ is the pushforward induced by the counit $h_{\#}\mathbf{1}_Z \cong h_{\#}h^*\mathbf{1}_X \rightarrow \mathbf{1}_X$: 
\[(6')(-) = \varepsilon_{(h_{\#},h^*)}(-) \wedge_X - .\]
The maps $(1)$ and $(1')$ are identical. The diagram
\[ \xymatrix{ & \bullet \ar[ld]_{(2)} \ar[rd]^{(2')} \\ \bullet \ar[rr]_{(a_1)} && \bullet   }\]
commutes, because $h_{\#} \cong f_{\#} \circ g_{\#}$, where 
\begin{multline*}
(a_1): SH(Y)[g_{\#}(Th_Z(h^*\xi)\wedge_ZTh_Z(-T_{fg})),E_Y^{2r-2(d+e),r-(d+e)} \wedge S^{p,q}] \stackrel{\cong}{\rightarrow} \\ SH(X)[h_{\#}(Th_Z(h^*\xi)\wedge_ZTh_Z(-T_{fg})),E_X^{2r-2(d+e),r-(d+e)} \wedge S^{p,q}]
\end{multline*}
is the natural isomorphism induced from the adjunction $(f_{\#},f^*)$. Indeed, let 
\[\alpha_1: Th_Z(h^*\xi)\wedge_Z Th(-T_{fg}) \rightarrow E_Z^{2r-2(d+e),r-(d+e)} \wedge S^{p,q} \]
be a morphism. Then one has
\[(2')(\alpha_1) = \varepsilon_{(h_{\#},h^*)}(-) \circ h_{\#}(\alpha_1),\]
and
\[(a_1)\circ (2)(\alpha_1) = (\varepsilon_{(f_{\#},f^*)} \circ f_{\#}(-)) \circ (\varepsilon_{(g_{\#},g^*)} \circ g_{\#}(-))(\alpha_1). \]
So we have $(2') = (a_1) \circ (2)$.
Consider the pentagon
\[\xymatrix{ & \bullet \ar[rr]^{(a_1)} \ar[ld]_{(3)} && \bullet \ar[rd]^{(3')} \\ \bullet \ar[rrd]_{(a_2)} &&&& \bullet  \\ && \bullet \ar[rru]_{(a_3)} }\]
where 
\begin{multline*}
(a_2): SH(Y)[Th_Y(f^*\xi)\wedge_Y g_{\#}Th_Z(-T_{fg}),E_Y^{2r-2(d+e),r-(d+e)} \wedge S^{p,q}] \stackrel{\cong}{\rightarrow} \\ SH(X)[f_{\#}(Th_Y(f^*\xi)\wedge_Y g_{\#}Th_Z(-T_{fg})),E_X^{2r-2(d+e),r-(d+e)} \wedge S^{p,q}]
\end{multline*}
is the natural isomorphism induced by the adjunction $(f_{\#},f^*)$ and 
\begin{multline*}
(a_3): SH(X)[f_{\#}(Th_Y(f^*\xi)\wedge_Y g_{\#}Th_Z(-T_{fg})),E_X^{2r-2(d+e),r-(d+e)} \wedge S^{p,q}] \stackrel{\cong}{\rightarrow} \\ SH(X)[Th_X(\xi) \wedge_X h_{\#}Th(-T_{fg}),E_X^{2r-2(d+e),r-(d+e)} \wedge S^{p,q}]
\end{multline*}
is the natural isomorphism given by the projection formula $Pr_{\#}^*(f)$. We remind the reader that the isomorphism in the projection formula $Pr_{\#}^*(f)$ is given by the composition: 
\[f_{\#}(M\wedge_Y f^*N) \rightarrow f_{\#}(f^*f_{\#}(M) \wedge_Y f^*N) \simeq f_{\#}f^*(f_{\#}M\wedge_X N) \rightarrow f_{\#} \wedge_X N.\]
The pentagon commutes since isomorphism induced by the projection formula $Pr_{\#}^*(h)$ is the composing of isomorphisms coming from projection formulas $Pr_{\#}^*(g)$ and $Pr_{\#}^*(f)$. Indeed, let 
\[\alpha_2: g_{\#}(Th_Z(h^*\xi)\wedge_ZTh_Z(-T_{fg})) \rightarrow E_Y^{2r-2(d+e),r-(d+e)} \wedge S^{p,q}\]
be a morphism. We have 
\begin{multline*}
(3')(-) \circ (a_1) (\alpha_2) = (3')(-)(\varepsilon_{(f_{\#},f^*)} \circ f_{\#}(\alpha_2)) = Pr_{\#}^*(h)(-)(\varepsilon_{(f_{\#},f^*)} \circ f_{\#}(\alpha_2)) \\ = \varepsilon_{(h_{\#},h^*)}(-) \circ h_{\#}(\eta_{(f_{\#},f^*)}(-) \wedge_Z \id)(\varepsilon_{(f_{\#},f^*)} \circ f_{\#}(\alpha_2)),
\end{multline*}
and
\begin{multline*}
(a3)(-)\circ (a2)(-) \circ (3)(\alpha_2) = (a_3)(-)\circ (a_2)(-) \circ Pr_{\#}^*(g)(\alpha_2) = \\ (a3)(-)\circ (a2)(-) \circ \varepsilon_{(g_{\#},g^*)}(-) \circ g_{\#}(\eta_{g_{\#},g^*)}(-) \wedge_Z \id)(\alpha_2) \\ = (a3)(-) \circ \varepsilon_{(f_{\#},f^*)}(-) \circ f_{\#}(-) \circ \varepsilon_{(g_{\#},g^*)}(-) \circ g_{\#}(\eta_{g_{\#},g^*)}(-) \wedge_Z \id)(\alpha_2) \\ = \varepsilon_{(f_{\#},f^*)}(-) \circ f_{\#}(-)(\eta_{(f_{\#},f^*)}(-) \wedge_Y \id) \circ \varepsilon_{(f_{\#},f^*)}(-) \circ f_{\#}(-) \circ \varepsilon_{(g_{\#},g^*)}(-) \circ g_{\#}(\eta_{g_{\#},g^*)}(-)\wedge_Z \id)(\alpha_2).  
\end{multline*}
So we have 
\[(a3)(-) \circ (a2)(-) \circ(3)(\alpha_2) = (3')(-)\circ(a_1)(\alpha_2).\]
For any $K \in SH(Z)$ one has commutative diagram (see \cite[\S 1.4.2, \S 1.5]{Ay08} and \cite[Rem. 2.4.52]{CD10})
\[\xymatrix{h_{\#}K \ar@{=}[rrrr] \ar[ddd] &&&& f_{\#}g_{\#}K \ar[d]^{\cong} \\ &&&& f_!(Th_Y(T_f)\wedge_Yg_!(Th_Z(T_g)\wedge_Z K) \ar[d]^{\cong} \\ &&&& f_!g_!(g^*Th_Y(T_f) \wedge_Z Th_Z(T_g)\wedge_Z K) \ar[d]_{e_{\sigma}}^{\cong} \\ h_!(Th_Z(T_{fg})\wedge_Z K) \ar@{=}[rrrr] &&&& f_!g_!(Th_Z(T_{fg})\wedge_Z K) }\]
Now we take $K = \mathbf{1}_Z$ and dualize $D_X(-)$ the commutative diagram above. One has 
\[D_X(h_!(Th_Z(T_{fg})) \cong D_X(h_{\#}\mathbf{1}_Z) \cong h_{\#}Th_Z(-T_{fg}),\]
\[D_X(f_!Th_Y(T_f)) \cong D_X(f_{\#}\mathbf{1}_Y) \cong f_{\#}Th_Y(-T_f),\]
and 
\[D_Y(g_!Th_Z(T_g)) = D_Y(g_{\#}\mathbf{1}_Z) \cong g_{\#}Th_Z(-T_g).\]
So we can conclude that 
\[(6')\circ (5') \circ (4') \circ (a_3) \circ (a_2) = (13)\circ (12) \circ \cdots \circ (5) \circ(4).\]
Putting all together we have 
\[(6')\circ \cdots \circ (1') = (13)\circ \cdots \circ (1), \]
which means the pushforward on twisted $E$-cohomology satisfies $(f\circ g)_* = f_* \circ g_*$ up to a natural isomorphism induced by the natural $2$-isomorphism
\[e_{\sigma}^{\vee}: Th_Z^{-1}(s_Z,p_{T_{fg}}) \stackrel{\cong}{\longrightarrow} Th_Z^{-1}(s_Z,p_{T_g})\circ g^*Th_Y^{-1}(s_Y,p_{T_g}).\]
\end{proof}
Now we follow a suggestion by M. Levine to make a refinement to the result of Voevodsky in \ref{ThmVoeThom}, since as pointed out by M. Levine it is not enough to use the identities in $K_0(-)$ to constructs maps between twisted $\mathbf{E}$-cohomology groups.  
\begin{prop}(A refinement of Voevodsky's theorem) \label{refinement1}
Let $X \in SmProj(k)$ of dimension $d_X$, where $k$ is a field. After fixing an embedding $X \inj \P^d$ there exists a vector bundle $V_X$ on $X$ of rank $d^2+2d-d_X$, such that one has a specific isomorphism between objects in the Picard category of virtual bundles $\sV(X)$ on $X$: 
\[V_X \oplus T_X \cong \sO_X^{d^2+2d}.\]
\end{prop}
\begin{proof}
Case 1: $X = \P^d$. One has an exact sequence 
\[0 \rightarrow \sO_{\P^d} \rightarrow \sO_{\P^d}(1)^{\oplus (d+1)} \rightarrow T_{\P^d} \rightarrow 0.\]
By taking dual one also has 
\[0 \rightarrow \Omega_{\P^d} \rightarrow \sO_{\P^d}(-1)^{\oplus (d+1)} \rightarrow \sO_{\P^d} \rightarrow 0.\]
There are two isomorphisms between objects in $\sV(X)$:
\[\sO_{\P^d} \oplus T_{\P^d} \cong \sO_{\P^d}(1)^{\oplus (d+1)}\]
and 
\[\Omega_{\P^d} \oplus \sO_{\P^d} \cong \sO_{\P^d}(-1)^{\oplus (d+1)}.\] 
Define 
\[V_{\P^d} \stackrel{defn}{=} \Omega_{\P^d} \oplus (\Omega_{\P^d} \otimes T_{\P^d}).\]
As the Picard category $\sV(X) = V(Vect(X))$ (the category of virtual objects associated to the exact category of vector bundles on $X$) has not just $\oplus$, but also a biexact functor 
\[-\otimes- : \sV(X) \times \sV(X) \rightarrow \sV(X),\]
which is distributive (\cite{Del87}), one has an isomorphism in $\sV(X)$: 
\[(\Omega_{\P^d} \oplus \sO_{\P^d}) \otimes (\sO_{\P^d} \oplus T_{\P^d}) \cong \Omega_{\P^d} \oplus \sO_{\P^d} \oplus (\Omega_{\P^d} \otimes T_{\P^d}) \oplus T_{\P^d} \cong \sO_{\P^d}^{\oplus (d^2 + 2d + 1)}.\]
This implies that we have an isomorphism in $\sV(X)$:
\[V_{\P^d} \oplus T_{\P^d} \cong \sO_{\P^d}^{d^2+2d}.\]

Case 2: $X$ is smooth projective. Let $i : X \inj \P^d$ be a closed embedding. One define 
\[V_X \stackrel{defn}{=} N_{X/\P^d} \oplus i^*(V_{\P^d}),\]
where $N_{X/\P^d}$ denotes the normal bundle of $X$ in $\P^d$. 
In $\sV(X)$ one has an isomorphism between objects 
\[N_{X/\P^d} \oplus i^*(V_{\P^d} \oplus T_{\P^d}) \cong N_{X/\P^d} \oplus \sO_X^{d^2 + 2d}.\]
From the exact sequence 
\[0 \rightarrow T_X \rightarrow i^*T_{\P^d} \rightarrow N_{X/\P^d} \rightarrow 0\]
one has an isomorphism in $\sV(X)$:
\[T_X \oplus N_{X/\P^d} \cong i^*T_{\P^d}.\]
This implies that we have a isomorphism in $\sV(X)$:
\[N_{X/\P^d} \oplus i^*V_{\P^d} \oplus T_X \oplus N_{X/\P^d} \cong N_{X/\P^d} \oplus \sO_X^{d^2 + 2d}.\]  
This implies that we have a specific isomorphism in $\sV(X)$: 
\[V_X \oplus T_X \cong \sO_{X}^{d^2+2d}.\]
\end{proof}
P. Hu in \cite{Hu05} didn't check if her construction is the same as the construction of Voevodsky. We notice that the map constructed by Voevodsky \cite[Thm. 2.11]{Voe03}
\[T^{\wedge n_X+d_X} \rightarrow Th(V_X)\]
is first of all only in $\mathbf{Ho}_{\A^1,+}(k)$ and secondly very difficult to follow. We will take the refinement $V_X \oplus T_X \stackrel{\cong}{\longrightarrow} \sO_X^{d^2+2d}$ in $\sV(X)$ and construct the Pontryagin-Thom collapse map  
\[PTV: \mathbb{S}^0 \rightarrow \Sigma^{\infty}_{T,+}Th(V_X) \wedge S^{-2(d^2+2d),-(d^2+2d)}\]
by unpacking Hu's construction. Firstly, for a projective smooth $k$-variety $i : X \inj \P^d$, we have by definition 
\[V_X = N_{X/\P^d} \oplus i^*V_{\P^d} \cong N_{X/V_{\P^d}}.\] 
If $PTV$ is already for $\P^d$ constructed, then $PTV$ for $X$ is defined by the composition 
\begin{multline*}
\mathbb{S}^0 \longrightarrow \Sigma^{\infty}_{T,+} Th(V_{\P^d}) \wedge S^{-2(d^2+2d),-(d^2+2d)} \stackrel{q}{\longrightarrow} \Sigma^{\infty}_{T,+} Th(V_X) \wedge S^{-2(d^2+2d),-(d^2+2d)},
\end{multline*}
where $q$ is the quotient map 
\[Th(V_{\P^d}) \rightarrow V_{\P^d}/(V_{\P^d} - X) \stackrel{\cong}{\longrightarrow} Th(N_{X/V_{\P^d}}) = Th(V_X).\]
The isomorphism $V_{\P^d}/(V_{\P^d}-X) \stackrel{\cong}{\longrightarrow} Th(N_{X/V_{\P^d}})$ is the homotopy purity isomorphism (\cite[\S 3 Thm.2.23]{MV01}). For $X = \P^1$ one has a commutative diagram in $\mathbf{Ho}_{\A^1,+}(k)$ (\cite[pp. 9]{Hu05}) 
\[\xymatrix{((X \times X)-\Delta_X)_+ \ar[rr] \ar[d]_{pr_1}^{\cong} && (X \times X)_+ \ar[r]^{g} & Th(T_X) \\ X_+ \ar[rru]_{f}}\]
because $pr_1: (X \times X)- \Delta_X \rightarrow X$ is an affine bundle. So one has a cofiber sequence in $\mathbf{Ho}_{\A^1,+}(k)$
\[X_+ \stackrel{f}{\longrightarrow} (X\times X)_+ \stackrel{g}{\longrightarrow} Th(T_X).\]
For a vector bundle $\xi$ on $X$ one has 
\[Th(pr_1^*\xi/X \times X) = Th(\xi/X) \wedge X_+.\]
One the other hand one has commutative diagram (\cite[(3.13)]{Hu05})
\[\xymatrix{Th(pr_1^*\xi/X\times X - \Delta_X) \ar[d]^{\cong}_{pr_1} \ar[r] & Th(pr_1^*\xi/X\times X) \ar[d]^{pr_1} \\ Th(\xi/X) \ar@{=}[r] & Th(\xi/X)  }\]
So one obtains a cofiber sequence in $\mathbf{Ho}_{\A^1,+}(k)$
\[Th(\xi/X) \stackrel{f^{\xi}}{\longrightarrow} Th(\xi/X) \wedge X_+ \stackrel{g^{\xi}}{\longrightarrow} Th(T_X \oplus \xi).\]
Now we take $\xi = V_X$ and by the refinement $V_X \oplus T_X \stackrel{\cong}{\longrightarrow} \sO_X^{d^2+2d}$ in $\sV(X)$ we have then a cofiber sequence 
\[Th(V_X) \stackrel{f^{V_X}}{\longrightarrow} Th(V_X) \wedge X_+ \stackrel{g^{V_X}}{\longrightarrow} Th(V_X\oplus T_X) \cong S^{2(d^2+2d),(d^2+2d)} \wedge X_+.\]
This gives rise to a map in $SH(k)$:
\[\varepsilon: \Sigma^{\infty}_{T,+}Th(V_X) \wedge X_+ \wedge S^{-2(d^2+2d),-(d^2+2d)} \stackrel{g^{V_X}}{\longrightarrow} \Sigma^{\infty}_{T,+}X \rightarrow \mathbb{S}^0.\]
To construct the Pontryagin-Thom collapse map 
\[PTV: \mathbb{S}^0 \rightarrow \Sigma^{\infty}_{T,+}Th(V_X) \wedge S^{-2(d^2+2d),-(d^2+2d)},\]
such that the composition 
\begin{multline*}
g^{V_X} \circ (PTV \wedge \id) : \mathbb{S}^0 \wedge \Sigma^{\infty}_{T,+}X  \stackrel{PTV \wedge \id}{\longrightarrow} \Sigma^{\infty}_{T,+}Th(V_X) \wedge X_+ \wedge S^{-2(d^2+2d),-(d^2+2d)} \\ \stackrel{g^{V_X}}{\longrightarrow} \Sigma^{\infty}_{T,+}X
\end{multline*}
is the identity $\id_{\Sigma^{\infty}_{T,+}X}$ in $SH(k)$, it is enough to construct a map 
\[PTV: \mathbb{S}^0 \rightarrow \Sigma^{\infty}_{T,+}Th(V_X) \wedge S^{-2(d^2+2d),-(d^2+2d)},\]
such that the composition $\varepsilon \circ (PTV \wedge \id)$ is the collapse map $\Sigma^{\infty}_{T,+}X \rightarrow \mathbb{S}^0$, because $g^{V_X}$ is the composition 
\begin{multline*}
\Sigma^{\infty}_{T,+}Th(V_X) \wedge X_+ \wedge S^{-2(d^2+2d),-(d^2+2d)} \stackrel{\id \wedge \Delta}{\longrightarrow} \Sigma^{\infty}_{T,+}Th(V_X) \wedge X_+ \wedge X_+ \wedge S^{-2(d^2+2d),-(d^2+2d)} \\ \stackrel{\varepsilon \wedge \id}{\longrightarrow} \Sigma^{\infty}_{T,+}X 
\end{multline*}
By adjunction $\varepsilon$ gives us a map 
\[\lambda_X : \Sigma^{\infty}_{T,+}Th(V_X) \wedge S^{-2(d^2+2d),-(d^2+2d)} \rightarrow D_k(\Sigma^{\infty}_{T,+}(X)) = \underline{\Hom}(\Sigma^{\infty}_{T,+}X,\mathbb{S}^0).\]
We remind the reader that P. Hu started with $X = \P^1$ before \cite[Lem. 3.8]{Hu05}, since she wanted to prove some particular results for projective quadric. For general $X$ and a vector bundle $\xi$ on $X$ one still has the map 
\[g: (X\times X)_+ \rightarrow Th(T_X) = Th(N_{X/X\times X})\]
and hence a map
\[g^{\xi}: Th(pr_1^*\xi/X\times X) = Th(\xi/X) \wedge X_+ \rightarrow Th(T_X \oplus \xi)\]
and hence by applying $\xi = V_X$ one has a map 
\[\lambda_X : \Sigma^{\infty}_{T,+}Th(V_X) \wedge S^{-2(d^2+2d),-(d^2+2d)} \rightarrow \underline{\Hom}(\Sigma^{\infty}_{T,+}X, \mathbb{S}^0).\]
Now we consider the linear embedding $i: \P^{d} \inj \P^{d+1}$. By construction 
\[V_{\P^{d+1}} = N_{\P^{d}/\P^{d+1}} \oplus i^*(V_{\P^d})\] and so the diagram 
\begin{equation} \label{eqhu1}
\xymatrix{\Sigma^{\infty}_{T,+}Th(V_{\P^{d+1}}) \wedge S^{-2((d+1)^2+2(d+1)),-((d+1)^2+2(d+1))} \ar[d]_{q^{\P^{d+1}}_{\P^d}} \ar[r]^{\quad \quad \quad \quad \quad \quad \quad \quad \quad \lambda_{\P^{d+1}}} & D(\P^{d+1}_+) \ar[d]^{D(i)} \\ \Sigma^{\infty}_{T,+}Th(V_{\P^d}) \wedge S^{-2(d^2+2d),-(d^2+2d)} \ar[r]_{\quad \quad \quad \quad \quad \quad \quad \lambda_{\P^d}} & D(\P^d_+) }
\end{equation}
commutes, since it is adjoint to the commutativity of the diagram 
\[\xymatrix{ Th(V_{\P^{d+1}}) \wedge \P^d_+ \wedge S^{\star,\star} \ar[d]_{q^{\P^{d+1}}_{\P^d} \wedge \id} \ar[r]^{\id \wedge i_+} & Th(V_{\P^{d+1}}) \wedge \P^{d+1} \wedge S^{\star,\star} \ar[d]^{g^{V_{\P^{d+1}}}} \\ Th(V_{\P^d}) \wedge \P^d_+ \wedge S^{*,*} \ar[r]_{\quad \quad \quad i_+ \circ g^{V_{\P^d}}} & \P^{d+1}_+   }\]
where we write $q^{\P^{d+1}}_{\P^d}$ for the quotient map and the last commutative diagram is obtained by $V_{\P^{d+1}}$-Thomification (i.e. we apply $g^{V_{\P^{d+1}}}$ on $\P^{d+1} \times \P^{d+1} \rightarrow Th(T_{\P^{d+1}})$) of the commutative diagram 
\[\xymatrix{(\P^{d+1}\times \P^d)_+ \ar[r]^i \ar[d] & (\P^{d+1} \times \P^{d+1})_+ \ar[dd]  \\ (\P^{d+1}\times \P^d)_+/((\P^{d+1}\times \P^d)-(\P^d \times \P^d))_+ \ar[d] \\ (\P^{d+1} \times \P^d)_+/((\P^{d+1}\times \P^d)-\Delta_{\P^d})_+ \ar[r] & (\P^{d+1} \times \P^{d+1})_+/((\P^{d+1}\times \P^{d+1})-\Delta_{\P^{d+1}})_+    }\]
Consider the composition 
\[(\P^{d+1}-\P^d) \times \P^{d+1} \rightarrow \P^{d+1} \times \P^{d+1} \rightarrow (\P^{d+1}\times \P^{d+1})/(\P^{d+1} \times \P^{d+1} - \Delta_{\P^{d+1}}). \] 
$(\P^{d+1}-\P^d)\times \P^d$ is mapped to $(\P^{d+1}\times \P^{d+1})- \Delta_{\P^{d+1}}$. So the composition above induces a map
\[Th(j^*V_{\P^{d+1}}) \wedge (\P^{d+1}/\P^d)_+ \wedge S^{\star,\star} \rightarrow \P^{d+1}_+,\]
where $j: (\P^{d+1}-\P^d) \inj \P^{d+1}$ denotes the open immersion. After composing with the collapse map $\P^{d+1}_+ \rightarrow \mathbb{S}^0$ and taking adjoint one obtains a map 
\[\lambda_{\P^{d+1}/\P^d}: Th(j^*V_{\P^{d+1}}) \wedge S^{\star,\star} \rightarrow D((\P^{d+1}/\P^d)_+).\]
By construction there is a commutative diagram in $SH(k)$:
\begin{equation}\label{eqhu2}
\xymatrix{Th(j^*V_{\P^{d+1}}) \wedge S^{\star,\star} \ar[rr]^{\lambda_{\P^{d+1}/\P^d}} \ar[d]_{Th(j)} && D((\P^{d+1}/\P^d)_+) \ar[d]^{D(p)} \\ Th(V_{\P^{d+1}}) \wedge S^{\star,\star} \ar[rr]_{\lambda_{\P^{d+1}}} && D(\P^{d+1}_+) }
\end{equation}  
Now the Claim $1$ in the proof of \cite[Lem. 3.8]{Hu05} implies that there is a morphism of distinguished triangles given by the commutative diagrams \ref{eqhu1} and \ref{eqhu2}. 
\[\xymatrix{ S^{-2d,-d} \ar[d] \ar[r]^{\cong} & S^{-2d,-d} \ar[d] \\ Th(V_{\P^{d+1}}) \wedge S^{\star,\star} \ar[r]^{\quad \quad  \lambda_{\P^{d+1}}} \ar[d] & D(\P^{d+1}_+) \ar[d] \\ Th(V_{\P^d}) \wedge S^{*,*} \ar[r]_{\quad \quad \quad \lambda_{\P^d}} & D(\P^d_+)  }\]
So by induction on $d$ one can conclude that 
\[\lambda_{\P^d} : Th(V_{\P^d}) \wedge S^{-2(d^2+2d),-(d^2+2d)} \rightarrow D(\P^d_+)\]
is an isomorphism in $SH(k)$ for all $d \geq 0$. Now we can construct the Pontryagin-Thom collapse map 
\[PTV: \mathbb{S}^0 \rightarrow D(\P^d_+) \stackrel{\lambda_{\P^d}^{-1}}{\longrightarrow} Th(V_{\P^d}) \wedge S^{-2(d^2+2d),-(d^2+2d)}.\] 
If $X \inj \P^d$ is a smooth projective $k$-variety, we define $PTV$ for $X$ as the composition of $PTV$ for $\P^d$ with the quotient map 
\[Th(V_d) \wedge S^{-2(d^2+2d),-(d^2+2d)} \rightarrow Th(V_X) \wedge S^{-2(d^2+2d),-(d^2+2d)}.\]
So by construction we have:
\begin{prop}\label{refinement2}
Let $X$ be a smooth projective $k$-scheme. After fixing an embedding $i : X \inj \P^d$, there is a commutative diagram in $SH(k)$
\[\xymatrix{\mathbb{S}^0 \ar[rr]^{PTV \quad \quad \quad \quad \quad \quad} \ar[rrd]_{PTH} && \Sigma^{\infty}_{T,+}Th(V_X) \wedge S^{-2(d^2+2d),-(d^2+2d)} \ar[d]^{\cong} \\ && \Sigma^{\infty}_{T,+} Th(-T_X) }\]
where $PTH: \mathbb{S}^0 \rightarrow \Sigma^{\infty}_{T,+}Th(-T_X)$ is the map constructed in \cite[Lem. 3.18]{Hu05}.
\end{prop}
\begin{prop}\label{isothomtiwst}
Let $f: Y \rightarrow X$ be a projective smooth morphism of projective smooth $k$-schemes of relative dimension $d = d_Y - d_X$. After fixing an embedding $X \inj \P^N$, there is an isomorphism 
\[tt_E^Y: E^{p,q}(Y,f^*V_X - T_f) \cong E^{p+2n_Y,q+n_Y}(Th(V_Y)), \]
where $V_X$ and $V_Y$ are vector bundles on $X$ and $Y$ of rank $n_X$ and $n_Y$ as in theorem \ref{ThmVoeThom} respectively with the refinement in the Proposition \ref{refinement1}. Moreover, the isomorphism $tt_E^Y$ is independent from the choice of the projective embeddings up to a unique canonical isomorphism.  
\end{prop}
\begin{proof}
Let us denote by $s_Y$ the $0$-section of the vector bundle $p_{f^*V_X}: f^*V_X \rightarrow Y$ and by $s'_Y : Y \rightarrow T_f$ the $0$-section of the relative tangent bundle. Let $\sV(X)$ and $\sV(Y)$ be the categories of virtual bundles on $X$ and $Y$ respectively. We have 
\begin{multline*}
E^{p,q}(Y,f^*V_X-T_f) \stackrel{defn}{=} SH(Y)[\mathbf{1}_Y,Th_Y(T_f) \wedge_Y s_Y^!p_{f^*V_X}^*E_Y^{2(n_X+d_X-d_Y),n_X+d_X-d_Y} \wedge S^{p,q}] \cong \\ SH(Y)[p_{f^*V_X\#}s_{Y!}Th_Y(-T_f),E_Y^{2(n_X+d_X-d_Y),n_X+d_X-d_Y} \wedge S^{p,q}] \cong \\ SH(Y)[Th_Y(f^*V_X)\wedge_YTh_Y(-T_f),E_Y^{2(n_X+d_X-d_Y),n_X+d_X-d_Y} \wedge S^{p,q}] \cong \\ SH(Y)[Th_Y(f^*V_X-T_f),E_Y^{2(n_X+d_X-d_Y),n_X+d_X-d_Y} \wedge S^{p,q}],
\end{multline*}
where we write $Th_Y(f^*V_X - T_f)$ for the Thom transformation 
\[Th_Y(f^*V_X - T_f) = Th(s_Y,p_{f^*V_X}) \circ Th^{-1}(s'_Y,p_{T_f})(\mathbf{1}_Y).\]
Let $\pi_Y: Y \rightarrow \Spec k$ be the structure morphism of $Y$. By the adjunction $(\pi_{Y\#},\pi_Y^*)$ and since $E_Y = \pi_Y^* E$, we have then 
\begin{multline*}
SH(Y)[Th_Y(f^*V_X-T_f),E_Y^{2(n_X+d_X-d_Y),n_X+d_X-d_Y} \wedge S^{p,q}] \cong \\ SH(k)[Th(f^*V_X - T_f),E^{2(n_X+d_X-d_Y),n_X+d_X-d_Y} \wedge S^{p,q}],
\end{multline*}
which comes from the fact that (cf. \cite[Thm. 1.5.9]{Ay08} and \cite[Rem. 1.5.10]{Ay08}):
\[\pi_{Y\#}(Th(s_Y,p_{f^*V_X}) \circ Th^{-1}(s'_Y,T_f))(\mathbf{1}_Y) \cong Th(f^*V_X-T_f).\]
Now we apply the Voevodsky's theorem \ref{ThmVoeThom} with a refinement as in Proposition \ref{refinement1}. After fixing an embedding $X \inj \P^N$ we have in $\sV(X)$: 
\[V_X \oplus T_X \cong \sO_X^{N^2+2N}.\]
Since $f$ is projective, there is a factorization
\[\xymatrix{Y \ar[rr]^f \ar[rd] && X \\ & \P^M \times_k X \ar[ru]_{pr}  }\]
where $Y \inj \P^M \times_k X$ is a closed immersion. We take then the Segre embedding 
\[Y \inj \P^M \times \P^N \inj \P^{(N+1)(M+1)-1}.\]
and apply the Proposition \ref{refinement1}, so we have in $\sV(Y)$ a specific isomorphism  
\[V_Y \oplus T_Y \cong \sO_Y^{((N+1)(M+1)-1)^2 + 2((N+1)(M+1)-1)}.\]
One has a functor \cite[\S 4]{Del87}
\[f^*: \sV(X) \rightarrow \sV(Y).\]
Since $f$ is smooth we have an exact sequence
\[0 \rightarrow f^*T_X \rightarrow T_Y \rightarrow T_f \rightarrow 0,\]
which gives rise to an isomorphism in $\sV(Y)$ 
\[f^*T_X \oplus T_f \cong T_Y.\]
So we have then in $\sV(Y)$ an isomorphsim
\[f^*V_X - T_f \stackrel{\cong}{\longrightarrow} V_Y + \sO_Y^{(n_X+d_X)-(n_Y+d_Y)},\]
where $-$ means $+$ the opposite object as explained in \cite{Del87} and 
\[n_X+d_X = N^2 + 2N \]
\[n_Y + d_Y = ((N+1)(M+1)-1)^2 + 2((N+1)(M+1)-1).\]
Now since $Th$ defines a functor (cf. \cite[Def. 4.1.2]{Rio10})
\[Th: \sV(Y) \rightarrow SH(k),\]
where $\sV(Y)$ is the category of virtual bundles on $Y$, we can conclude that there is canonical isomorphism
\[Th(f^*V_X - T_f) \cong Th(V_Y) \wedge S^{2(n_X+d_X)-2(n_Y+d_Y),(n_X+d_X)-(n_Y+d_Y)},\]
where the right hand side is by \ref{refinement2} canonical isomorphic to $D(Y_+) \wedge S^{2(n_X+d_X),(n_X+d_X)}$. So we can conclude that there is an isomorphism 
\[tt_E^Y: E^{p,q}(Y,f^*V_X-T_f) \cong SH(k)[Th(V_Y),E^{2n_Y,n_Y}\wedge S^{p,q}] = E^{p+2n_Y,q+n_Y}(Th(V_Y)).\]
Now we have to show that this isomorphism is independent from the projective embeddings up to a unique canonical isomorphism. Let $Y \inj \P^{N'}$ be any closed embedding. Then we still have a canonical isomorphism
\[Th(f^*V'_{X}-T_f) \cong D(Y_+) \wedge S^{-2*,-*} \]
As $D(Y_+)$ is unique up to a canonical isomorphism we can conclude that $tt_E^Y$ is independent from the choice of the embeddings up to a unique canonical isomorphism.  
\end{proof}
\begin{rem}{\rm
As pointed out by M. Levine, one can simplify the arguments in the Proposition above by using the maps 
\[\mathbb{S}^0 \rightarrow S^{-d^2-2d} \wedge \Sigma^{\infty}_{T,+}Th(V_X) \wedge \Sigma^{\infty}_{T,+}X \]
and
\[S^{-d^2-2d} \wedge \Sigma^{\infty}_{T,+}Th(V_X) \wedge \Sigma^{\infty}_{T,+}X \rightarrow \mathbb{S}^0,\]
which rigidify the situation considerably. 
} 
\end{rem}
\begin{rem}{\rm
We will show later that with the refinement of $Th(V_X)$ as in \ref{refinement2} the isomorphism $tt_E^Y$ is natural in sense that it is compatible with duality. 
}
\end{rem}
\begin{rem}{\rm
The isomorphism $tt_E^Y(V_Y)$ in \ref{isothomtiwst} is a natural candidate for a replacement of the twisted Thom isomorphism $th_E$ in \ref{isothomE} in case of $E$-cohomology twisted by formal difference of vector bundles. But we should remind the reader that we can only compute for a very particular case, namely $\xi = V_X$, where $V_X$ is the vector bundle as in Voevodsky's theorem \ref{ThmVoeThom}.  
}
\end{rem}
Now we can compare: 
\begin{cor}\label{corSurprising}
Let $f: Y\rightarrow X$ be a smooth projective morphism of projective smooth $k$-schemes. There is an isomorphism up to a natural isomorphisms induced by the natural canonical isomorphism between duals 
\[E^{p,q}(Y,f^*V_X-T_f) \stackrel{tt_E^Y}{\rightarrow} E^{p+2n_Y,q+n_Y}(Th(V_Y)) \stackrel{th_E^Y(V_Y)^{-1}}{\longrightarrow} E^{p,q}(Y,V_Y).\]
\end{cor}
\begin{proof}
This is a consequence of \ref{isothomE} and \ref{isothomtiwst}. 
\end{proof}
\begin{rem}{\rm
The Corollary \ref{corSurprising} is a surprising fact to us. At a first glance we have the impression that $E^{p,q}(Y,f^*V_X - T_f)$ should depend relatively wrt. $f$ and $X$. At the end it turns out that $E^{p,q}(Y,f^*V_X-T_f)$ is isomorphic to $E^{p,q}(Y,V_Y)$, which depends absolutely only on $Y$. But it is clear that this is not the case for a general vector bundle $\xi$.  
}
\end{rem}
Let $f:Y\rightarrow X$ be a smooth projective morphism of smooth projective $k$-schemes of dimension $d_Y$ and $d_X$ respectively. By Atiyah-Spanier-Whitehead duality and by the Proposition \ref{refinement2} we obtain its dual morphism in $SH(k):$
\[f^{\vee}: X^{\vee} = \Sigma^{\infty}_{T,+}Th(V_X) \wedge S^{-2(n_X+d_X),-(n_X+d_X)} \rightarrow Y^{\vee} = \Sigma^{\infty}_{T,+}Th(V_Y)\wedge S^{-2(n_Y+d_Y),-(n_Y+d_Y)},\]
where $V_X$ and $V_Y$ are vector bundles on $X$ and $Y$ of rank $n_X$ and $n_Y$ as in theorem \ref{ThmVoeThom} with a refinement in \ref{refinement1} respectively. By taking pullback of this map on $E$-cohomology and appyling Thom isomorphism we obtain a pushforward 
\begin{multline*}
E^{p+2d_Y,q+d_Y}(Y,V_Y) \stackrel{th_E}{\cong} E^{p+2(n_Y+d_Y),q+(n_Y+d_Y)}(Th(V_Y)) \stackrel{(f^{\vee})^*}{\longrightarrow} \\ E^{p+2(n_X+d_X),q+(n_X+d_X)}(Th(V_X)) \stackrel{th_E}{\cong} E^{p+2d_X,q+d_X}(X,V_X).
\end{multline*}
We show that the two pushforwards are the same and the isomorphism $tt^Y_E$ is natural in the sense that it is compatible with the duality in the following:  
\begin{prop}\label{propthEpushforward}
Let $f: Y \rightarrow X$ be a smooth projective $k$-morphism between smooth projective $k$-schemes. One has a commutative diagram up to natural isomorphisms induced by natural $2$-isomorphisms and the natural canonical isomorphism between duals 
\[\xymatrix{E^{p+2d_Y,q+d_Y}(Y,f^*V_X - T_f) \ar[r]^{f_*} \ar[d]^{\cong}_{tt_E^Y} & E^{p+2d_X,q+d_X}(X,V_X) \ar[d]^{\cong}_{th^X_E(V_X)} \\ E^{p+2(n_Y+d_Y),q+(n_Y+d_Y)}(Th(V_Y)) \ar[r]_{(f^{\vee})^*} & E^{p+2(n_X+d_X),q+(n_X+d_X)}(Th(V_X)) }\]
\end{prop}
\begin{proof}
Let us denote by $p_{V_X}: V_X \rightarrow X$ the duality vector bundle on $X$ (cf. \ref{ThmVoeThom}) with the zero-section $s_X: X \rightarrow V_X$ and $s_Y: Y \rightarrow f^*V_X$ the $0$-section of the pullback bundle $p_{f^*V_X}: f^*V_X \rightarrow Y$.  By construction, the first pushforward map is the following composition:
\begin{multline*}
f_*: E^{p+2d_Y,q+d_Y}(Y,f^*V_X-T_f) \stackrel{defn}{=} SH(Y)[\mathbf{1}_Y,Th_Y(T_f) \wedge_Y s_Y^!p_{f^*V_X}^*E_Y^{2(n_X+d_X),(n_X+d_X)} \wedge S^{p,q}] \\ \stackrel{(1)}{\cong} SH(Y)[p_{f^*V_X\#}s_{Y!}Th_Y(-T_f),E_Y^{2(n_X+d_X),(n_X+d_X)} \wedge S^{p,q}] \stackrel{(2)}{\cong} \\ SH(Y)[Th_Y(f^*V_Y)\wedge_Y Th_Y(-T_f),E_Y^{2(n_X+d_X),(n_X+d_X)} \wedge S^{p,q}] \stackrel{(3)}{\cong} \\ SH(X)[f_{\#}(Th_Y(f^*V_Y)\wedge_Y Th_Y(-T_f)),E_X^{2(n_X+d_X),(n_X+d_X)} \wedge S^{p,q}] \stackrel{(4)}{\cong} \\ SH(X)[Th_X(V_X) \wedge_X f_{\#}Th_Y(-T_f),E_X^{2(n_X+d_X)} \wedge S^{p,q}] \stackrel{(5)}{\cong} \\ SH(X)[Th_X(V_X)\wedge_X D_X(f_{\#}\mathbf{1}_Y),E_X^{2(n_X+d_X),(n_X+d_X)}\wedge S^{p,q}] \stackrel{(6)}{\cong} \\ SH(X)[Th_X(V_X),f_{\#}(\mathbf{1}_Y) \wedge_X E_X^{2(n_X+d_X),(n_X+d_X)} \wedge S^{p,q}] \stackrel{(7)}{\rightarrow} \\ SH(X)[Th_X(V_X),E_X^{2(n_X+d_X),(n_X+d_X)} \wedge S^{p,q}] \stackrel{th_E^X(V_X)^{-1}}{\cong} E^{p+2d_X,q+d_X}(X,V_X),
\end{multline*}    
where $(1)$ is the natural isomorphism induced by adjunction of Thom transformations, $(2)$ is the natural isomorphism given by composing Thom transformations, $(3)$ is the natural isomorphism given by the adjunction $(f_{\#},f^*)$
\[(3)(-) = \varepsilon_{(f_{\#},f^*)}(-)\circ f_{\#}(-),\]
$(4)$ is the natural isomorphism given by projection formula $Pr_{\#}^*(f)$
\[(4)(-) = \varepsilon_{(f_{\#},f^*)}(-) \circ f_{\#}(\eta_{(f_{\#},f^*)}(-) \wedge_Y \id)(-) ,\]
$(5)$ is the natural isomorphism induced by $f_{\#}Th_Y(-T_f) \cong D_X(f_{\#}\mathbf{1}_Y)$, $(6)$ is the natural isomorphism induced by adjunction of duality in $SH(X)$:
\[(6)(-) = (\id_{D_X(f_{\#}\mathbf{1}_Y)} \wedge_X \ - ) \circ (coev_{f_{\#}\mathbf{1}_Y} \wedge \id_{-}),\]
and finally $(7)$ is the pushforward induced by the counit $\eta_{(f_{\#},f^*)}: f_{\#}(\mathbf{1}_Y) \cong f_{\#}f^*(\mathbf{1}_X) \rightarrow \mathbf{1}_X$:
\[(7)(-) = \varepsilon_{(f_{\#},f^*)}(-) \wedge_X - .\]
The last isomorphism is the inverse of the twisted Thom isomorphism. So we have
\[th_E^X(V_X) \circ f_* = (7)\circ \cdots \circ (1).\]
The map $tt_E^Y$ is the following composition 
\begin{multline*}
tt_E^Y: \\ E^{p+2d_Y,q+d_Y}(Y,f^*T_X-T_f) \stackrel{defn}{=} SH(Y)[\mathbf{1}_Y,Th_Y(T_f) \wedge_Y s_Y^!p_{f^*V_X}^*E_Y^{2(n_X+d_X),(n_X+d_X)} \wedge S^{p,q}] \\ \stackrel{(1')}{\cong} SH(Y)[p_{f^*V_X\#}s_{Y!}Th_Y(-T_f),E_Y^{2(n_X+d_X),(n_X+d_X)} \wedge S^{p,q}] \stackrel{(2')}{\cong} \\ SH(Y)[Th_Y(f^*V_X)\wedge_Y Th_Y(-T_f),E_Y^{2(n_X+d_X),(n_X+d_X)} \wedge S^{p,q}] \stackrel{(3')}{\cong} \\ SH(Y)[Th_Y(f^*V_X - T_f),E_Y^{2(n_X+d_X),(n_X+d_X)} \wedge S^{p,q}] \stackrel{(4')}{\cong} \\ SH(Y)[Th_Y(V_Y) \wedge S^{2(n_X+d_X)-2(n_Y+d_Y),(n_X+d_X)-(n_Y+d_Y)},E_Y^{2(n_X+d_X),(n_X+d_X)} \wedge S^{p,q}] \stackrel{(5')}{\cong} \\ SH(Y)[Th_Y(V_Y),E_Y^{2(n_Y+d_Y),(n_Y+d_Y)} \wedge S^{p,q}] \stackrel{(6')}{\cong} E^{p+2(n_Y+d_Y),q+(n_Y+d_Y)}(Th(V_Y)),
\end{multline*}
where $(1') = (1)$ , $(2') = (2)$, $(3')$ is the natural isomorphism induced by composing Thom transformations, $(4')$ is induced by the isomorphism in $\sV(Y)$: 
\[f^*V_Y - T_f \cong V_Y + \sO_Y^{(n_X+d_X)-(n_Y+d_Y)},\]
$(5')$ is the cancellation in $SH(Y)$ and finally $(6')$ is the natural isomorphism induced by the adjunction $(\pi_{Y\#},\pi_Y^*)$ with $\pi_Y: Y \rightarrow \Spec k$ is the structure morphism of $Y$: 
\[(6')(-) = \varepsilon_{(\pi_{Y\#},\pi_Y^*)}(-)\circ \pi_{Y\#}(-). \]  
Let us consider the diagram 
\[\xymatrix{Y \ar[rr]^f \ar[rd]_{\pi_Y} && X \ar[ld]^{\pi_X} \\ & \Spec k    }\]
We have 
\begin{multline*}
[SH(Y)[Th_Y(f^*V_X - T_f), E_Y^{2(n_X+d_X),(n_X+d_X)} \wedge S^{p,q}] \stackrel{\cong}{\longrightarrow} \\ SH(X)[f_{\#}Th_Y(f^*V_X-T_f),E_X^{2(n_X+d_X),(n_X+d_X)} \wedge S^{p,q}],
\end{multline*}
where the natural isomorphism is induced by the adjunction $(f_{\#},f^*)$ as $f$ is smooth: 
\[\varepsilon_{(f_{\#},f^*)} \circ f_{\#}(-).\]
By the projection formula $Pr_{\#}^*(f)$ we have a natural isomorphism
\begin{multline*}
SH(X)[f_{\#}Th_Y(f^*V_X-T_f),E_X^{2(n_X+d_X),(n_X+d_X)} \wedge S^{p,q}] \stackrel{\cong}{\longrightarrow} \\ 
SH(X)[Th_X(V_X) \wedge_X f_{\#}Th_Y(-T_f),E_X^{2(n_X+d_X)} \wedge S^{p,q}],
\end{multline*}
which is explicitly written as 
\[\varepsilon_{(f_{\#},f^*)}(-)\circ f_{\#}(\eta_{(f_\#,f^*)}(-) \wedge_Y \id)(-).\]
But since $f$ is smooth and projective 
\[\mathfrak{p}_f: f_\#Th_Y(-T_f) \stackrel{\cong}{\longrightarrow} f_*(\mathbf{1}_Y) \cong f_*f^*(\mathbf{1}_X). \]
Thanks to Proposition \ref{refinement2} the composition $(f^{\vee})^* \circ tt^Y_E$ on $E$-cohomology is nothing but just the composition of natural isomorphisms above with the map induced by the unit $\mathbf{1}_X \rightarrow f_*f^*\mathbf{1}_X$ of the adjunction $(f^*,f_*)$: 
\begin{multline*}
 \circ \eta_{(f^*,f_*)}: SH(X)[Th_X(V_X) \wedge_X f_*f^*(\mathbf{1}_X),E_X^{2(n_X+d_X),(n_X+d_X)} \wedge S^{p,q}] \longrightarrow \\ SH(X)[Th_X(V_X),E_X^{2(n_X+d_X),(n_X+d_X)} \wedge S^{p,q}],
\end{multline*}
then followed by the natural isomorphism induced from the adjunction $(\pi_{X\#},\pi_X^*)$:
\begin{multline*}
\varepsilon_{(\pi_{X\#},\pi_X^*)}\circ \pi_{X\#} : SH(X)[Th_X(V_X),E_X^{2(n_X+d_X),(n_X+d_X)} \wedge S^{p,q}] \stackrel{\cong}{\longrightarrow} \\ SH(k)[\pi_{X\#}Th_X(V_X),E^{2(n_X+d_X),(n_X+d_X)} \wedge S^{p,q}] = E^{p+2(n_X+d_X),q+(n_X+d_X)}(Th(V_X)). 
\end{multline*}
Indeed, by the very construction of the $6$ operations formalism \cite[Thm. 4.5.23]{Ay08}, the stabilization functor $\Sigma^{\infty}_{T,+}: Sm/k \rightarrow SH(k)$ induces a morphism in $SH(k)$:
\[\Sigma^{\infty}_{T,+}(f): \Sigma^{\infty}_{T,+}Y \rightarrow \Sigma^{\infty}_{T,+}X,\]
which can be understood as a morphism $\pi_{Y\#}(\mathbf{1}_Y) \rightarrow \pi_{X\#}(\mathbf{1}_X)$, which in turn is the composition $\pi_{X\#}\circ \varepsilon_{(f_{\#},f^*)}(\mathbf{1}_X)$. In terms of six operations and by the Proposition \ref{refinement2}, the dual objects in $SH(k)$ are: 
\[Th(V_X) \wedge S^{-2(n_X+d_X),-(n_X+d_X)} \stackrel{\cong}{\longrightarrow} X^{\vee} = D_k(\pi_{X\#}\mathbf{1}_X),\]
and
\[Th(V_Y)\wedge S^{-2(n_Y+d_Y),-(n_Y+d_Y)} \stackrel{\cong}{\longrightarrow} Y^{\vee} = D_k(\pi_{Y\#}\mathbf{1}_Y).\]
Hence the pullback map $(f^{\vee})^*$ is just the pullback $SH(k)[-,E^{*,*}\wedge S^{*,*}]$ of the map 
\[D_k(\pi_{X\#} \circ \varepsilon_{(f_{\#},f^*)}).\]
We have to check that 
\[(7)\circ \cdots \circ (3) = D_k(\pi_{X\#}\circ \varepsilon_{(f_{\#},f^*)})\circ (6') \circ \cdots \circ (3'),\]
which means that we have to check 
\begin{multline*}
D_k(\pi_{X\#} \circ \varepsilon_{(f_{\#},f^*)}(-)) \circ \varepsilon_{(\pi_{Y\#},\pi_Y^*)}\circ \pi_{Y\#}(-) = (\varepsilon_{(f_{\#},f^*)}(-) \wedge_X - ) \circ (\id_{D_X(f_{\#}\mathbf{1}_Y)} \wedge_X - ) \circ \\ (coev_{f_{\#}\mathbf{1}_Y} \wedge_X \id_{-}) \circ \varepsilon_{(f_{\#},f^*)}(-) \circ f_{\#}(\eta_{(f_{\#},f^*)}(-) \wedge_Y \id)(-) \circ \varepsilon_{(f_{\#},f^*)} \circ f_{\#}(-).
\end{multline*}
But this is clear, since for a smooth projective morphism $\pi : T \rightarrow S$ one has a natural $2$-isomorphism 
\[D_S(\pi_*(-)) \stackrel{\cong}{\longrightarrow} \pi_*D_T(- \wedge_T Th_T(T_f)),\]
which is the composition 
\begin{multline*}
D_S(f^*(-)) = \underline{\Hom}_S(\pi_*(-),\mathbf{1}_S) \stackrel{\cong}{\longrightarrow} \underline{\Hom}_S(\pi_!(-),\mathbf{1}_S) \stackrel{\cong}{\longrightarrow} \pi_*\underline{\Hom}_T((-),\pi^!\mathbf{1}_S) \stackrel{\cong}{\longrightarrow} \\ \pi_*\underline{\Hom}_T((-),\pi^*\mathbf{1}_S \wedge_T Th_T(-T_f)) \stackrel{\cong}{\longrightarrow} \pi_*\underline{\Hom}_T((-)\wedge_T Th_T(T_f),\mathbf{1}_T) = \pi_*D_T(-\wedge_T Th_T(T_f)).
\end{multline*}
The equality, which we need to check above, follows simply from this fact and from the fact that, we have a natural $2$-isomorphism 
\[f_{\#} Th_Y(s_Y,p_{f^*V_X}) \stackrel{\cong}{\longrightarrow} Th_X(s_X,p_{V_X})f_{\#}\] 
as $f$ is assumed to be smooth (cf. \cite[Thm. 1.5.9]{Ay08}). 
\end{proof}
Let us construct the pullback for twisted $E$-cohomology of formal difference of vector bundles along a cartesian square. Let
\[\xymatrix{Y' \ar[d]_g \ar[r]^v & Y \ar[d]^f \\ X' \ar[r]_u & X}\]
be a cartesian square of projective smooth $k$-schemes, where $f$ is smooth projective of relative dimension $d = \dim(Y)-\dim(X)$ and $u$ is any morphism. Let $pi_{\xi}: \xi \rightarrow X$ be a vector bundle of rank $r$ and denote by $s_Y : Y \rightarrow f^*\xi$ the $0$-section of the pullback bundle $p_{f^*\xi}: f^*\xi \rightarrow Y$. Consider 
\[E^{p,q}(Y,f^*\xi-T_f) \stackrel{defn}{=} SH(Y)[\mathbf{1}_Y,Th_Y(T_f) \wedge_Y s_Y^!p_{f^*\xi}^*E_Y^{2(r-d),(r-d)} \wedge S^{p,q}].\]
By adjunction of Thom transformation we have 
\begin{multline*}
SH(Y)[\mathbf{1}_Y,Th_Y(T_f) \wedge_Y s_Y^!p_{f^*\xi}^*E_Y^{2(r-d),(r-d)}] \stackrel{\cong}{\longrightarrow} \\ SH(Y)[Th_Y(f^*\xi),Th_Y(T_f) \wedge_Y E_Y^{2(r-d),(r-d)}],
\end{multline*}
where the isomorphism is 
\[ev_{Th_Y(f^*\xi)} \circ (\id_{Th_Y(f^*\xi)} \wedge -) \]
By applying the functor $v^*: SH(Y) \rightarrow SH(Y')$ we have an induced map 
\[E^{p,q}(Y,f^*\xi-T_f) \rightarrow SH(Y')[v^*Th_Y(f^*\xi),v^*(Th_Y(T_f) \wedge_Y E_Y^{2(r-d),(r-d)}\wedge S^{p,q})].\]
We have $v^*Th_Y(f^*\xi) \cong Th_{Y'}(v^*f^*\xi) = Th_{Y'}(g^*u^*\xi)$ as $\xi$ is an actual bundle. Since $v^*$ is a  monoidal functor, so we have 
\[v^*(Th_Y(T_f)\wedge_Y E_Y^{2(r-d),(r-d)}) \cong v^*Th_Y(T_f) \wedge_{Y'} E_{Y'}^{2(r-d),(r-d)}.\]
But we have $v^*Th_Y(T_f) \cong Th_{Y'}(v^*T_f)$, since $T_f$ is an actual bundle. By \cite[16.5.12.2]{EGA4} one has $v^*T_f \cong T_g$, so $v^*Th_Y(T_f) \cong Th_{Y'}(T_g)$. So we obtain the pullback map for twisted $\mathbf{E}$-cohomology of formal difference of vector bundle
\[E^{p,q}(Y,f^*\xi - T_f) \rightarrow E^{p,q}(Y',v^*f^*\xi-T_g) = E^{p,q}(Y',g^*u^*\xi-T_g).\]
\begin{prop}
Consider the composition of cartesian squares of smooth projective $k$-schemes
\[\xymatrix{Y'' \ar[r]^{v'} \ar[d]_h & Y' \ar[r]^{v} \ar[d]_g & Y \ar[d]^f \\ X'' \ar[r]_{u'} & X' \ar[r]_{u} & X}\]
where $f$ is a smooth projective morphism, $u$ and $u'$ are  morphisms. Let $\xi$ be a vector bundle on $X$. Then up to natural isomorphisms induced by the natural $2$-isomorphism induced by $(-\circ -)^* \stackrel{\cong}{\longrightarrow} (-)^* \circ (-)^*$ one has 
\[(v \circ v')^* = v'^* \circ v^*: E^{p,q}(Y,f^*\xi-T_f) \rightarrow E^{p,q}(Y'',v'^*v^*f^*\xi - T_f) = E^{p,q}(Y'',h^*u'^*u^*\xi - T_g). \]
\end{prop}
\begin{proof}
Obvious. 
\end{proof}
\begin{prop}(projective smooth base change)
\label{projsmbasechange}
Consider a cartesian square of projective smooth $k$-schemes 
\[\xymatrix{Y' \ar[d]_g \ar[r]^v & Y \ar[d]^f \\ X' \ar[r]_u & X}\]
where $f$ is smooth projective of relative dimension $d = \dim(Y) - \dim(X)$ and $u$ is a morphism. Let $\xi/X$ be a vector bundle of rank $r$. One has a commutative diagram up to natural isomorphisms induced by natural $2$-isomorphisms
\[\xymatrix{E^{p,q}(Y,f^*\xi-T_f) \ar[r]^{f_*} \ar[d]_{v^*} & E^{p-2d,q-d}(X,\xi) \ar[d]^{u^*} \\ E^{p,q}(Y',g^*u^*\xi - T_g) \ar[r]_{g_*}  & E^{p-2d,q-d}(X',u^*\xi) }\] 
\end{prop}
\begin{proof}
It is quite straightforward. We write $s_Y: Y \rightarrow f^*\xi$ and $s_{Y'} : Y \rightarrow v^*f^*\xi$ for the $0$-sections of the vector bundles $p_{f^*\xi}: f^*\xi \rightarrow Y$ and $p_{v^*f^*\xi}: v^*f^*\xi \rightarrow Y'$ respectively. $u^*f_*$ is the following composition 
\begin{multline*}
E^{p,q}(Y,f^*\xi-T_f) \stackrel{defn}{=} SH(Y)[\mathbf{1}_Y,Th_Y(T_f) \wedge_Y s_Y^!p_{f^*\xi}^*E_Y^{2(r-d),(r-d)} \wedge S^{p,q}] \stackrel{\alpha_1}{\cong} \\ SH(Y)[Th_Y(f^*\xi) \wedge_Y Th_Y(-T_f),E_Y^{2(r-d),(r-d)} \wedge S^{p,q}] \stackrel{\varepsilon_{(f_{\#},f^*)}\circ f_{\#}}{\cong} \\ SH(X)[f_{\#}(Th_Y(f^*\xi)\wedge_Y Th_Y(-T_f)),E_X^{2(r-d),(r-d)} \wedge S^{p,q}] \stackrel{Pr_{\#}^*(f)}{\cong} \\ SH(X)[Th_X(\xi)\wedge_X f_{\#}Th_Y(-T_f),E_X^{2(r-d),(r-d)}\wedge S^{p,q}] \stackrel{\alpha_2}{\cong} \\ SH(X)[Th_X(\xi),f_{\#}(\mathbf{1}_Y)\wedge_X E_X^{2(r-d),(r-d)} \wedge S^{p,q}] \stackrel{-\circ\varepsilon_{(f_{\#},f^*)}}{\longrightarrow} SH(X)[Th_X(\xi),E_X^{2(r-d),(r-d)} \wedge S^{p,q}] \\ \stackrel{u^*(-)}{\longrightarrow} SH(X')[u^*Th_X(\xi),E_{X'}^{2(r-d),(r-d)} \wedge S^{p,q}] \stackrel{\alpha_3}{\longrightarrow} SH(X')[Th_{X'}(u^*\xi),E_{X'}^{2(r-d),(r-d)}\wedge S^{p,q}] \\ \stackrel{th_E^{X'}(u^*\xi)^{-1}}{\cong} E^{p-2d,q-d}(X',u^*\xi),
\end{multline*}
where 
\[\alpha_1(-) = ev_{Th_Y(f^*\xi)} \circ (\id_{Th_Y(f^*\xi)} \wedge -) \circ ev_{Th_Y(T_f)} \circ (\id_{Th_Y(T_f)} \wedge -),\]
and 
\[\alpha_2(-) = (\id_{f_{\#}\mathbf{1}_Y} \wedge -) \circ (coev_{f_{\#}Th_Y(-T_f)} \wedge \id).\]
$\alpha_3$ is the natural isomorphism
\[\alpha_3 = Ex_*^*(\Delta_b) \circ Ex_{\#}^*(\Delta_a)^{-1} \]
$\Delta_a$ is the Cartesian square
\[\xymatrix{u^*\xi \ar[r]^{u_{\xi}} \ar[d]_{p_{u^*\xi}} & \xi \ar[d]^{p_{\xi}} \\ X' \ar[r]_{u} & X}\]
\[Ex_{\#}^*(\Delta_a)^{-1} : u^*p_{\xi \#} \stackrel{\cong}{\longrightarrow} p_{u^*\xi \#}u_{\xi}^*.\]
$\Delta_b$ is the Cartesian square
\[\xymatrix{X' \ar[d]_u \ar[r]^{s_{X'}} & u^*\xi \ar[d]^{u_{\xi}} \\ X \ar[r]_{s_X} & \xi}\]
\[Ex_*^*(\Delta_b) : u^*_{\xi}s_{X!} \cong u^*_{\xi}s_{X*} \stackrel{\cong}{\longrightarrow} s_{X'*}u^* \cong s_{X'*}u^*.\]
$g_*v^*$ is the following composition
\begin{multline*}
E^{p,q}(Y,f^*\xi-T_f) \stackrel{defn}{=} SH(Y)[\mathbf{1}_Y,Th_Y(T_f) \wedge_Y s_Y^!p_{f^*\xi}^*E_Y^{2(r-d),(r-d)} \wedge S^{p,q}] \stackrel{\beta_1}{\longrightarrow} \\ SH(Y)[Th_Y(f^*\xi), Th_Y(T_f) \wedge_{Y} E_{Y}^{2(r-d),(r-d)} \wedge S^{p,q}] \stackrel{v^*(-)}{\longrightarrow} \\ SH(Y')[v^*Th_Y(f^*\xi),v^*(Th_Y(T_f)\wedge_{Y}E_{Y}^{2(r-d),(r-d)} \wedge S^{p,q})] \stackrel{\beta_2}{\cong} \\ SH(Y')[Th_{Y'}(g^*u^*\xi) \wedge_{Y'}Th_{Y'}(-T_g),E_{Y'}^{2(r-d),(r-d)} \wedge S^{p,q}] \stackrel{\varepsilon_{(g_{\#},g^*)}\circ g_{\#}}{\cong} \\ SH(X')[g_{\#}(Th_{Y'}(g^*u^*\xi) \wedge_{Y'}Th_{Y'}(-T_g)),E_{X'}^{2(r-d),(r-d)} \wedge S^{p,q}] \stackrel{Pr_{\#}^*(g)}{\cong} \\ SH(X')[Th_{X'}(u^*\xi) \wedge_{X'} g_{\#}Th_{Y'}(-T_g)),E_{X'}^{2(r-d),(r-d)} \wedge S^{p,q}] \stackrel{\beta_3}{\cong} \\ SH(X')[Th_{X'}(u^*\xi),g_{\#}(\mathbf{1}_{Y'}) \wedge_{X'}E_{X'}^{2(r-d),(r-d)} \wedge S^{p,q}] \stackrel{\circ \varepsilon_{(g_{\#},g^*)}}{\longrightarrow} \\ SH(X')[Th_{X'}(u^*\xi),E_{X'}^{2(r-d),(r-d)} \wedge S^{p,q}] \stackrel{th_E^{X'}(u^*\xi)^{-1}}{\cong} E^{p-2d,q-2d}(X',u^*\xi),
\end{multline*}
where $\beta_1$ is the natural isomorphism 
\[\beta_1(-) = ev_{Th_Y(f^*\xi)} \circ (\id_{Th_Y(f^*\xi)} \wedge -),\]
and 
\[\beta_2(-) = \beta_2' \circ  Ex_{*}^*(\Delta_4) \circ Ex_{\#}^*(\Delta_3)^{-1} \circ  Ex_*^*(\Delta_2) \circ Ex_{\#}^*(\Delta_1)^{-1}.\]
$\Delta_1$ is the Cartesian square 
\[\xymatrix{g^*u^*\xi =  v^*f^* \xi \ar[r]^{\quad \quad v_{\xi}} \ar[d]_{p_{v^*f^*}\xi = p_{g^*u^*\xi}} & f^*\xi \ar[d]^{p_{f^*}\xi} \\ Y' \ar[r]_v & Y}\]
\[Ex_{\#}^*(\Delta_1)^{-1}: v^*p_{f^*\xi} \stackrel{\cong}{\longrightarrow} p_{v^*f^*\xi\#}v_{\xi}^*.\]
$\Delta_2$ is the Cartesian square
\[\xymatrix{Y' \ar[r]^{s_{Y'}} \ar[d]_v & v^*f^*\xi \ar[d]^{v_{\xi}} \\ Y \ar[r]_{s_Y} & f^*\xi}\]
\[Ex_*^*(\Delta_2): v_{\xi}^*s_{Y!} \cong v_{\xi}^*s_{Y*} \stackrel{\cong}{\longrightarrow} s_{Y'*}v^* \cong s_{Y'!}v^*. \]
$\Delta_3$ is the Cartesian square 
\[\xymatrix{T_g \cong v^* T_f \ar[r]^{\quad \quad v_{T_f}} \ar[d]_{p_{v^*T_f} = p_{T_g}} & T_f \ar[d]^{p_{T_f}} \\ Y' \ar[r]_v & Y}\]
\[Ex_{\#}^*(\Delta_3)^{-1}: v^*p_{T_f} \stackrel{\cong}{\longrightarrow} p_{v^*T_f\#}v_{T_f}^* \cong p_{T_g\#}v_{T_f}^*.\]
$\Delta_4$ is the Cartesian square
\[\xymatrix{Y' \ar[r]^{s_{Y'/T_g} \quad \quad} \ar[d]_v & v^*T_f \cong T_g \ar[d]^{v_{T_f}} \\ Y \ar[r]_{s_Y/T_f} & T_f}\]
\[Ex_*^*(\Delta_4): v_{T_f}^*s_{Y/T_f!} \cong v_{T_f}^*s_{T_f*} \stackrel{\cong}{\longrightarrow}s_{Y'/T_g*}v_{T_f}^* \cong s_{Y'/T_g!}v_{T_f}^*.\]
$\beta_2'$ is the natural isomorphism 
\[\beta_2'(-) = ev_{Th_{Y'}(-T_g)} \circ (\id_{Th_{Y'}(-T_g)} \wedge -).\] 
$\beta_3$ is the natural isomorphism 
\[\beta_3 = (\id_{g\#\mathbf{1}_{Y'}} \wedge - ) \circ (coev_{g_{\#}Th_{Y'}(-T_g)} \wedge \id_{-} ).\]
Gathering all together we have to check the following equality up to natural $2$-isomorphisms: 
\begin{multline*}
th_E^{X'}(u^*\xi)^{-1} \circ Ex^*_*(\Delta_b) \circ Ex^*_{\#}(\Delta_a)^{-1} \circ \epsilon_{(f_{\#},f^*)} \circ u^* \circ (\id_{f\#\mathbf{1}_Y} \wedge_X - ) \circ (coev_{f_{\#}Th_Y(-T_f)} \wedge_X \id_{-}) \\ \circ Pr_{\#}^*(f) \circ \varepsilon_{(f_{\#},f^*)} \circ f_{\#} \circ ev_{Th_Y(f^*\xi)} \circ (\id_{Th_Y(f^*\xi)} \wedge_Y - ) \circ ev_{Th_Y(T_f)} \circ (\id_{Th_Y(T_f} \wedge_Y - ) = \\ th_E^{X'}(u^*\xi)^{-1} \circ \varepsilon_{(g_{\#},g^*)} \circ (\id_{g_{\#}\mathbf{1}_{Y'}} \wedge_{X'} - ) \circ (coev_{g_{\#}Th_{Y'}(-T_g)} \wedge_{X'} \id_{-}) \circ Pr_{\#}^*(f) \circ \varepsilon_{(g_{\#},g^*)} \circ g_{\#} \circ \\ ev_{Th_{Y'}(-T_g)} \circ (\id_{Th_{Y'}(T_g)} \wedge_{Y'} - ) \circ Ex_{*}^*(\Delta_4) \circ Ex_{\#}^*(\Delta_3)^{-1}\circ Ex^*_*(\Delta_2) \circ Ex^*_{\#}(\Delta_1)^{-1} \\ \circ v^* \circ ev_{Th_Y(f^*\xi)} \circ (\id_{Th_Y(f^*\xi)} \wedge - ). 
\end{multline*}
This equality can be chased step by step by using the natural $2$-isomorphism 
\[f_{\#}u^* \stackrel{\cong}{\longleftarrow} g_{\#}v^*,\]
which is the following composition 
\[g_{\#}v^* \stackrel{\eta_{(f_{\#},f^*)}}{\longrightarrow} g_{\#}v^*f^*f_{\#} \cong g_{\#}(f\circ v)^*f_{\#} = g_{\#}(u \circ g)^*f_{\#} \cong g_{\#}g^*u^*f_{\#} \stackrel{\varepsilon_{(g_{\#},g^*)}}{\longrightarrow} u^*f_{\#}\]
and also the coherence of the exchange transformations.
\end{proof}
Now we construct the exceptional pullback for twisted $E$-cohomology. We keep the notation as above and let $i: T \inj Y$ be a regular embedding, where $T$ is a smooth $k$-scheme. Let $N_{T/Y}$ be the normal bundle of $T$ in $Y$. Let $Bl_T(Y)$ be the blow-up of $X$ with the center $Z$. Similarly, $Bl_{T \times \{0 \} }(Y \times \A^1)$ is the blow-up of $Y \times \A^1$ with the center $T \times \{0\}$. The deformation space is the $k$-scheme 
\[D_T(Y) \stackrel{defn}{=} Bl_{T \times \{0\}}(Y \times \A^1) - Bl_T(Y).\]
Note that $D_T(T) = T \times \A^1$ is a closed subscheme of $D_T(Y)$. The scheme $D_T(Y)$ is fibred over $\A^1$. The flat morphism
\[\pi: D_T(Y) \rightarrow \A^1 \]
has $\pi^{-1}(1) = Y$ and $\pi^{-1}(0) = N_{T/Y}$. One has a deformation diagram of closed pairs
\[(Y,T) \stackrel{\sigma_1}{\longrightarrow} (D_T(Y),T \times \A^1) \stackrel{\sigma_0}{\longleftarrow} (N_{T/Y},T). \]
The homotopy purity theorem of Morel-Voevodsky \cite[\S 3 Thm. 2.23]{MV01} states 
\[Y/Y-T \stackrel{\sigma_{1*}}{\longrightarrow} D_T(Y)/D_T(Y) - T \times \A^1 \stackrel{\sigma_{0*}}{\longleftarrow} Th(N_{T/Y})\]
are isomorphism in $\mathbf{Ho}_{\A^1,+}(k)$, which is generalized to motivic categories in \cite[Thm. 2.4.35]{CD10}. Consider now the adjunction 
\[i_! : SH(T) \rightleftarrows SH(Y) : i^!. \]
Let 
\[\xymatrix{T \ar[r]^{i} \ar[d]_g & Y \ar[d]^f \\ S \ar[r]_k & X }\]
be a cartesian square of smooth projective $k$-schemes, where $f$ is smooth projective of relative dimension $d = \dim(Y) - \dim(X)$, $k$ and $i$ are regular embeddings. Let $\xi$ be a vector bundle of rank $r$ on $X$. We define the exceptional pullback of twisted $E$-cohomology along a regular embedding $i: T \inj Y$ as the following composition:
\begin{multline*}
i^!: E^{p,q}(Y,f^*\xi-T_f) \stackrel{defn}{=} SH(Y)[\mathbf{1}_Y,Th_Y(T_f) \wedge_Y s_Y^!p_{f^*\xi}^* E_Y^{2(r-d),(r-d)} \wedge S^{p,q}] \stackrel{i^!(-)}{\longrightarrow} \\ SH(T)[i^!\mathbf{1}_Y,i^!Th_Y(T_f) \wedge_Y s_Y^!p_{f^*\xi}^*E_Y^{2(r-d),(r-d)} \wedge S^{p,q}] \stackrel{\varepsilon_{(i_!,i^!)} \circ i_!}{\cong} \\ SH(Y)[i_!i^!(\mathbf{1}_Y),Th_Y(T_f)\wedge_Y s_Y^!p_{f^*\xi}^*E_Y^{2(r-d),(r-d)} \wedge S^{p,q}] \stackrel{\circ Ex^{*!}}{\longrightarrow} \\ SH(Y)[i_!i^*(\mathbf{1}_Y),Th_Y(T_f)\wedge_Y s_Y^!p_{f^*\xi}^*E_Y^{2(r-d),(r-d)} \wedge S^{p,q}] \cong \\ SH(Y)[i_*\mathbf{1}_T,Th_Y(T_f)\wedge_Ys_Y^!p_{f^*\xi}^*E_Y^{2(r-d),(r-d)} \wedge S^{p,q}] \stackrel{i^*(-)}{\longrightarrow} \\ SH(T)[i^*i_*\mathbf{1}_T,i^*(Th_Y(T_f)\wedge_Y s_Y^!p^*_{f^*\xi}E_Y^{2(r-d),(r-d)} \wedge S^{p,q})] \cong \\ 
SH(T)[\mathbf{1}_T,i^*(Th_Y(T_f) \wedge_Y s^!_Yp^*_{f^*\xi}E_Y^{2(r-d),(r-d)}\wedge S^{p,q})] \stackrel{Ex^{!*}\circ Ex^{!*}}{\longrightarrow} \\
SH(T)[\mathbf{1}_T,Th_T(i^*T_f)\wedge_T s_T^!p_{i^*f^*\xi}^*E_T^{2(r-d),(r-d)} \wedge S^{p,q}] = E^{p,q}(T,i^*f^*\xi-T_g).
\end{multline*}
\begin{prop}
Let 
\[\xymatrix{T' \ar[d]_h \ar[r]^{i'} & T \ar[r]^i \ar[d]_g & Y \ar[d]^f \\ S' \ar[r]_{k'} & S \ar[r]_k & X  }\]
be a chain of cartesian squares of smooth projective $k$-schemes, where $f$ is smooth projective, $i,i',k,k'$ are regular embeddings. Then we have $(i\circ i')^! = i^! \circ i'^!$ up to natural isomorphisms induced by natural $2$-isomorphisms .
 \end{prop}
\begin{proof}
Obvious.
\end{proof}
By using deformation to the cone as discussed above, one can prove the following result. However, we will not need this result, so we just omit the proof.  
\begin{prop}
Consider a cartesian square of projective smooth $k$-schemes
\[\xymatrix{T \ar[r]^{i} \ar[d]_g & Y \ar[d]^f \\ S \ar[r]_k & X  }\]
where $f$ is smooth projective of relative dimension $d = \dim(Y)- \dim(X)$ and $k$ and $i$ are regular embeddings. Let $p_{\xi}: \xi \rightarrow X$ be a vector bundle of rank $r$. One has a commutative diagram up to a natural isomorphism
\[\xymatrix{ E^{p,q}(Y,f^*\xi-T_f) \ar[d]_{i^!} \ar[r]^{f_*} & E^{p-2d,q-d}(X,\xi) \ar[d]^{k^!} \\ E^{p,q}(T,g^*k^*\xi - T_g) \ar[r]_{g_*} & E^{p-2d,q-d}(S,k^*\xi) }\] 
\end{prop}
Let $p_{\xi}: \xi \rightarrow X$ and $p_{\xi'}: \xi' \rightarrow X$ be two vector bundles of rank $r$ and $r'$ resp. on $X$ with the zero sections $s: X \rightarrow \xi$ and $s':X \rightarrow \xi'$ respectively. Let $s'': X \rightarrow \xi \oplus  \xi'$ to be the zero section of the bundle $\xi \oplus \xi'$. We define the cup product 
\[\cup_E: E^{p,q}(X,\xi) \otimes E^{p',q'}(X,\xi') \rightarrow E^{p+p',q+q'}(X,\xi \oplus \xi')\]
as follow: Given morphisms in $SH(X)$ 
\[\alpha: \mathbf{1}_X \rightarrow s^!p_{\xi}^*E_X^{2r,r} \wedge S^{p,q}\]
and 
\[\beta: \mathbf{1}_X \rightarrow s'^!p_{\xi'}^*E_X^{2r',r'} \wedge S^{p',q'}.\]
Then 
\begin{multline*}
\alpha \cup_E \beta = \mu_E \circ (\alpha \wedge_X \beta):\mathbf{1}_X = \mathbf{1}_X \wedge_X \mathbf{1}_X \stackrel{\alpha \wedge_X \beta}{\longrightarrow} s^!p_{\xi}^*E_X^{2r,r} \wedge_X S^{p,q} \wedge s'^!p_{\xi'}^*E_X^{2r',r'} \wedge S^{p',q'}  \cong \\ \cong s''^!p_{\xi \oplus \xi'}^* E_X \wedge_X E_X \wedge S^{p+p'+2r+2r',q+q'+r+r'} \stackrel{\mu_E}{\rightarrow} s''^!p_{\xi \oplus \xi'}^* E_X^{2(r+r'),r+r'} \wedge S^{p+p',q+q'}.
\end{multline*}
\begin{rem}\label{remthE} {\rm
If $f: T \rightarrow S$ is a morphism of finite type between schemes, then we have 
\[f^*(E \wedge_S^{\mathbb{L}} F) = f^*E \wedge_T^{\mathbb{L}} f^*F.\]
}
\end{rem}
\begin{prop}(projection formula) \label{projectionformula}
Let $f : Y \rightarrow X$ be a smooth projective morphism of smooth projective $k$-schemes of relative dimension $d = \dim(Y)-\dim(X)$. Let $\xi$ and $\xi'$ be two vector bundles on $X$. Let $a \in E^{p,q}(X,\xi)$ and  $b \in E^{p',q'}(Y,f^*\xi'-T_f)$. Then one has up to natural isomorphisms
\[f_*(f^*a \cup_{E} b) = a \cup_{E} f_*b \] 
in $E^{p+p'-2d,q+q'-d}(X,\xi \oplus \xi')$. 
\end{prop}
\begin{proof}
This follows from the projective smooth base change \ref{projsmbasechange} by standard argument. Consider the commutative diagram 
\[\xymatrix{Y \ar[d]_f \ar[r]^{\Gamma_f \, \, \, \, \, \,} & Y \times_k X \ar[d]^{f \times \id} \\ X \ar[r]_{\Delta_X \, \, \, \,} & X\times_k X }\]
We have 
\[\Delta_X^*(f\times \id)_* = f_*\Gamma_f^* = f_*\Delta_X^*(f \times \id)^*,\]
where $\Gamma_f = (f \times \id)\Delta_X$. 
\end{proof}
\begin{prop}
Let $\mathbf{E} \in SH(k)$ be a motivic ring spectrum. Let $X,Y,Z \in SmProj(k)$. Let $\alpha \in \mathbf{E}^{2d_Y,d_Y}(X\times Y,pr_Y^{XY*}V_Y)$ and $\beta \in \mathbf{E}^{2d_Z,d_Z}(Y\times Z, pr_Z^{YZ*}V_Z)$, where $V_Y$ and $V_Z$ are the vector bundles given in theorem \ref{ThmVoeThom} with a refinement in \ref{refinement1}. Then we have up to natural isomorphisms 
\[pr_{XZ*}^{XYZ}(pr_{XY}^{XYZ*} \alpha \cup_{\mathbf{E}} pr_{YZ}^{XYZ*} \beta) \in \mathbf{E}^{2d_Z,d_Z}(X\times Z,pr_Z^{XZ*}V_Z).\] 
\end{prop}
\begin{proof}
This follows from our construction of pullback, pushforward and cup product and the projections fit to the following commutative diagram 
\[\xymatrix{X \times Z \ar@/^/[drrrrrr]^{pr^{XZ}_Z} \ar@/_/[rrddd]_{pr^{XZ}_X} \\ && X \times Y \times Z \ar[llu]^{pr^{XYZ}_{XZ}} \ar[d]_{pr^{XYZ}_{XY}} \ar[rr]^{pr^{XYZ}_{YZ}} && Y \times Z \ar[d]^{pr^{YZ}_Y} \ar[rr]_{pr^{YZ}_Z} && Z \\ && X \times Y \ar[d]^{pr^{XY}_X} \ar[rr]_{pr^{XY}_Y} && Y \\ && X     }\]
\end{proof}
\begin{prop}
Let $\mathbf{E} \in SH(k)$ be a motivic ring spectrum. Let $X,Y,Z,W \in SmProj(k)$. Let $\alpha \in \mathbf{E}^{2d_Y,d_Y}(X\times Y,pr_Y^{XY*}V_Y)$, $\beta \in \mathbf{E}^{2d_Z,d_Z}(Y\times Z, pr_Z^{YZ*}V_Z)$ and $\gamma \in \mathbf{E}^{2d_W,d_W}(Z \times W, pr_W^{ZW*}V_W)$. Let's denote 
\[\beta \circ \alpha = pr_{XZ*}^{XYZ}(pr_{XY}^{XYZ*} \alpha \cup_{\mathbf{E}} pr_{YZ}^{XYZ*} \beta),\]
and similarly for $\gamma \circ \beta$. Then $\circ$ is associative up to natural isomorphisms induced by $2$-isomorphisms.
\end{prop}
\begin{proof}
We have 
\begin{multline*}
\gamma \circ (\beta \circ \alpha) \stackrel{(1)}{=} pr_{XW*}^{XZW}(pr^{XZW*}_{XZ}(pr^{XYZ}_{XZ*}(pr^{XYZ*}_{XY}\alpha \cup_{\mathbf{E}} pr^{XYZ*}_{YZ}\beta))\cup_{\mathbf{E}} pr^{XZW*}_{ZW}\gamma) \stackrel{(2)}{=} \\ pr^{XZW}_{XW*}(pr^{XYZW}_{XZW*}(pr^{XYZW*}_{XYZ}(pr^{XYZ*}_{XY}\alpha \cup_{\mathbf{E}} pr^{XYZ*}_{YZ}\beta))\cup_{\mathbf{E}} pr^{XZW*}_{ZW} \gamma) \stackrel{(3)}{=} \\ pr^{XZW}_{XW*}(pr^{XYZW}_{XZW*}((pr^{XYZW*}_{XY}\alpha \cup_{\mathbf{E}} pr^{XYWZ*}_{YZ}\beta)\cup_{\mathbf{E}}pr^{XYZW*}_{XZW}pr^{XZW*}_{ZW}\gamma)) \stackrel{(4)}{=} \\ pr^{XYZW}_{XW*}(pr^{XYZW*}_{XY}\alpha \cup_{\mathbf{E}}(pr^{XYZW*}_{YZ}\beta \cup_{\mathbf{E}} pr^{XYZW*}_{ZW}\gamma)),
\end{multline*}
where $(1)$ is the definition, $(2)$ follows from smooth projective base change:
\[pr^{XZW}_{XW*}pr^{XYZ*}_{XZ} = pr^{XYZW}_{XZW*}pr^{XYZW*}_{XYZ},\]
$(3)$ follows from the compatibility of pullback and $\cup_{\mathbf{E}}$ (\ref{remthE}), functoriality of pullback (\ref{functorialityPullback}) and the projection formula (\ref{projectionformula}), $(4)$ follows from functoriality of pullback (\ref{functorialityPullback}) and pushforward (\ref{functorialitypushforward}), $(5)$ follows from the associativity of $\cup_{\mathbf{E}}$, which is a consequence of our requirement that $\mathbf{E}$ is a motivic ring spectrum (see the beginning of \S 4.1). Symmetrically, the last expression is exactly $(\gamma \circ \beta) \circ \alpha$.    
\end{proof}
\begin{defn}{\rm
Let $\mathbf{E} \in SH(k)$ be a motivic ring spectrum. We define the category of twisted $\mathbf{E}$-correspondences $\widetilde{Corr}_{\mathbf{E}}(k)$ to be the category, whose objects are 
\[Obj(\widetilde{Corr}_{\mathbf{E}}(k)) = Obj(SmProj(k))\] and morphisms are given by 
\[\widetilde{Corr}_{\mathbf{E}}(k)(X,Y) = \mathbf{E}^{2d_Y,d_Y}(X \times Y,pr_Y^{XY*}V_Y),\]
where $V_Y/Y$ is the vector bundle given in the theorem \ref{ThmVoeThom}. Given $\alpha \in \mathbf{E}^{2d_Y,d_Y}(X\times Y,pr_Y^{XY*}V_Y)$ and $\beta \in \mathbf{E}^{2d_Z,d_Z}(Y\times Z, pr_Z^{YZ*}V_Z)$ we define their composition to be 
\[\beta \circ \alpha = pr_{XZ*}^{XYZ}(pr_{XY}^{XYZ*} \alpha \cup_{\mathbf{E}} pr_{YZ}^{XYZ*} \beta),\]
which is associative up to  natural isomorphisms. 
}
\end{defn}
\begin{prop}
Let $\mathbf{E} \in SH(k)$ be a motivic ring spectrum. Let $X,Y,Z \in SmProj(k)$. Let $\alpha \in \mathbf{E}^{2(n_Y+d_Y),(n_Y+d_Y)}(X_+ \wedge Th(V_Y))$ and $\beta \in \mathbf{E}^{2(n_Z+d_Z),n_Z+d_Z}(Y_+ \wedge Th(V_Z))$. Then the following composition 
\begin{multline*}\beta \circ^{\dagger} \alpha: \Sigma^{\infty}_{T,+}X \wedge Th(V_Z) \stackrel{coev_Y}{\longrightarrow} \Sigma^{\infty}_{T,+}X \wedge Y \wedge Th(V_Y) \wedge S^{-2(n_Y+d_Y),-(n_Y+d_Y)} \wedge Th(V_Z) \\ \stackrel{\tau}{\rightarrow} \Sigma^{\infty}_{T,+} X \wedge Th(V_Y) \wedge Y \wedge Th(V_Z) \wedge S^{-2(n_Y+d_Y),-(n_Y+d_Y)} \stackrel{\alpha \wedge \beta}{\longrightarrow} \mathbf{E} \wedge^{\mathbb{L}}_S \mathbf{E} \wedge S^{2(n_Z+d_Z),(n_Z+d_Z)} \\ \stackrel{\mu_{\mathbf{E}}}{\longrightarrow} \mathbf{E} \wedge S^{2(n_Z+d_Z),(n_Z+d_Z)}   
\end{multline*}
lies in $\mathbf{E}^{2(n_Z+d_Z),(n_Z+d_Z)}(X_+ \wedge Th(V_Z))$, where 
\[coev_Y: \mathbb{S}^0 \rightarrow \Sigma^{\infty}_{T,+} Y \wedge Th(V_Y) \wedge S^{-2(n_Y+d_Y),-(n_Y+d_Y)}\]
is the coevaluation map of the Atiyah-Spanier-Whitehead duality on $Y$.
\end{prop}
\begin{proof}
Trivial. 
\end{proof}
\begin{prop}
Let $\mathbf{E} \in SH(k)$ be a motivic ring spectrum. Let $X,Y,Z,W \in SmProj(k)$. Let $\alpha \in \mathbf{E}^{2(n_Y+d_Y),(n_Y+d_Y)}(X_+ \wedge Th(V_Y))$, $\beta \in \mathbf{E}^{2(n_Z+d_Z),(n_Z+d_Z)}(Y_+ \wedge Th(V_Z))$ and $\gamma \in \mathbf{E}^{2(n_W+d_W),(n_W+d_W)}(Z_+ \wedge Th(V_W))$. Let us denote by $\beta \circ^{\dagger} \alpha$ for the composition of the above proposition and similarly for $\gamma \circ^{\dagger} \beta$. Then $\circ^{\dagger}$ is associative and unital. 
\end{prop}
\begin{proof}
We have that $\gamma \circ^{\dagger} (\beta \circ^{\dagger} \alpha)$ is the following composition by definition:
\begin{multline*}
\gamma \circ^{\dagger} (\beta \circ^{\dagger} \alpha) : \Sigma^{\infty}_{T,+}X \wedge Th(V_W) \stackrel{coev_Z}{\longrightarrow} \Sigma^{\infty}_{T,+} X \wedge Z \wedge Th(V_Z) \wedge S^{-2(n_Z+d_Z),-(n_Z+d_Z)} \wedge Th(V_W) \\ \stackrel{\tau}{\longrightarrow} \Sigma^{\infty}_{T,+}X \wedge Th(V_Z) \wedge Z \wedge Th(V_W) \wedge S^{-2(n_Z+d_Z),-(n_Z+d_Z)} \stackrel{(\beta \circ^{\dagger} \alpha) \wedge \gamma}{\longrightarrow} \mathbf{E} \wedge \mathbf{E} \wedge S^{2(n_W+d_W),(n_W+d_W)} \\ \stackrel{\mu_{\mathbf{E}}}{\longrightarrow} \mathbf{E} \wedge S^{2(n_W+d_W),(n_W+d_W)},
\end{multline*}
which can be rewritten as 
\begin{multline*}
\gamma \circ^{\dagger} (\beta \circ^{\dagger} \alpha): \Sigma^{\infty}_{T,+}X \wedge Th(V_W) \stackrel{coev_Z}{\longrightarrow} \Sigma^{\infty}_{T,+}X \wedge Z \wedge Th(V_Z) \wedge S^{-2(n_Z+d_Z),-(n_Z+d_Z)} \wedge Th(V_W) \\ \stackrel{\tau}{\longrightarrow} \Sigma^{\infty}_{T,+}X \wedge Th(V_Z) \wedge Z \wedge Th(V_W) \wedge S^{-2(n_Z+d_Z),-(n_Z+d_Z)} \stackrel{coev_Y\wedge \gamma}{\longrightarrow} \\ \Sigma^{\infty}_{T,+}X \wedge Y \wedge Th(V_Y) \wedge S^{-2(n_Y+d_Y),-(n_Y+d_Y)} \wedge Th(V_Z) \wedge \mathbf{E} \wedge S^{2(n_W+d_W)-2(n_Z+d_Z),(n_W+d_W)-(n_Z+d_Z)} \\ \stackrel{\tau}{\longrightarrow} \Sigma^{\infty}_{T,+}X \wedge Th(V_Y) \wedge Y \wedge Th(V_Z) \wedge \mathbf{E} \wedge S^{2(n_W+d_W)-2(n_Y+d_Y)-2(n_Z+d_Z),(n_W+d_W)-(n_Z+d_Z)-(n_Y+d_Y)} \\ \stackrel{\alpha \wedge \beta}{\longrightarrow} \mathbf{E} \wedge \mathbf{E} \wedge \mathbf{E} \wedge S^{2(n_W+d_W),(n_W+d_W)} \stackrel{\mu_{\mathbf{E}} \wedge \id_{\mathbf{E}}}{\longrightarrow} \mathbf{E} \wedge \mathbf{E} \wedge S^{2(n_W+d_W),(n_W+d_W)} \stackrel{\mu_{\mathbf{E}}}{\longrightarrow} \mathbf{E} \wedge S^{2(n_W+d_W),(n_W+d_W)}.
\end{multline*}
The composition $(\gamma \circ^{\dagger} \beta) \circ^{\dagger} \alpha$ is by definition:  
\begin{multline*}
(\gamma \circ^{\dagger} \beta) \circ^{\dagger} \alpha: \Sigma^{\infty}_{T,+} X \wedge Th(V_W) \stackrel{coev_Y}{\longrightarrow} \Sigma^{\infty}_{T,+}X \wedge Y \wedge Th(V_Y) \wedge S^{-2(n_Y+d_Y),-(n_Y+d_Y)} \wedge Th(V_W) \\ \stackrel{\tau}{\rightarrow} \Sigma^{\infty}_{T,+}X \wedge Th(V_Y) \wedge Y \wedge Th(V_W) \wedge S^{-2(n_Y+d_Y),-(n_Y+d_Y)} \stackrel{\alpha \wedge (\gamma \circ^{\dagger} \beta)}{\longrightarrow} \mathbf{E}\wedge \mathbf{E} \wedge S^{2(n_W+d_W),(n_W+d_W)} \\ \stackrel{\mu_{\mathbf{E}}}{\longrightarrow} \mathbf{E} \wedge S^{2(n_W+d_W),(n_W+d_W)},
\end{multline*}
which can be rewritten as
\begin{multline*}
(\gamma \circ^{\dagger} \beta) \circ^{\dagger} \alpha: \Sigma^{\infty}_{T,+} X \wedge Th(V_W) \stackrel{coev_Y}{\longrightarrow} \Sigma^{\infty}_{T,+}X \wedge Y \wedge Th(V_Y) \wedge S^{-2(n_Y+d_Y),-(n_Y+d_Y)} \wedge Th(V_W) \\ \stackrel{\tau}{\rightarrow} \Sigma^{\infty}_{T,+}X \wedge Th(V_Y) \wedge Y \wedge Th(V_W) \wedge S^{-2(n_Y+d_Y),-(n_Y+d_Y)} \stackrel{\alpha \wedge coev_Z}{\longrightarrow} \\ \mathbf{E} \wedge \Sigma^{\infty}_{T,+}Y \wedge Z \wedge Th(V_Z) \wedge S^{-2(n_Z+d_Z),-(n_Z+d_Z)} \wedge Th(V_W) \stackrel{\tau}{\longrightarrow} \\ \Sigma^{\infty}_{T,+}Y \wedge Th(V_Z) \wedge Z \wedge Th(V_W) \wedge S^{-2(n_Z+d_Z),(n_Z+d_Z)} \wedge \mathbf{E} \stackrel{\gamma \wedge \beta}{\longrightarrow} \mathbf{E} \wedge \mathbf{E} \wedge \mathbf{E} \wedge S^{2(n_W+d_W),(n_W+d_W)} \\ \stackrel{\mu_{\mathbf{E}} \wedge \id_{\mathbf{E}}}{\longrightarrow} \mathbf{E} \wedge \mathbf{E} \wedge S^{2(n_W+d_W),(n_W+d_W)} \stackrel{\mu_{\mathbf{E}}}{\longrightarrow}\mathbf{E} \wedge S^{2(n_W+d_W),(n_W+d_W)}. 
\end{multline*}
Both $\gamma \circ^{\dagger}(\beta \circ^{\dagger}\alpha)$ and $(\gamma \circ^{\dagger} \beta) \circ^{\dagger} \alpha$ are equal to the following composition
\begin{multline*}
\Sigma^{\infty}_{T,+}X \wedge Th(V_W) \stackrel{coev_Y \wedge coev_Z}{\longrightarrow} \\ \Sigma^{\infty}_{T,+}X \wedge Y \wedge Th(V_Y) \wedge S^{-2(n_Y+d_Y),-(n_Y+d_Y)} \wedge Z \wedge Th(V_Z) \wedge S^{-2(n_Z+d_Z),-(n_Z+d_Z)} \wedge Th(V_W) \\ \stackrel{\tau}{\longrightarrow} \Sigma^{\infty}_{T,+}X \wedge Th(V_Y) \wedge Y \wedge Th(V_Z) \wedge Z \wedge Th(V_W) \wedge S^{-2(n_Y+d_Y+n_Z+d_Z),-(n_Y+d_Y+n_Z+d_Z)} \\ \stackrel{\alpha \wedge \beta \wedge \gamma}{\longrightarrow} \mathbf{E} \wedge \mathbf{E} \wedge \mathbf{E} \wedge S^{2(n_W+d_W),(n_W+d_W)} \stackrel{\mu_{\mathbf{E}} \wedge \id_{\mathbf{E}}}{\longrightarrow} \mathbf{E} \wedge \mathbf{E} \wedge S^{2(n_W+d_W),(n_W+d_W)} \\ \stackrel{\mu_{\mathbf{E}}}{\longrightarrow} \mathbf{E} \wedge S^{2(n_W+d_W),(n_W+d_W)}.
\end{multline*}
\end{proof}
\begin{defn}{\rm
Let $\mathbf{E} \in SH(k)$ be a motivic ring spectrum. We define the category of Thom-$\mathbf{E}$-correspondences $Corr_{\mathbf{E}}(k)^{\dagger}$ to be the category, whose objects are 
\[Obj(Corr_{\mathbf{E}}(k)^{\dagger} = Obj(SmProj(k))\]
and morphisms are given by 
\[Corr_{\mathbf{E}}(k)^{\dagger}(X,Y) = \mathbf{E}^{2(n_Y+d_Y),n_Y+d_Y}(X_+ \wedge Th(V_Y)),\]
where $V_Y/Y$ is the duality vector bundle of rank $n_Y$. Given
\[\alpha \in \mathbf{E}^{2(n_Y+d_Y),(n_Y+d_Y)}(X_+ \wedge Th(V_Y))\]
and 
\[\beta \in \mathbf{E}^{2(n_Z+d_Z),(n_Z+d_Z)}(Y_+ \wedge Th(V_Z)),\] 
we define their composition to be 
\begin{multline*}\beta \circ^{\dagger} \alpha: \Sigma^{\infty}_{T,+}X \wedge Th(V_Z) \stackrel{coev_Y}{\longrightarrow} \Sigma^{\infty}_{T,+}X \wedge Y \wedge Th(V_Y) \wedge S^{-2(n_Y+d_Y),-(n_Y+d_Y)} \wedge Th(V_Z) \\ \stackrel{\tau}{\rightarrow} \Sigma^{\infty}_{T,+} X \wedge Th(V_Y) \wedge Y \wedge Th(V_Z) \wedge S^{-2(n_Y+d_Y),-(n_Y+d_Y)} \stackrel{\alpha \wedge \beta}{\longrightarrow} \mathbf{E} \wedge^{\mathbb{L}}_S \mathbf{E} \wedge S^{2(n_Z+d_Z),(n_Z+d_Z)} \\ \stackrel{\mu_{\mathbf{E}}}{\longrightarrow} \mathbf{E} \wedge S^{2(n_Z+d_Z),(n_Z+d_Z)},   
\end{multline*}
where 
\[coev_Y: \mathbb{S}^0 \rightarrow \Sigma^{\infty}_{T,+} Y \wedge Th(V_Y) \wedge S^{-2(n_Y+d_Y),-(n_Y+d_Y)}\]
is the coevaluation map of the Atiyah-Spanier-Whitehead duality on $Y$.
}
\end{defn}
As we may write $X_+ \wedge Th(V_Y) = Th(pr_{Y}^{XY*}V_Y)$, we have then the pullback map 
\[pr_{XY}^{XYZ*}: \mathbf{E}^{2(n_Y+d_Y),(n_Y+d_Y)}(Th(pr_{Y}^{XY*}V_Y)) \rightarrow \mathbf{E}^{2(n_Y+d_Y),(n_Y+d_Y)}(Th(pr_{XY}^{XYZ*}V_Y)).\]
Similarly
\[pr_{YZ}^{XYZ*} : \mathbf{E}^{2(n_Z+d_Z),(n_Z+d_Z)}(Th(pr_Z^{YZ*}V_Z)) \rightarrow \mathbf{E}^{2(n_Z+d_Z),(n_Z+d_Z)}(Th(pr_{YZ}^{XYZ*}V_Z)).\]
By taking cup product 
\begin{multline*}
-\cup_{\mathbf{E}}- : \mathbf{E}^{2(n_Y+d_Y),(n_Y+d_Y)}(Th(pr_{XY}^{XYZ*}V_Y)) \otimes \mathbf{E}^{2(n_Z+d_Z),(n_Z+d_Z)}(Th(pr_{YZ}^{XYZ*}V_Z)) \rightarrow \\ \mathbf{E}^{2(n_Y+d_Y+n_Z+d_Z),(n_Y+d_Y+n_Z+d_Z)}(Th(pr_{XY}^{XYZ*}V_Y) \wedge Th(pr_{YZ}^{XYZ*}V_Z)),
\end{multline*}
and applying the pushforward $pr_{XZ}^{XYZ*} = (coev_Y)^*$ we see that 
\begin{multline*}
pr_{XZ*}^{XYZ}(pr_{XY}^{XYZ*}\alpha \cup_{\mathbf{E}} pr_{YZ}^{XYZ*}\beta) \in \mathbf{E}^{2(n_Z+d_Z),(n_Z+d_Z)}(Th(pr_{Z}^{XZ*}V_Z)) = \\ \mathbf{E}^{2(n_Z+d_Z),(n_Z+d_Z)}(X_+ \wedge Th(V_Z)).
\end{multline*}
\begin{prop}\label{composecompare}
Let $\mathbf{E} \in SH(k)$ be a motivic ring spectrum. Let $X,Y,Z \in SmProj(k)$ of dimension $d_X,d_Y,d_Z$ respectively. Let $\alpha \in \mathbf{E}^{2(n_Y+d_Y),(n_Y+d_Y)}(X_+ \wedge Th(V_Y))$ and $\beta \in E^{2(n_Z+d_Z),(n_Z+d_Z)}(Y_+ \wedge Th(V_Z))$. Then the composition $\beta \circ^{\dagger} \alpha$ in $Corr_{\mathbf{E}}(k)^{\dagger}$ satisfies 
\[\beta \circ^{\dagger} \alpha = pr_{XZ*}^{XYZ}(pr_{XY}^{XYZ*}\alpha \cup_{\mathbf{E}} pr_{YZ}^{XYZ*}\beta),\]
where $pr_{XZ*}^{XYZ} = coev_Y^*$. 
\end{prop}
\begin{proof}
Trivial. 
\end{proof}
\begin{thm}(Comparison)
Let $\mathbf{E} \in SH(k)$ be a motivic ring spectrum. There is an equivalence of categories up to a natural $2$-isomorphism
\[\widetilde{Corr}_{\mathbf{E}}(k) \stackrel{\simeq}{\rightarrow} Corr_{\mathbf{E}}(k)^{\dagger} \stackrel{\simeq}{\rightarrow} Corr_{\mathbf{E}}(k).  \]
\end{thm}
\begin{proof}
We have the following association 
\[\widetilde{Corr}_{\mathbf{E}}(k) \rightarrow Corr_{\mathbf{E}}(k)^{\dagger} \rightarrow Corr_{\mathbf{E}}(k), \quad X \mapsto X \mapsto X\]
and 
\begin{multline*}
\widetilde{Corr}_{\mathbf{E}}(k)(X,Y) \stackrel{defn}{=} \mathbf{E}^{2d_Y,d_Y}(X\times Y,pr_Y^{XY*}V_Y) \stackrel{th_\mathbf{E}}{\cong} \mathbf{E}^{2(d_Y+n_Y),d_Y+n_Y}(X_+ \wedge Th(V_Y))  \stackrel{defn}{=} \\ = SH(k)[\Sigma^{\infty}_{T,+}X \wedge Th(V_Y),\mathbf{E}^{2(d_Y+n_Y),d_Y+n_Y}] = Corr_{\mathbf{E}}(k)(X,Y)^{\dagger} \\ \stackrel{D}{\cong} SH(k)[\Sigma^{\infty}_{T,+}X,\Sigma^{\infty}_{T,+}Y \wedge^{\mathbb{L}}_S \mathbf{E}] \stackrel{defn}{=} Corr_{\mathbf{E}}(k)(X,Y),
\end{multline*}
where $th_{\mathbf{E}}$ denotes the twisted Thom isomorphism and $D$ is the isomorphism induced by duality. It remains to check that the composition law in $\widetilde{Corr}_{\mathbf{E}}(k)$ is compatible with the composition law in $Corr_{\mathbf{E}}(k)^{\dagger}$ and $Corr_{\mathbf{E}}(k)$ via $th_{\mathbf{E}}$ and $D$ respectively. The compatibility of composition laws via the Thom isomorphism $th_{\mathbf{E}}$ follows from the Propositions \ref{propthEpullback}, \ref{propthEpushforward} and \ref{composecompare}. Now given a cohomology class $\alpha \in \mathbf{E}^{2(d_Y+n_Y),(d_Y+n_Y)}(Th_{XY}(pr_Y^{XY*}V_Y))$
we obtain its pullback by 
\[\xymatrix{\Sigma^{\infty}_{T,+}Th_{XY}(pr_Y^{XY*}V_Y) \ar[r]^{\alpha} & \mathbf{E}^{2(d_Y+n_Y),(d_Y+n_Y)} \\ \Sigma^{\infty}_{T,+}Th_{XYZ}(pr_Y^{XYZ*}V_Y) \ar[u] \ar[ru]_{pr_{XY}^{XYZ*}\alpha} }\]
Similarly, given  $\beta \in \mathbf{E}^{2(d_Z+n_Z),(d_Z+n_Z)}(Th_{YZ}(pr_Z^{YZ*}V_Z))$ we obtain its pullback by 
\[\xymatrix{\Sigma^{\infty}_{T,+}Th_{YZ}(pr_Z^{YZ*}V_Z) \ar[r]^{\beta} & E^{2(d_Z+n_Z),(d_Z+n_Z)} \\ \Sigma^{\infty}_{T,+}Th_{XYZ}(pr_Z^{XYZ*}V_Z) \ar[u] \ar[ru]_{pr_{YZ}^{XYZ*}\beta} } \]
So $pr_{XY}^{XYZ*}\alpha \cup_{\mathbf{E}} pr_{YZ}^{XYZ}\beta$ is the following composition 
\begin{multline*}
pr_{XY}^{XYZ*}\alpha \cup_E pr_{YZ}^{XYZ}\beta: \Sigma^{\infty}_{T,+}Th_{XYZ}(pr_Y^{XYZ*}V_Y) \wedge \Sigma^{\infty}_{T,+}Th_{XYZ}(pr_Z^{XYZ*}V_Z) \stackrel{-\wedge-}{\longrightarrow} \\ \mathbf{E} \wedge \mathbf{E} \wedge S^{2(d_Y+d_Z+n_Y+n_Z),d_Y+d_Z+n_Y+n_Z} \stackrel{\mu_\mathbf{E}}{\rightarrow} \mathbf{E} \wedge S^{2(d_Y+d_Z+n_Y+n_Z),d_Y+d_Z+n_Y+n_Z},
\end{multline*}
which corresponds to the morphism 
\[\Sigma^{\infty}_{T,+}Th_{XYZ}(-pr_Y^{XYZ*}T_Y) \wedge \Sigma^{\infty}_{T,+}Th_{XYZ}(-pr_Z^{XYZ*}T_Z) \rightarrow \mathbf{E}.\]
By definition the composition $\beta \circ \alpha$ as composition of $\mathbf{E}$-correspondences is given by the composition 
\[\Sigma^{\infty}_{T,+}Th_{XZ}(-pr_Z^{XZ*}T_Z) \rightarrow \Sigma^{\infty}_{T,+}Th_{XYZ}(-pr_Y^{XYZ*}T_Y) \wedge \Sigma^{\infty}_{T,+}Th_{XYZ}(-pr_Z^{XYZ*}T_Z) \rightarrow \mathbf{E},\]
where the first map by construction is given as 
\[\xymatrix{\Sigma^{\infty}_{T,+}Th_{XZ}(-pr_Z^{XZ*}T_Z) \ar[r] \ar@{=}[d] &  \Sigma^{\infty}_{T,+}Th_{XYZ}(-pr_Y^{XYZ*}T_Y) \wedge \Sigma^{\infty}_{T,+}Th_{XYZ}(-pr_Z^{XYZ*}T_Z) \ar@{=}[d] \\ \Sigma^{\infty}_{T,+} X \wedge (\Sigma^{\infty}_{T,+}Z)^{\vee} \ar[r]_{\id_X \wedge coev_Y \wedge \id_{Z^{\vee}}\quad \quad \quad \quad} & \Sigma^{\infty}_{T,+}X \wedge \Sigma^{\infty}_{T,+}Y \wedge (\Sigma^{\infty}_{T,+}Y)^{\vee} \wedge (\Sigma^{\infty}_{T,+}Z)^{\vee} }\] 
This implies that the composition $\beta \circ \alpha$ as $\mathbf{E}$-correspondences is the same as 
\[\xymatrix{X \wedge Y^{\vee} \wedge Y \wedge Z \ar[r] & \mathbf{E} \wedge \mathbf{E} \ar[r]^{\mu_\mathbf{E}} & \mathbf{E} \\ X \wedge Z^{\vee} \ar[u]_{coev_Y} \ar[urr]  }\]
This shows that the composition laws of $\widetilde{Corr}_{\mathbf{E}}(k)$ and $Corr_{\mathbf{E}}(k)$ are compatible. 
\end{proof}
\section{Proof of theorem \ref{mainthm}}
\subsection{Homotopy $t$-structure}
We recall in this section the notion homotopy $t$-structure in terms of generators (cf. \cite{Ay08}). Let $k$ be a field. The subcategory $SH(k)_{\geq n}$ is generated under homotopy colimits and extensions by 
\[\{S^{p,q} \wedge \Sigma^{\infty}_{\P^1}(X_+) | X \in Sm/k, p-q \geq n \},\]
where $S^{p,q} = S^{p-q}_s \wedge S_t^{q}$ denotes the motivic spheres. We set 
\[SH(k)_{\leq n} = \{ E \in SH(k)| [F,E] = 0, \forall F \in SH(k)_{\geq n+1} \}\]
The bigraded motivic homotopy sheaves are defined as 
\[\underline{\pi}_{p,q}^{st\A^1}(E) = a_{Nis}(U \mapsto SH(k)[S^{p,q}\wedge \Sigma^{\infty}_{\P^1}(U_+),E].\]
We let 
\[\underline{\pi}^{st\A^1}_p(E)_n \stackrel{defn}{=} a_{Nis}(U \mapsto SH(k)[\Sigma^{\infty}_{\P^1}U_+,S^{n-p,n} \wedge E]) \]
For a fix $p \in \Z$, $\pi^{st\A^1}_p(E)_*$ is considered as an abelian $\Z$-graded sheaf. An abelian Nisnevich sheaf $F \in Sh_{Nis}(Sm/k)$ is called strictly $\A^1$-invariant, if the map induced by the projection $U \times \A^1 \rightarrow U$: 
\[H^i_{Nis}(U,F) \rightarrow H^i_{Nis}(U \times \A^1,F)\]
is an isomorphism $\forall U \in Sm/k$ and $\forall i \geq 0$. For an abelian Nisnevich sheaf $F \in Sh_{Nis}(Sm/k)$ we will denote by 
\[F_{-1}(U) = \Ker(F(U\times_k \G_m) \rightarrow F(X)),\]
where the map is induced by the unit section of $\G_m$.
\begin{defn}(Morel). {\rm A homotopy module is a pair $(F_*,\varepsilon_*)$, where $F$ is a strictly $\A^1$-invariant $\Z$-graded abelian Nisnevich sheaf with 
\[\varepsilon_n : F_n \stackrel{\cong}{\longrightarrow} (F_{n+1})_{-1}. \]
}
\end{defn}    
The following description of the homotopy $t$-structure is a consequence of F. Morel's stable $\A^1$-connectivity result (see for instance \cite{Mor04a}): 
\begin{thm}[F. Morel]
Let $k$ be field. 
\begin{enumerate}
\item The triple $(SH(k),SH(k)_{\geq 0},SH(k)_{\leq 0})$ is a $t$-structure on $SH(k)$.
\item The heart of the homotopy $t$-structure $\pi_*^{\A^1}(k) = SH(k)_{\geq 0} \cap SH(k)_{\leq 0}$ is identified with the category of homotopy modules. 
\item The homotopy $t$-structure is non-degenerated in the sense that for any $U \in Sm/k$ and any $E \in SH(k)$, one has the morphism 
\[[\Sigma^{\infty}_{\P^1}(U_+),E_{\geq n}] \rightarrow [\Sigma^{\infty}_{\P^1}(U_+),E]\]
is an isomorphism for $n \leq 0$ and the morphism 
\[[\Sigma^{\infty}_{\P^1}(U_+),E] \rightarrow [\Sigma^{\infty}_{\P^1}(U_+),E_{\leq n}]\] 
is an isomorphism for $n > \dim (U)$.
\end{enumerate}
\end{thm}
By sending $E \mapsto E_{\geq n}$ and $E \mapsto E_{\leq n-1}$ respectively, one has the following adjunctions respectively: 
\[i_{\geq n} : SH(k)_{\geq n} \leftrightarrows SH(k) : \tau_{\geq n}, \quad \tau_{\leq n-1} : SH(k) \leftrightarrows SH(k)_{\leq n-1} : i_{\leq n -1},\]
where we denote by $i_{\geq n}$ and $i_{\leq n -1}$ the inclusion functors. We denote by
\[H : \pi_*^{\A^1}(k) \rightarrow SH(k)\]
the inclusion functor. For a homotopy module $F_* \in \pi_*^{\A^1}(k)$ we will call $H(F_*)$ the Eilenberg-Maclane spectrum associated to $F_*$. Let $S$ be now a Noetherian scheme of finite Krull dimension. We recall the rationally splitting of $SH(S)_{\Q}$ constructed by F. Morel (see \cite[\S 16.2]{CD10}). The permutation isomorphism 
\[\tau: \Sigma^{\infty}_{\P^1,+} \G_{m,\Q} \wedge \Sigma^{\infty}_{\P^1,+}\G_{m,\Q} \rightarrow  \Sigma^{\infty}_{\P^1,+} \G_{m,\Q} \wedge \Sigma^{\infty}_{\P^1,+}\G_{m,\Q}\]
satisfies $\tau^2 = 1$. This defines an element $e \in End_{SH(S)_{\Q}}(\mathbf{1}_{\Q})$, such that $e^2 = 1$. So we may define 
\[e_+ = \frac{e-1}{2}, \quad e_- = \frac{e+1}{2}.\]
Remark that $e_+$ and $e_-$ are idempotents. Hence we can define $\mathbf{1}_{\Q+} = \im(e_+)$ and $\mathbf{1}_{\Q-} = \im(e_-)$. For any spectrum $\sE \in SH(S)_{\Q}$, one defines $\sE_+ = \mathbf{1}_{\Q+} \wedge \sE$ and $\sE_- = \mathbf{1}_{\Q-} \wedge \sE$. This leads to a splitting of stable homotopy category 
\[SH(S)_{\Q+} \times SH(S)_{\Q-} \stackrel{\cong}{\longrightarrow} SH(S)_{\Q}, \quad (\sE_+,\sE_-) \mapsto \sE_+ \wedge \sE_-\]
Let us assume now $S = \Spec k$. The algebraic Hopf fibration is the map 
\[\A^2_k - \{0\} \rightarrow \P^1_k, \quad (x,y) \mapsto [x:y].\]
This gives us the stable Hopf map in $SH(k)$
\[\eta : \Sigma^{\infty}_{T,+}\G_m \rightarrow \mathbb{S}^0_k.\] 
Remark that from \cite[6.2.1]{Mor04a} one has a homotopy fiber sequence in $SH(k)$:
\[\Sigma^{\infty}_{T,+}(\A^2_k-\{0\}) \stackrel{S^{2,1}\wedge \eta}{\longrightarrow} \Sigma^{\infty}_{T,+}\P_k^1 \stackrel{\Sigma^{\infty}_{T,+}(i)}{\longrightarrow} \Sigma^{\infty}_{T,+}\P_k^2,\]
where $i: \P^1_k \inj \P^2_k$ is the linear embedding. Following \cite{Mor12} we define the Milnor-Witt $K$-theory of a field $F$ without any assumption on $char(F)$: 
\begin{defn}{\rm
Let $F$ be a field. $K^{MW}_*(F)$ is the $\Z$-graded associative unital ring freely generated by the symbols $[u]$, where $u \in F^{\times}$ is of degree $1$ and a symbol $\eta$ of degree $-1$ subject to the relation 
\begin{enumerate}
\item $[u]\cdot[1-u] = 0, \forall u \in F^{\times} - \{1\}.$
\item $[uv] = [u]+[v] + \eta \cdot [u] \cdot [v], \forall (u,v) \in (F^{\times})^2.$
\item $\eta \cdot [u] = [u] \cdot \eta, \forall u \in F^{\times}$. 
\item Define $h \stackrel{defn}{=} \eta \cdot [-1] + 2$. Then $\eta \cdot h = 0$. 
\end{enumerate} 
}
\end{defn}
Let $GW(F)$ be the Grothendieck-Witt ring of non-degenerate bilinear symmetric forms over $F$, where addition is given by orthogonal sum $\oplus$ and multiplication is given by tensor product $\otimes$. There is a surjective ring homomorphism 
\[rk: GW(F) \surj \Z, \quad Q \mapsto rk(Q).\]
The fundamental ideal is defined as 
\[I(F) \stackrel{defn}{=} \Ker (rk: GW(F) \surj \Z).\]
Denote by $I^n(F)$ the $n$-th power of $I(F)$. If $n \leq 0$ one sets $I^n(F) = W(F)$, where $W(F)$ is the Witt ring over $F$. Remark that $W(F) = GW(F)/(h)$, where $(h)$ is the ideal generated by hyperbolic spaces. By \cite[Lem. 3.10]{Mor12} there is a ring isomorphism 
\[GW(F) \stackrel{\cong}{\longrightarrow} K_0^{MW}(F), \quad \langle u \rangle \mapsto 1 + \eta \cdot [u].\]
Let $K^M_*(F)$ be the Milnor $K$-theory
\[K^M_*(F) \stackrel{defn}{=} Tens^*(F^{\times})/\langle u \otimes (1-u) \rangle.\]
There is a graded surjective homomorphism 
\[U : K^{MW}_*(F) \surj K^M_*(F), \quad [u] \mapsto \{u\}, \, \eta \mapsto 0.\]
In fact, one can show that for each $n$ there is a pullback diagram 
\[\xymatrix{K_n^{MW}(F) \ar[d]_U \ar[r] & I^n(F) \ar[d] \\ K^M_n(F) \ar[r] & I^n(F)/I^{n+1}(F)}\]
Following \cite[\S 3.2]{Mor12} we let $\underline{\mathbf{K}}^{MW}_n$ be the $n$-th Milnor-Witt sheaf, which is a strictly $\A^1$-invariant sheaf on $(Sm/k)_{Nis}$. In \cite[p. 437]{Mor04a} Morel showed that one can define a homotopy module $\underline{\mathbf{K}}^{MW}_*$ asscociated to the Milnor-Witt $K$-theory and in fact one has an isomorphism between homotopy modules 
\[\pi^{st\A^1}_0(\mathbb{S}^0)_* \cong \underline{\mathbf{K}}^{MW}_*. \]
The homotopy module $\underline{\mathbf{W}}_*$ is defined by setting every terms to be the unramified Witt sheaf $\mathbf{W} = a_{Nis}(U \mapsto W(U) = W(k(U)))$ and all the maps $\varepsilon_n$ are identity.  
\begin{lem}\label{corHKW}
Let $k$ be a field. Let $H\underline{\mathbf{K}}^{MW}_{*,\Q}$ be the Eilenberg-Maclane spectrum associated to the Milnor-Witt $K$-theory homotopy module $\underline{\mathbf{K}}^{MW}_{*,\Q}$. There exists a strict motivic ring spectrum $\hat{H}\underline{\mathbf{K}}^{MW}_{*,\Q} \in \mathbf{Spect}^{\Sigma}_{T}(k)_{\Q}$, which is isomorphic to $H\underline{\mathbf{K}}^{MW}_{*,\Q}$ in $SH(k)_{\Q}$.   
\end{lem}
\begin{proof}
We have a splitting 
\[H\underline{\mathbf{K}}^{MW}_{*,\Q} = H(\underline{\mathbf{K}}^M_{*,\Q}) \vee H(\underline{\mathbf{W}}_{*,\Q}).\] 
The result of D\'eglise \cite[Cor. 4.1.7]{Deg13} asserts that $H(\underline{\mathbf{K}}^M_*)$ is a strict $H\Z$-module, where $H\Z$ denotes the motivic cohomology spectrum. The construction in \cite[\S 4]{ALP15} shows that the cofibrant replacement $H(\underline{\mathbf{W}}_{*,\Q})^{cof}$ is a commutative monoid object in $\mathbf{Spect}^{\Sigma}_T(k)_{\Q}$, which is isomorphic to $H(\underline{\mathbf{W}}_{*,\Q})$ and $\mathbb{S}^0_k[\eta^{-1}]$. These imply that
\[\hat{H}\underline{\mathbf{K}}^{MW}_{*,\Q} = H\underline{\mathbf{K}}_{*,\Q}^M \vee H(\underline{\mathbf{W}}_{*,\Q})^{cof}\]
is also a strict motivic ring spectrum in $\mathbf{Spect}^{\Sigma}_T(k)_{\Q}$, which is isomorphic to $H\underline{\mathbf{K}}^{MW}_{*,\Q}$ in $SH(k)_{\Q}$. 
\end{proof}
\begin{defn}{\rm 
Let $k$ be a field. We define the category of pure Chow-Witt motives to be 
\[CHW(k)_{\Q} = \mathbf{Mot}_{\hat{H}\underline{\mathbf{K}}^{MW}_{*,\Q}}(k)\]
}
\end{defn}
\begin{cor}\label{corCHW}
Let $k$ be a field. There is a functor 
\[CHW(k)_{\Q} \rightarrow SH(k)_{\Q}.\]
\end{cor}
\begin{proof}
This is a consequence of the Corollary \ref{corMotE} and Lemma \ref{corHKW}.
\end{proof}
\subsection{Isomorphism between $\Hom$-groups} 
Let $k$ be a field. In this section we prove that one has a fully faithful embedding 
\[CHW(k)_{\Q} \rightarrow D_{\A^1,gm}(k)_{\Q}.\]
Remark that one has the equivalences of categories: 
\[\mathbf{StHo}_{\A^1,S^1}(k)_{\Q} \cong D_{\A^1}^{eff}(k)_{\Q}, \quad \mathbf{StHo}_{\A^1,\P^1}(k)_{\Q} \cong D_{\A^1}(k)_{\Q}.\]
For $E \in SH(k)$ we define its stable $\A^1$-cohomology as 
\[H^{p,q}_{st\A^1}(E,\Z) = SH(k)(E,S^{p,q}).\]
We denote by $SH(k)_{\Q}$ the localization of $SH(k)$. One has an adjunction 
\[\mathbb{L} L_{\Q}: SH(k) \rightleftarrows SH(k)_{\Q} : \R U,\]
which is induced by the Quillen adjunction 
\[L_{\Q} : \mathbf{Spect}_T^{\Sigma}(k) \leftrightarrows \mathbf{Spect}_T^{\Sigma}(k)_{\Q} : U,\]
where $U: \mathbf{Spect}_T^{\Sigma}(k)_{\Q} \rightarrow \mathbf{Spect}_T^{\Sigma}(k)$ is the forgetful functor by considering 
\[E_{\Q} = E \wedge \mathbf{1}_{\Q} = E \wedge hocolim (\mathbb{S}^0 \stackrel{2}{\rightarrow} \mathbb{S}^0 \stackrel{3}{\rightarrow} \mathbb{S}^0 \stackrel{4}{\rightarrow} \cdots )\]
as a symmetric motivic $T$-spectrum in $\mathbf{Spect}_T^{\Sigma}(k)$. For a motivic spectrum $E \in SH(k)$ we define its rational stable $\A^1$-cohomology as 
\[H^{p,q}_{st\A^1}(E,\Q) = SH(k)_{\Q}(E_{\Q},S^{p,q}_{\Q}) = SH(k)(E_{\Q},S^{p,q}_{\Q}).\]
We remark that by \cite[Lem. B2]{Lev13} if $E$ is a compact object in $SH(k)$ then one has an isomorphism 
\[H^{p,q}_{st\A^1}(E,\Q) = H^{p,q}_{st\A^1}(E,\Z)\otimes \Q.\]
Similarly, we define the motivic cohomology of $E$ as $H^{p,q}_M(E,\Z) = SH(k)(E,H\Z \wedge S^{p,q})$. If $F_* \in \pi_0^{\A^1}(k)_*$ is a homotopy module, then the $HF_*$-cohomology of $E$ is defined as 
\[H(F_*)^{p,q}(E) = SH(k)(E,S^{p,q} \wedge HF_*)\]
and if $E = \Sigma^{\infty}_{T,+} \sX$, where $\sX \in Spc(k)_+$ is a $k$-space (eg. Thom spaces), then the later cohomology is $H^{p-q}_{Nis}(\sX,F_q)$, where this cohomology is defined as 
\[H^{p-q}_{Nis}(\sX,F_q) = \mathbf{Ho}_{\A^1,+}(k)[\sX,K(F_q)[p-q]],\] where $K(-)$ denotes the Eilenberg-Maclane functor. 
\begin{thm}\label{thmsecondapproach}
Let $k$ be a field and $E = \Sigma^{\infty}_{T,+}Th(V/X)$ be the Thom spectrum of a vector bundle $V$ on a smooth $k$-scheme $X$. Let $\mathbb{S}^0$ be the motivic sphere spectrum. There exists a canonical isomorphism 
\begin{equation*}
\varphi: H^{2p,p}_{st\A^1}(Th(V/X),\Q) \stackrel{\cong}{\longrightarrow} H^{p}_{Nis}(Th(V/X),\underline{\mathbf{K}}^{MW}_p)_{\Q}, 
\end{equation*}
where $\varphi$ is induced by the unit $\varphi_{MW}: \mathbb{S}^0 \rightarrow H\underline{\mathbf{K}}^{MW}_*$. 
\end{thm}
\begin{proof}
By stable $\A^1$-connectivity theorem of Morel \cite{Mor05} the motivic sphere spectrum $\mathbb{S}^0$ is $-1$-connective. So we have a distinguished triangle 
\[(\mathbb{S}^0)_{\geq 1} \rightarrow \mathbb{S}^0 \rightarrow H \pi_0(\mathbb{S}^0)_* \stackrel{+1}{\rightarrow}.\]
By the computation of Morel we have $\pi_0(\mathbb{S}^0)_* = \underline{\mathbf{K}}^{MW}_*$. So after smashing with $S^{2p,p}_{\Q}$ we obtain a distinguished triangle 
\[(\mathbb{S}^0)_{\geq 1} \wedge S^{2p,p}_{\Q} \rightarrow S^{2p,p}_{\Q} \rightarrow H\underline{\mathbf{K}}^{MW}_* \wedge S^{2p,p}_{\Q} \stackrel{+1}{\rightarrow}.\]
By taking $[Th(V/X),-]$ we have a long exact sequence 
\begin{multline*}
\cdots \rightarrow [Th(V/X),(\mathbb{S}^0)_{\geq 1} \wedge S^{2p,p}_{\Q}] \rightarrow [Th(V/X),S^{2p,p}_{\Q}] \stackrel{\varphi}{\rightarrow} [Th(V/X),H\underline{\mathbf{K}}^{MW}_* \wedge S^{2p,p}_{\Q}] \rightarrow \\ \rightarrow [Th(V/X),(\mathbb{S}^0)_{\geq 1} \wedge S^{2p+1,p}_{\Q}] \rightarrow \cdots 
\end{multline*}
Now we have $(\mathbb{S}^0)_{\geq 1} \wedge \mathbb{S}^0_{\Q} = (\mathbb{S}^0_{\Q})_{\geq 1}$. By the work of  C. D. Cisinski, F. D\'eglise (\cite{CD10}) and the work of  A. Ananyevskiy, M. Levine, I. Panin (\cite{ALP15}) we have 
\[\mathbb{S}^0_{\Q} = H\Q \vee H\underline{\mathbf{W}}_{*,\Q}.\] 
This implies $(\mathbb{S}^0_{\Q})_{\geq 1} = (H\Q)_{\geq 1}$. The motivic cohomology spectrum $H\Q$ is also $-1$-connective, so we have a distinguished triangle 
\[(H\Q)_{\geq 1} \rightarrow H\Q \rightarrow H \pi_0(H\Q)_* \stackrel{+1}{\rightarrow}.\]
The homotopy module $\pi_0(H\Q)_*$ is $\underline{\mathbf{K}}^M_{*,\Q}$. We have (by \cite[Cor. 19.2]{MVW06} and by purity)
\[H\Q^{2p,p}(Th(V/X)) \stackrel{\cong}{\longrightarrow} H\underline{\mathbf{K}}_*^{M,2p,p}(Th(V/X))_{\Q}.\]
Now from the splitting $\mathbb{S}^0_{\Q} = H\Q \vee H\underline{\mathbf{W}}_{*,\Q}$ the map $\varphi$ take the form:
\begin{multline*}
[Th(V/X),S^{2p,p}_{\Q}] \cong H\Q^{2p,p}(Th(V/X)) \oplus H^p_{Nis}(Th(V/X),\underline{\mathbf{W}}_{\Q}) \stackrel{\varphi}{\rightarrow} \\ H^{p}_{Nis}(Th(V/X),\underline{\mathbf{K}}^M_p)_{\Q} \oplus H^p_{Nis}(Th(V/X),\mathbf{W}_{\Q}) \cong [Th(V/X),H\underline{K}^{MW}_* \wedge S^{2p,p}_{\Q}]. 
\end{multline*}
This implies that $\varphi$ is a canonical isomorphism.  
\end{proof}
\begin{cor}\label{corsecondapproach}
Let $k$ be a field. The functor constructed in \ref{corCHW} 
\[CHW(k)_{\Q} \rightarrow SH(k)_{\Q}\]
is fully faithful 
\end{cor}
\begin{proof}
By definition $CHW(k)_{\Q}$ is the smallest pseudo-abelian full subcategory of the homotopy category $Ho_k(\hat{H}\underline{\mathbf{K}}^{MW}_{*,\Q}-Mod)$ generated as an additive category by 
\[\{\Sigma^{\infty}_{T,+}X \wedge^{\mathbb{L}}_S \hat{H}\underline{\mathbf{K}}^{MW}_{*,\Q}| X \in SmProj(k)\}.\]
The adjunction 
\[-\wedge^{\mathbb{L}}_S \hat{H}\underline{\mathbf{K}}^{MW}_{*,\Q}: SH(k)_{\Q} \rightleftarrows Ho_k(\hat{H}\underline{\mathbf{K}}^{MW}_{*,\Q}-Mod): RU \] 
gives us a natural isomorphism
\begin{multline*}
Ho_k(\hat{H}\underline{\mathbf{K}}^{MW}_{*,\Q}-Mod)(\Sigma^{\infty}_{T,+} X \wedge^{\mathbb{L}}_S \hat{H}\underline{\mathbf{K}}^{MW}_{*,\Q}, \Sigma^{\infty}_{T,+} Y \wedge^{\mathbb{L}}_S \hat{H}\underline{\mathbf{K}}^{MW}_{*,\Q}]  \cong \\ SH(k)_{\Q}[\Sigma^{\infty}_{T,+}X,\Sigma^{\infty}_{T,+}Y \wedge^{\mathbb{L}}_S \hat{H}\underline{\mathbf{K}}^{MW}_{*,\Q}]. 
\end{multline*}
By duality \ref{refinement2} we have 
\[SH(k)_{\Q}[\Sigma^{\infty}_{T,+}X,\Sigma^{\infty}_{T,+}Y \wedge^{\mathbb{L}}_S \hat{H}\underline{\mathbf{K}}^{MW}_{*,\Q}] \cong H^{n_Y+d_Y}_{Nis}(Th(V_Y) \wedge X_+,\underline{\mathbf{K}}^{MW}_{n_Y+d_Y})_{\Q},\]
where $d_Y = \dim (Y)$, $V_Y$ is the duality vector bundle given in the theorem \ref{ThmVoeThom} and $n_Y = \rank (V_Y)$. 
The corollary follows now from the Theorem \ref{thmsecondapproach}.
\end{proof}
\section{Appendix}
In this appendix we simply recollect some facts and definitions in model categories. All the results are well-known and classical (see \cite{Q67}, \cite{Hir03}, \cite{Hov99}). 
\subsection{Model Categories}
\begin{defn}{\rm
A model category $\sM$ is a category with three classes of morphisms 
\[(Fib(\sM), Cof(\sM), W(\sM))\]
called fibrations, cofibrations and weak equivalences, such that: 
\begin{enumerate}
\item $\sM$ is closed under small limits and colimits.
\item If $f, g \in Mor(\sM)$ are composable and two out of $f, g, g\circ f$ are in $W(\sM)$, so is the third one.
\item Given a commutative diagram
\[\xymatrix{ A \ar[r] \ar@{>->}[d]_i & X \ar@{->>}[d]^p \\ B \ar[r] \ar@{-->}[ru] & Y }\] 
where $i \in Cof(\sM)$, $p \in Fib(\sM)$ and either $i$ or $p$ is in $W(\sM)$, then there exists a morphism $B \rightarrow X$ making the diagram commutative. 
\item $W(\sM), Cof(\sM)$ and $Fib(\sM)$ are closed under retracts.
\item Given any morphism $f : X \rightarrow Y$ in $Mor(\sM)$, there exist two functorial factorizations
\[\xymatrix{ & Z \ar@{->>}[rd] \\ X \ar@{>->}[ru]^{\simeq} \ar[rr]^f \ar@{>->}[rd] && Y \\ & W \ar@{->>}[ru]_{\simeq} }\] 
\end{enumerate}
}
\end{defn}
The first axiom implies that there exist an initial object $\varnothing$ and a final object $\star$. We say $\sM$ is pointed if $\varnothing \stackrel{\cong}{\longrightarrow} \star$.
\begin{defn}{\rm
Let $X \in Obj(\sM)$ be an object. $X$ is called cofibrant if the natural morphism $\varnothing \rightarrow X$ is in $Cof(\sM)$. $X$ is called fibrant if the natural morphism $X \rightarrow *$ is in $Fib(\sM)$.
} 
\end{defn}
Let $i: A \rightarrow B$ and $p: X \rightarrow Y$ be two morphisms in $Mor(\sM)$. We say $i$ has left lifting property wrt. $p$ or $p$ has right lifting property wrt. $i$, if for every solide commutative diagram 
\[\xymatrix{ A \ar[r] \ar@{>->}[d]_i & X \ar@{->>}[d]^p \\ B \ar[r] \ar@{-->}[ru] & Y }\] 
the dotted morphism exists and makes the diagram commutative. Given two morphisms $Mor(\sM) \ni f: A \rightarrow B$ and $Mor(\sM) \ni g: C \rightarrow D$, we say $f$ is a retract of $g$, if there is a commutative diagram
\[\xymatrix{A \ar[r] \ar[d]_f & C \ar[r] \ar[d]_g & A \ar[d]^f \\ B \ar[r] & D \ar[r] & B }\]
where the horizontal composites are identities. Given an object $X \in Obj(\sM)$, the factorization axiom tells us that we can factor 
\[\xymatrix{\varnothing \ar@{>->}[r] & X^{cof} \ar@{->>}[r]^{\simeq} & X,}\]  
where $X^{cof}$ is cofibrant. We call $X^{cof}$ a cofibrant replacement of $X$. Similarly, we can factor
\[\xymatrix{X \ar@{>->}[r]^{\simeq} & X^{fib} \ar@{->>}[r] & \star,}\]
where $X^{fib}$ is fibrant. We call $X^{fib}$ a fibrant replacement of $X$.
\begin{defn}{\rm
Let $\sM, \sN$ be two model categories. A functor $F: \sM \rightarrow \sN$ is called a left Quillen functor, if it has a right adjoint $G: \sN \rightarrow \sM$ and 
\begin{enumerate}
\item If $i \in Cof(\sM)$, then $F(i) \in Cof(\sN)$.
\item If $j \in Cof(\sM) \cap W(\sM)$, then $F(j) \in Cof(\sN) \cap W(\sM)$. 
\end{enumerate}
The right adjoint $G: \sN \rightarrow \sM$ is called a right Quillen functor and the adjunction 
\[F : \sM \leftrightarrows \sN : G\]
is called a Quillen adjunction. 
}
\end{defn}
\begin{defn}{\rm
Let 
\[F: \sM \leftrightarrows \sN : G\]
be a Quillen adjunction. $F$ is called a left Quillen equivalence, if for every cofibrant object $X \in Obj(\sM)$ and every fibrant object $Y \in Obj(\sN)$ one has the following: A morphism $f : X \rightarrow GY$ is in $W(\sM)$ iff its adjoint $g = \varepsilon_{(F,G)}\circ F(f) : FX \rightarrow Y$ is in $W(\sN)$. $G$ is called then a right Quillen equivalence. The adjunction 
\[F: \sM \leftrightarrows \sN : G\]
is called a Quillen equivalence. 
}
\end{defn}
\begin{defn}{\rm
Let $X \in Obj(\sM)$ be an object in a model category $\sM$. The cylinder object for $X$ is an object $Cyl(X)$, such that we have a factorization 
\[\xymatrix{X \coprod X \ar[rr]^{\nabla} \ar@{>->}[d]_i && X \\ Cyl(X) \ar[rru]^{\simeq}_s }\]
where $i \in Cof(\sM)$ and $s \in W(\sM)$.
}
\end{defn}
\begin{defn}{\rm
Let $X \in Obj(\sM)$ be an object in a model category $\sM$. A path object for $X$ is an object $\sP(X)$, such that we have a factorization
\[\xymatrix{X \ar[rr]^{\Delta} \ar[d]_r^{\simeq} && X \times X \\ \sP(X) \ar@{->>}[rru]_p }\]
where $r \in W(\sM)$ and $p \in Fib(\sM)$. 
}
\end{defn}
\begin{defn}{\rm
Let $f,g : X \rightrightarrows Y$ be two morphisms in $Mor(\sM)$ of a model category $\sM$. $f$ is left homotopic to $g$ if there is a cylinder object $Cyl(X)$ for $X$, such that we have a factorization 
\[\xymatrix{X \coprod X \ar@{>->}[d]_i \ar[rr]^{(f,g)} && Y \\ Cyl(X) \ar[rru]_{LH} }\]
The map $LH$ is called a left homotopy from $f$ to $g$. 
}
\end{defn}
\begin{defn}{\rm
Let $f,g : X \rightrightarrows Y$ be two morphisms in $Mor(\sM)$ of a model category $\sM$. $f$ is right homotopic to $g$ if there is a path object $\sP(Y)$ for $Y$, such that we have a factorization 
\[\xymatrix{X \ar[d]_{RH} \ar[rr]^{(f,g)} && Y \times Y \\ \sP(Y) \ar@{->>}[rru]_p }\]
The map $RH$ is called a right homotopy from $f$ to $g$.
}
\end{defn}
\begin{defn}
Let $f,g : X \rightrightarrows Y$ be two morphisms in $Mor(\sM)$ of a model category $\sM$. $f$ is homotopic to $g$ if $f$ is left and right homotopic to $g$.
\end{defn}
\begin{thm}(Quillen \cite[I.1 Thm. 1]{Q67}). Let $\sM$ be a model category. There exists a category $Ho(\sM) = \sM[W(\sM)^{-1}]$, which is called the homotopy category of $\sM$, where
\begin{enumerate}
\item $Obj(Ho(\sM)) = Obj(\sM)$.
\item $Ho(\sM)(X,Y) = \pi((X^{cof})^{fib},(Y^{cof})^{fib}),$, where $\pi$ denotes the set of homotopy classes and the composition law is induced by the composition law of $\sM$.  
\end{enumerate}
\end{thm}
\begin{thm}(Quillen \cite[I.4 Thm. 3]{Q67}). 
Let 
\[F : \sM \leftrightarrows \sN : G\]
be a Quillen adjunction. Then $(F,G)$ induces an adjunction of homotopy categories
\[\mathbb{L}F : Ho(\sM) \leftrightarrows Ho(\sN) : \mathbb{R}G.\]
\end{thm}
\begin{defn}{\rm
\begin{enumerate}
\item Let $\sM$ be a model category. $\sM$ is left proper, if in any pushout diagram
\[\xymatrix{A \ar@{>->}[d]_i \ar[r]^{h_1}_{\simeq} & X \ar[d] \\ B \ar[r]_{h_2} & Y }\] 
where $i \in Cof(\sM)$ and $h_1 \in W(\sM)$, so $h_2 \in W(\sM)$.
\item Let $\sM$ be a model category. $\sM$ is right proper, if in any pullback diagram 
\[\xymatrix{A \ar[d] \ar[r]^{h_1}_{\simeq} & X \ar@{->>}[d] \\ B \ar[r]_{h_2} & Y }\] 
where $p \in Fib(\sM)$ and $h_2\in W(\sM)$, so $h_1 \in W(\sM)$.
\item $\sM$ is proper, if it is left and right proper. 
\end{enumerate}
}
\end{defn}
Let $\mathbf{\Delta}$ denote the category, whose objects are ordered finite sets
\[ \underline{n} = \{0 < 1 < \cdots < n\}, n \geq 0\]
and
\[Mor(\mathbf{\Delta})(\underline{m},\underline{n}) = \{f : \underline{m} \rightarrow \underline{n}| i \leq j \implies f(i) \leq f(j)\}.\]
There are cofaces $\delta^i : \underline{n} \rightarrow \underline{n+1}$ and codegeneracies $\sigma^i : \underline{n+1} \rightarrow \underline{n}$ defined by 
\[\delta^i(j) = \begin{cases}j, & \mbox{if}\, \, j < i \\ j+1, & \mbox{if}\, \, j \geq i \end{cases} \]
\[\sigma^i(j) = \begin{cases} j, & \mbox{if} \, \, j \leq i \\ j-1, & \mbox{if} \, \, j > i \end{cases} \]
Cofaces and codegeneracies are generators for the maps in $\mathbf{\Delta}$. They satisfy a list of relations (cf. \cite[\S 8]{Weib94}). Now one defines the category of simplicial sets as
\[\mathbf{SSets} \stackrel{defn}{=} \mathbf{\Delta}^{op}(Sets).\]
So simplicial sets are just presheaves of sets on $\mathbf{\Delta}$. For a general category $\sA$ the category of simplicial objects and cosimplicial objects in $\sA$ are defined to be $\mathbf{\Delta}^{op}(\sA)$ and $\mathbf{\Delta}(\sA)$ respectively. Let $\mathbf{Top}$ be the category of compactly generated Hausdorff topological spaces. The geometric realization functor is defined by 
\[R : \mathbf{SSets} \rightarrow \mathbf{Top}, \quad X \mapsto R(X) = \int^{\underline{n}} X(\underline{n}) \times \Delta^n,\]
where $\Delta^n$ is the presheaf $Mor(\mathbf{\Delta})(-,\underline{n})$. There is an adjunction 
\[R : \mathbf{SSets} \leftrightarrows \mathbf{Top} : S,\]
where $S$ is the singular functor 
\[S(T): \mathbf{\Delta}^{op} \rightarrow Sets, \quad \underline{n} \mapsto \mathbf{Top}(R(\Delta^n),T).\]
Here $R(\Delta^n)$ is 
\[R(\Delta^n) = \{(x_0,\cdots,x_n) \in \R^{n+1}| x_i \geq 0, \sum_{i=0}^n x_i = 1\}.\]
\begin{thm}(Quillen \cite[II.3 Thm. 3]{Q67}). The category $\mathbf{SSets}$ has a model category structure. 
\end{thm}
\begin{defn}{\rm
Let $\sM$ be a category. $\sM$ is called simplicial if there is a functor 
\[\sM^{op} \times \sM \rightarrow \mathbf{SSets}, \quad (X,Y) \mapsto \mathbf{SSMap}(X,Y),\] 
such that
\begin{enumerate}
\item $\mathbf{SSMap}(X,Y)_0 = \sM(X,Y)$. 
\item there exists a composition law
\[\circ : \mathbf{SSMap}(Y,Z) \times \mathbf{SSMap}(X,Y) \rightarrow \mathbf{SSMap}(X,Z),\]
which is compatible with the composition law in $\sM$.
\item There is a simplicial sets map $i_X : \star \rightarrow \mathbf{SSMap}(X,X), \forall X \in Obj(\sM)$, where the associativity of the composition law, right and left unit properties of $i_X$ follows from three commutative diagrams (\cite[Def. 9.1.2]{Hir03}).  
\end{enumerate}
}
\end{defn}
\begin{defn}{\rm
Let $\sM$ be a model category. $\sM$ is called a simplicial model category if $\sM$ is simplicial and 
\begin{enumerate}
\item $\forall X \in Obj(\sM)$ there is an adjunction 
\[X \otimes - : \mathbf{SSets} \leftrightarrows \sM : \mathbf{SSMap}(X,-),\]
which is compatible with the simplicial structure on $\sM$.
\item $\forall Y \in Obj(\sM)$ there is an adjunction 
\[Y^{-} : \mathbf{SSets} \leftrightarrows \sM^{op} : \mathbf{SSMap}(-,Y),\]
which is compatible with the simplicial structure on $\sM$.
\item For $Cof(\sM) \ni i : A \rightarrow B$ and $Fib(\sM) \ni p : X \rightarrow Y$ the map 
\[\xymatrix{\mathbf{SSMap}(B,X) \ar@{->>}[rr]^{(i^*,p_*) \quad \quad \quad \quad \quad \quad \quad} && \mathbf{SSMap}(A,X) \times_{\mathbf{SSMap}(A,Y)} \mathbf{SSMap}(B,Y)}\]
is in $Fib(\mathbf{SSets}$, which is also in $W(\mathbf{SSets})$, if either $i$ or $p$ is in $W(\sM)$. 
\end{enumerate}
}
\end{defn}
\begin{ex}
$\mathbf{SSets}$ has a canonical simplicial model category structure. $\mathbf{SSMap}(X,Y)$ is the simplicial set with 
\[\mathbf{SSMap}(X,Y)_n = \mathbf{SSets}(X \times \Delta^n, Y),\]
with faces and degeneracies induced from the cosimplicial object $\Delta^{\bullet}$. 
\end{ex}
\begin{prop}
Let $\sM$ be a simplicial model category. If $X$ is cofibrant and $Y$ is fibrant, then 
\[Ho(\sM)(X,Y) = \pi_0\mathbf{SSMap}(X,Y).\]
Consequently, for any objects $A, B \in Obj(\sM)$ one has 
\[Ho(\sM)(A,B) = \pi_0\mathbf{SSMap}((A^{cof})^{fib},(B^{cof})^{fib}).\]
\end{prop}
\subsection{Localization}
All model categories in this subsection are being considered simplicial. 
\begin{defn}{\rm
Let $\sM$ be a model category and $\sV$ be a class of morphisms in $Mor(\sM)$. A left localization of $\sM$ wrt. $\sV$ is a model category $L_{\sV}\sM$ together with a left Quillen functor $F : \sM \rightarrow L_{\sV}\sM$, such that:
\begin{enumerate}
\item The total left derived functor $\mathbb{L}F: Ho(\sM) \rightarrow Ho(L_{\sV}\sM)$ takes the images in $Ho(\sM)$ of elements in $\sV$ into isomorphisms in $Ho(L_{\sV}\sM)$.
\item If $\sN$ is a model category and $T: \sM \rightarrow \sN$ is a left Quillen functor such that $\mathbb{L}T: Ho(\sM) \rightarrow Ho(\sN)$ take the images in $Ho(\sM)$ of elements in $\sV$ into isomorphisms in $Ho(\sN)$, then there is a unique left Quillen functor $L_{\sV}\sM \rightarrow \sN$, such that
\[\xymatrix{\sM \ar[r]^{F} \ar[d]_T & L_{\sV}\sM \ar@{-->}[ld]^{\exists!} \\ \sN }\]  
\end{enumerate}
}
\end{defn}
\begin{defn}{\rm
Let $\sM$ be a model category and $\sV$ be a class of morphisms in $Mor(\sM)$. A right localization of $\sM$ wrt. $\sV$ is a model category $R_{\sV}\sM$ together with a right Quillen functor $G: \sM \rightarrow R_{\sV}\sM$, such that:
\begin{enumerate}
\item The total right derived functor $\R G: Ho(\sM) \rightarrow Ho(R_{\sV}\sM)$ takes the images in $Ho(\sM)$ of elements in $\sV$ into isomorphisms in $Ho(R_{\sV}\sM)$
\item If $\sN$ is a model category and $T: \sM \rightarrow \sN$ is a right Quillen functor such that $\R T$ takes the images in $Ho(\sM)$ of elements in $\sV$ into isomorphisms in $Ho(\sN)$, then there is a unique right Quillen functor $R_{\sV}\sM \rightarrow \sM$, such that:
\[ \xymatrix{ \sM \ar[r]^G \ar[d]_T & R_{\sV}\sM \ar@{-->}[ld]^{\exists!} \\ \sN }\]
\end{enumerate}
}  
\end{defn}
\begin{defn}{\rm
Let $\sM$ be a model category and $\sV$ a class of morphisms in $Mor(\sM)$. 
\begin{enumerate}
\item An object $X \in Obj(\sM)$ is called $\sV$-local if $X$ is fibrant and for every $f: A \rightarrow B$ in $\sV$, $\mathbf{SSMap}(B^{cof},X) \stackrel{\simeq}{\longrightarrow} \mathbf{SSMap}(A^{cof},X)$. 
\item A morphism $f : X \rightarrow Y$ in $Mor(\sM)$ is a $\sV$-local equivalence if for every $\sV$-local object $T$, $\mathbf{SSMap}(Y^{cof},T) \stackrel{\simeq}{\longrightarrow} \mathbf{SSMap}(X^{cof},T)$
\item $X \in Obj(\sM)$ is called $\sV$-colocal if $X$ is cofibrant and for every $f : A \rightarrow B$ in $\sV$, $\mathbf{SSMap}(X,A^{fib}) \stackrel{\simeq}{\longrightarrow} \mathbf{SSMap}(X,B^{fib})$. 
\item $Mor(\sM) \ni f : X \rightarrow$ is a $\sV$-colocal equivalence if for every $\sV$-colocal object $T$, $\mathbf{SSMap}(T,X^{fib}) \stackrel{\simeq}{\longrightarrow} \mathbf{SSMap}(T,Y^{fib})$.
\end{enumerate}
}
\end{defn}
\begin{defn}{\rm
Let $\sM$ be a model category and $\sV$ be a class of morphisms in $Mor(\sM)$. The left Bousfield localization (if it exists) of $\sM$ wrt. $\sV$ is a model category structure $L_{\sV}\sM$ on the underlying category $\sM$ with:
\begin{enumerate}
\item $W(L_{\sV}\sM)$ is the class of $\sV$-local equivalences of $\sM$. 
\item $Cof(L_{\sV}\sM) = Cof(\sM)$.
\item $Fib(L_{\sV}\sM) = RLP(Cof(\sM) \cap W(L_{\sV}\sM))$. 
\end{enumerate}
}
\end{defn}
\begin{defn}{\rm
Let $\sM$ be a model category and $\sV$ be a class of morphisms in $Mor(\sM)$. The right Bousfield localization (if it exists) of $\sM$ wrt. $\sV$ is a model category structure $R_{\sV}\sM$ on the underlying category $\sM$ with:
\begin{enumerate}
\item $W(R_{\sV}\sM)$ is the class of $\sV$-colocal equivalences of $\sM$.
\item $Fib(R_{\sV}\sM) = Fib(\sM)$.
\item $Cof(R_{\sV}\sM) = LLP(Fib(\sM) \cap W(R_{\sV}\sM))$.
\end{enumerate}
}
\end{defn}
\subsection{Symmetric motivic $T$-Spectra}
The reference for this subsection is \cite{Jar00}. Let $S$ be a Noetherian scheme of finite Krull dimension. Consider the category $Sm/S$ of smooth of finite type $S$-schemes. A symmetric $T$-spectrum is a collection $\{X_n\}_{n \geq 0}$, where $X_n \in \mathbf{\Delta}^{op}(PrSh_{Nis}(Sm/S))_+$, together with the left actions 
\[\Sigma_n \times X_n \rightarrow X_n,\]
where $\Sigma_n$ is the $n$-th symmetric group. There are the bonding maps 
\[\sigma_n : T \wedge X_n \rightarrow X_{n+1},\]
such that the interative  composition 
\[T^{\wedge m} \wedge X_n \rightarrow X_{n+m}\]
is $\Sigma_m \times \Sigma_n$-equivariant. A morphism between symmetric $T$-spectra is a family $\{f_n: X_n \rightarrow Y_n \}_{n \geq 0}$, where the following diagram
\[\xymatrix{T \wedge X_n \ar[d]_{\sigma_n} \ar[r]^{\id \wedge f_n} & T \wedge Y_n \ar[d]^{\sigma_n} \\ X_{n+1} \ar[r]_{f_{n+1}} & Y_{n+1}}\]
commutes and $f_n$ is $\Sigma_n$-equivariant $\forall n \geq 0$. The category of symmetric $T$-spectra is denoted by $\mathbf{Spect}_T^{\Sigma}(S)$. A symmetric sequence $X$ is a family $\{X_n| X_n \in \mathbf{\Delta}^{op}(PrSh_{Nis}(Sm/S))_+\}_{n \geq 0}$ with left actions 
\[\Sigma_n \times X_n \rightarrow X_n.\]
A morphism $f : X \rightarrow Y$ of symmetric sequences is a family $\{f_n : X_n \rightarrow Y_n\}$, where $f_n$ are $\Sigma_n$-equivariant $\forall n \geq 0$. We denote the category of symmetric sequences of pointed simplicial presheaves by $\mathbf{\Delta}^{op}(PrSh_{Nis}(Sm/S))_+^{\Sigma}$. Recall that there are  families of functors 
\[F_n : \mathbf{\Delta}^{op}(PrSh_{Nis}(Sm/S))_+ \rightarrow \mathbf{\Delta}^{op}(PrSh_{Nis}(Sm/S))_+^{\Sigma},\]
where 
\[(F_n(\sX))_m = \begin{cases}\star & \mbox{if} \, \, m \neq n \\ \bigvee_{\sigma \in \Sigma_n} \sX & \mbox{if} \, \, m = n \end{cases} \]
and 
\[Ev_n : \mathbf{\Delta}^{op}(PrSh_{Nis}(Sm/S))_+^{\Sigma} \rightarrow \mathbf{\Delta}^{op}(PrSh_{Nis}(Sm/S))_+, \quad X \mapsto X_n.\]
They are in fact adjoint to each other
\[F_n : \mathbf{\Delta}^{op}(PrSh_{Nis}(Sm/S))_+ \leftrightarrows \mathbf{\Delta}^{op}(PrSh_{Nis}(Sm/S))_+^{\Sigma} : Ev_n.\]
For two symmetric sequences $X$ and $Y$, their product is defined as 
\[(X \otimes Y)_n = \bigvee_{p+q=n}\Sigma_n \otimes_{\Sigma_p \times \Sigma_q}X_p \wedge Y_q .\]
The notation $\Sigma_n \otimes_{\Sigma_p \times \Sigma_q} X_p \wedge Y_q$ means: there is an action $\gamma$ of $\Sigma_p \times \Sigma_q$ on $X_p \wedge Y_q$ via the canonical embedding $\Sigma_p \times \Sigma_q \subset \Sigma_n$ and also another action $\gamma': \Sigma_p\times \Sigma_q \times (X_p \wedge Y_q) \rightarrow X_p \wedge Y_q $. We let 
\[\Sigma_n \otimes_{\Sigma_p \times \Sigma_q} X_p \wedge Y_q = eq[\gamma_{\sigma} - \gamma'_{\sigma}]_{\sigma \in \Sigma_p \times \Sigma_q}.\]  
Now one can define 
\[F_n^{\Sigma} : \mathbf{\Delta}^{op}(PrSh_{Nis}(Sm/S))_+ \rightarrow \mathbf{Spect}_T^{\Sigma}(S), \quad \sX \mapsto \mathbb{S}^0 \otimes F_n(\sX),\]
where $\mathbb{S}^0$ denotes the motivic sphere spectrum
\[\mathbb{S}^0 = (S_+,T \wedge S_+, T^{\wedge 2} \wedge S_+, \cdots) \]
$\Sigma_n$ acts on $\mathbb{S}^0$ by permuting the $T^{\wedge n}$ factors and $S_+$ is pointed by $S \coprod S$. One has an adjunction 
\[F_n^{\Sigma} : \mathbf{\Delta}^{op}(PrSh_{Nis}(Sm/S))_+ \leftrightarrows \mathbf{Spect}_T^{\Sigma}(S) : Ev_n.\] 
In fact, one has $F^{\Sigma}_0(S_+) = \mathbb{S}^0$. A symmetric $T$-spectrum $X$ can be understood as a symmetric sequence with a module structure $\sigma_X: \mathbb{S}^0 \otimes X \rightarrow X$ over the motivic sphere spectrum $\mathbb{S}^0$. Now we can define the smash product of symmetric $T$-spectra as
\[\xymatrix{X \wedge Y \stackrel{defn}{=} coeq(\mathbb{S}^0 \otimes X \otimes Y \ar@<-.5ex>[rr] \ar@<.5ex>[rr] && X \otimes Y),}\]
where the top map is $\sigma_X \otimes \id_Y$ and the bottom map is 
\[\mathbb{S}^0 \otimes X \otimes Y \stackrel{\tau}{\longrightarrow} X \otimes \mathbb{S}^0 \otimes Y \stackrel{\id_X \otimes \sigma_Y}{\longrightarrow} X \otimes Y.\] 
We just mention the following results of Jardine.
\begin{thm}(Jardine \cite[Thm. 4.2]{Jar00})
The category $\mathbf{Spect}_{T}^{\Sigma}(S)$ has a model category structure, which is proper and simplicial. 
\end{thm}
\begin{thm}(Jardine \cite[Prop. 4.19]{Jar00}). 
$(\mathbf{Spect}_T^{\Sigma}(S),\mathbb{S}^0,\wedge)$ is a symmetric monoidal model category. 
\end{thm}
Now we discuss a little bit about the Quillen adjunction 
\[-\wedge \mathbf{E} : \mathbf{Spect}_T^{\Sigma}(S) \leftrightarrows \mathbf{E}-Mod^{\Sigma} : U,\]
where $\mathbf{E}$ is a motivic strict ring spectrum (we always consider only commutative ring spectrum). On the level of the underlying categories the unit and counit of the adjunction are defined by 
\[\eta_X : X \cong \mathbb{S}^0 \wedge X \stackrel{\varphi_{\mathbf{E}}\wedge \id_X}{\longrightarrow} \mathbf{E} \wedge X = U(\mathbf{E} \wedge X),\]
and
\[\varepsilon_M : \mathbf{E} \wedge U(M) = \mathbf{E} \wedge M \stackrel{\gamma_M}{\longrightarrow} M.\]
By \cite[Prop. 4.19]{Jar00} the category $\mathbf{Spect}_T^{\Sigma}(S)$ satisfies the axiom in \cite[Def. 3.3]{SS00}. By \cite[Thm. 4.1]{SS00} one can conclude that the adjunction 
\[- \wedge \mathbf{E} : \mathbf{Spect}_T^{\Sigma}(S) \leftrightarrows \mathbf{E}-Mod^{\Sigma} : U\]
induces a model category structure on $\mathbf{E}-Mod^{\Sigma}$. It is clear that the forgetful functor 
\[U : \mathbf{E}-Mod^{\Sigma} \rightarrow \mathbf{Spect}^{\Sigma}_T(S) \]
is a right Quillen functor, because $Fib(\mathbf{E}-Mod^{\Sigma})$ and $Fib(\mathbf{E}-Mod^{\Sigma}) \cap W(\mathbf{E}-Mod^{\Sigma})$ are detected in $\mathbf{Spect}^{\Sigma}_T(S)$. So we can claim that the adjunction above is a Quillen adjunction. Since $\mathbf{E}$ is a commutative ring spectrum, $\mathbf{E}-Mod^{\Sigma}$ has the closed symmetric monoidal category structure induced by the one on $\mathbf{Spect}_T^{\Sigma}(S)$ by declaring:
\[-\wedge_{\mathbf{E}}- : \mathbf{E}-Mod^{\Sigma} \times \mathbf{E}-Mod^{\Sigma} \rightarrow \mathbf{E}-Mod^{\Sigma}, \quad (M,N) \mapsto M \wedge_{\mathbf{E}} N  \]      
and
\[\underline{\Hom}_{\mathbf{E}-Mod^{\Sigma}} : \mathbf{E}-Mod^{\Sigma} \times \mathbf{E}-Mod^{\Sigma} \rightarrow \mathbf{E}-Mod^{\Sigma}, \quad \underline{\Hom}_{\mathbf{E}-Mod^{\Sigma}}(M,N), \] 
where 
\[\xymatrix{M \wedge_{\mathbf{E}} N \stackrel{defn}{=} coeq(\mathbf{E} \wedge M \wedge N \ar@<-.5ex>[rr] \ar@<.5ex>[rr] && M \wedge N)}.\]
The top map is $\gamma_M \wedge \id$ and the bottom map is the composition 
\[\mathbf{E} \wedge M \wedge N \stackrel{\tau \wedge \id}{\longrightarrow} M \wedge \mathbf{E} \wedge N \stackrel{\id \wedge \gamma_N}{\longrightarrow} M \wedge N.\]
The internal Hom is defined as 
\[\xymatrix{\underline{\Hom}_{\mathbf{E}-Mod^{\Sigma}}(M,N) \stackrel{defn}{=} eq(\underline{\Hom}_{\mathbf{Spect}_T^{\Sigma}(S)}(M,N) \ar@<-.5ex>[rr] \ar@<.5ex>[rr] && \underline{\Hom}_{\mathbf{Spect}_T^{\Sigma}(S)}(\mathbf{E} \wedge M, N)},\]
where the top map is $\gamma_M^* = \circ \gamma_M$ and the bottom map is 
\[\gamma_{N*}: \underline{\Hom}_{\mathbf{Spect}_T^{\Sigma}(S)}(M,N) \stackrel{\mathbf{E} \wedge}{\longrightarrow} \underline{\Hom}_{\mathbf{Spect}_T^{\Sigma}(S)}(\mathbf{E} \wedge M, \mathbf{E} \wedge N) \stackrel{\gamma_N \circ}{\longrightarrow} \underline{\Hom}_{\mathbf{Spect}_T^{\Sigma}(S)}(\mathbf{E} \wedge M,N).\]
We should also mention the following theorem of Jardine:
\begin{thm}\cite[Thm. 4.31]{Jar00}
There is a Quillen equivalence 
\[V: \mathbf{Spect}_T(S) \leftrightarrows \mathbf{Spect}_T^{\Sigma}(S) : U,\]
where $\mathbf{Spect}_T(S)$ is the category of motivic $T$-spectra, $V$ is the symmetrization functor and $U$ is the forgetful functor. 
\end{thm}
We remind the reader that throughout this work we take the motivic stable homotopy category as
\[SH(k) = Ho(\mathbf{Spect}_{T}^{\Sigma}(k)).\]
The theorem of Jardine allows us to identify $SH(k)$ equivalently to the $\A^1$-stable homotopy category $SH^{\P^1}(k) \cong SH^T(k)$ of Morel constructed in \cite[Defn. 5.1, Rem. 5.1.10 and pp. 420]{Mor04a}, which is defined as the homotopy category of the motivic $T$-spectra $\mathbf{Spect}_T(k)$. Hence, we can use Morel computation of $\pi_0^{st\A^1}(\mathbb{S}^0)_*$ and his stable $\A^1$-connectivity result.    
\bibliographystyle{plain}

\begin{thebibliography}{999999999}
\bibitem[ABGHR14]{ABGHR14} M. Ando, A. J. Blumberg, D. Gepner, M. J. Hopkins, C. Rezk, An $\infty$-categorical approach to $R$-line bundles, $R$-module Thom spectra, and twisted $R$-homology, J Topology, 7 $(\mathbf{3})$, (2014).
\bibitem[ALP15]{ALP15} A. Ananyevskiy, M. Levine, I. Panin, Witt sheaves and the $\eta$-inverted sphere spectrum, arXiv:1504.04860v1 [math.AT], Preprint (2015). 
\bibitem[AH11a]{AH11a} A. Asok, C. Haesemeyer, Stable $\A^1$-homotopy and $R$-equivalence, J. Pur. Appl. Alg., 215 $\bf{10}$, 2469-2472, (2011).
\bibitem[AH11]{AH11} A. Asok, C. Haesemeyer, The 0-th stable $\A^1$-homotopy sheaf and quadratic zero cycles, arXiv:1108.3854v1 [math.AG], Preprint (2011).
\bibitem[Ay08]{Ay08} J. Ayoub, Les six op\'eration de Grothendieck et le formalisme des cycles \'evanescents dans le monde motivique I and II. Ast\'erisque $\bf{314}$, $\bf{315}$, (2008). 
\bibitem[CD10]{CD10} C. D. Cisinski, F. D\'eglise, Triangulated category of mixed motives, arXiv:0912.2110v3 [math.AG], Preprint (2012).
\bibitem[Deg13]{Deg13} F. D\'eglise, Orientable homotopy modules, American Journal of Math., $\mathbf{135}$ (2), pp. 519-560, (2013).
\bibitem[Del87]{Del87} P. Deligne, Le d\'eterminant de la cohomologie, Current trends in arithmetic algebraic geometry, (Arata, Calif., 1985), Contemp. Math., vol. $\bf{67}$, Amer. Math. Soc., Providence, RI, 1987, pp. 93-177.
\bibitem[EGA4]{EGA4} A. Grothendieck, J. Dieudonn\'e, \'El\'ements de g\'eom\'etrie alg\'ebrique. IV. \'Etude locale des sch\'emas et des morphismes de sch\'emas IV, Publ. Math. IH\'ES $\mathbf{20, 24, 28, 32}$ (1964-1967).
\bibitem[Fas07]{Fas07} J. Fasel, The Chow-Witt ring, Doc. Math., $\bf{12}$, 275-312, (2007).
\bibitem[Fas08]{Fas08} J. Fasel, Groupes de Chow-Witt, M\'em. Soc. Math. Fr. (NS.), $\bf{113}$, (2008).
\bibitem[Hir03]{Hir03} P. S. Hirschhorn, Model categories and their localizations, Math. Surveys and Monographs, vol. $\mathbf{99}$, Amer. Math. Soc. (2003)
\bibitem[Hov99]{Hov99} M. Hovey, Model categories, Math. Surveys and Monographs, vol. $\mathbf{63}$, Amer. Math. Soc. (1999) 
\bibitem[Hov01]{Hov01} M. Hovey, Spectra and symmetric spectra in general model categories, J. Pure Appl. Alg., 165 $\bold{(1)}$, 63-127, (2001). 
\bibitem[Hu05]{Hu05} P. Hu, On the Picard group of the stable $\A^1$-homotopy category, Topology, 44 $\bf{(3)}$, 609-640, (2005).
\bibitem[Jar00]{Jar00} J.F. Jardine, Motivic symmetric spectra, Doc. Math. $\bf{5}$, (2000), 445-553
\bibitem[Lev10]{Lev10} M. Levine, Slices and transfers, Doc. Math., 393-443, (2010).  
\bibitem[Lev13]{Lev13} M. Levine, Convergence of Voevodsky's slice tower, Doc. Math., $\mathbf{18}$, 907-941, (2013).
\bibitem[MVW06]{MVW06} C. Mazza, V. Voevodsky, C. Weibel, Lecture Notes on motivic cohomology, vol. 2, Clay Math. Mono., Amer. Math. Soc., Providence RI, (2006).  
\bibitem[Mor04a]{Mor04a} F. Morel, An introduction to $\A^1$-homotopy theory, Contemporary developments in algebraic $K$-theory, 357-441, ICTP Lect. Notes, $\bf{XV}$, Abdus Salam Int. Cent. Theoret. Phys., Trieste, (2004).
\bibitem[Mor04]{Mor04} F. Morel, On the motivic stable $\pi_0$ of the sphere spectrum, in Axiomatic, Enriched and Motivic Homotopy Theory, 219-260, J.P.C. Greenlees (ed.), (2004), Kluwer Academic Publishers.
\bibitem[Mor05]{Mor05} F. Morel, The stable $\A^1$-connectivity theorems, $K$-theory, $\bf{35}$, 1-68, (2005).
\bibitem[Mor12]{Mor12} F. Morel, $\A^1$-algebraic topology over a field, vol. 2052, Lect. Notes in Math., Springer Heidelberg, (2012).
\bibitem[MV01]{MV01} F. Morel, V. Voevodsky, $\A^1$-homotopy theory of schemes, IH\'ES Publ. Math., $\bf{(90)}$, 45-143, (1999).
\bibitem[NSO09]{NSO09} N. Naumann, M. Spitzweck, Paul Arne \O stv\ae r, Chern classes, $K$-theory and Landweber exactness over nonregular base schemes, in Motives and Algebraic Cycles: A Celebration in Honour of Spencer J. Bloch, Fields Inst. Comm., Vol. $\mathbf{56}$, (2009). 
\bibitem[Q67]{Q67} D. G. Quillen, Homotopical algebra, Lect. Notes in Math., No. $\mathbf{43}$, Springer, (1967).
\bibitem[Rio05]{Rio05} J. Riou, Dualit\'e de Spanier-Whitehead en g\'eom\'etrie alg\'ebrique, C. R. Math. Acad. Sci. Paris $\mathbf{340}$, no. $\mathbf{6}$, (2005).
\bibitem[Rio10]{Rio10} J. Riou, Algebraic $K$-theory, $\A^1$-homotopy and Riemann-Roch theorems, J. of Topology $\bf{3}$, (2010), pp. 229-264.
\bibitem[RO08]{RO08} Oliver R\"ondigs, Paul Arne \O stv\ae r, Modules over motivic cohomology, Adv. in Math. vol. $\mathbf{219}$, Issue $\mathbf{2}$, (2008). 
\bibitem[SS00]{SS00} S. Schwede, B. E. Shipley, Algebras and modules in monoidal model categories, Proc. London Math. Soc. $\mathbf{80}$, pp. 491-511, (2000).
\bibitem[Voe00]{Voe00} V. Voevodsky, Triangulated categories of motives over a field, in Cycles, transfers and motivic homology theories, Annals of Mathematics Studies, vol. 143, Princeton Univ. Press, (2000).
\bibitem[Voe03]{Voe03} V. Voevodsky, Motivic cohomology with $\Z/2$-coefficients, Publ. Math. IH\'ES (2003).
\bibitem[Weib94]{Weib94} C. A. Weibel, An introduction to homological algebra, Cambridge stud. in adv. Math. $\mathbf{38}$, (1994).
\end{thebibliography}
\renewcommand\refname{References}

\end{document}